\documentclass[12pt]{amsart}
\advance\textheight\topskip
\usepackage[OT2,T1]{fontenc}
\newcommand\cyr{%
  \renewcommand\rmdefault{wncyr}%
  \renewcommand\sfdefault{wncyss}%
  \renewcommand\encodingdefault{OT2}%
  \normalfont
  \selectfont}
\DeclareTextFontCommand{\textcyr}{\cyr}

\usepackage{times}
\usepackage[leqno]{amsmath}
\usepackage{amsfonts,amssymb,amsthm}%%
\allowdisplaybreaks
\usepackage[mathscr]{eucal}

\renewcommand{\geq}{\geqslant}
\renewcommand{\leq}{\leqslant}
\newcommand{\charn}[2]{\omega^{#1}_{#2}}
\newcommand{\mD}{\mathfrak{D}}

\newcommand{\THETA}{\boldsymbol{\mathit{\Theta}}}
\newcommand{\SLiiZ}{SL(2,\oZ)}
\newcommand{\vermaNplus}[1]{\module{N}_{#1}^+}
\newcommand{\vermaNminus}[1]{\module{N}_{#1}^-}

\newcommand{\mfrac}[2]{\mbox{\small$\displaystyle\frac{#1}{#2}$}}
\newcommand{\ffrac}[2]{\raisebox{.6pt}{\mbox{%
    \standardfootnotesize$\displaystyle\frac{#1\mathstrut}{#2\mathstrut}$}}}

\newcommand{\fbinom}[2]{\mbox{%
    \standardfootnotesize$\displaystyle\binom{#1}{#2}$}}

\def\Xint#1{\mathchoice
   {\XXint\displaystyle\textstyle{#1}}%
   {\XXint\textstyle\scriptstyle{#1}}%
   {\XXint\scriptstyle\scriptscriptstyle{#1}}%
   {\XXint\scriptscriptstyle
                      \scriptscriptstyle{#1}}%
   \!\int}
\def\XXint#1#2#3{{%
     \setbox0=\hbox{$#1{#2#3}{\int}$}
     \vcenter{\hbox{$#2#3$}}\kern-.5\wd0}}
\def\dashint{\Xint-}

\usepackage[curve,matrix,arrow,cmtip]{xy}
\usepackage{epic,eepic}
\usepackage{hyperref}

\numberwithin{equation}{section}
\makeatletter
\@addtoreset{equation}{section}

\def\@secnumfont{\bfseries}
\def\subsubsection{\@startsection{subsubsection}{3}%
  \z@{.5\linespacing\@plus.7\linespacing}{-.5em}%
  {\normalfont\bfseries}}
\def\paragraph{\@startsection{paragraph}{4}%
  \z@\z@{-\fontdimen2\font}%
  \normalfont\bfseries}
\def\subparagraph{\@startsection{subparagraph}{5}%
  \z@\z@{-\fontdimen2\font}%
  \normalfont\bfseries}

\renewenvironment{smallmatrix}{\null\,\vcenter\bgroup
  \Let@\restore@math@cr\default@tag
  \baselineskip9\ex@ \lineskip1.5\ex@ \lineskiplimit\lineskip
  \ialign\bgroup\kern-1pt\hfil$\m@th\scriptstyle##$\hfil&&\thickspace\thinspace\hfil
  $\m@th\scriptstyle##$\hfil\crcr
}{%
  \crcr\egroup\egroup\kern1pt%
}

\def\standardfootnotesize{\footnotesize}

\makeatother

\newcommand{\smatrix}[4]{\left(\begin{smallmatrix}#1&#2\\
      #3&#4\end{smallmatrix}\right)}

\newcommand{\pell}{\ell}
\newcommand{\rbar}[2]{{\boldsymbol{r}}^{\phantom{y}}_{\kern-1pt#1}[#2]}

\newcommand{\ldot}{\mathbin{\boldsymbol{.}}}

\newcommand{\chzero}{\textrm{\textup{R}}}
\newcommand{\schzero}{\textrm{\textup{sR}}}
\newcommand{\chzerotext}{Ramond}

\newcommand{\chhalf}{\textrm{\textup{NS}}}
\newcommand{\schhalf}{\textrm{\textup{sNS}}}
\newcommand{\chhalftext}{Neveu--Schwarz}

\newcommand{\acycle}{\boldsymbol{a}}
\newcommand{\bcycle}{\boldsymbol{b}}
\newcommand{\Jacobi}[2]{\genfrac{(}{)}{}{1}{#1}{#2}}
\newcommand{\acts}{\mathbin{\boldsymbol{.}}}
\newcommand{\upperH}{\mathfrak{h}} %{\mathbb{H}}
\newcommand{\mat}[1]{\mathsf{#1}}

\newcommand{\floor}[1]{\left\lfloor#1\right\rfloor}
\newcommand{\bundle}[1]{{\mathbb{V}_{#1}}}
\newcommand{\Erf}{\mathop{\mathrm{erf}}}
\newcommand{\Erfc}{\mathop{\mathrm{erfc}}}
\newcommand{\torus}[1]{\mathbb{T}_{#1}}
\newcommand{\lattice}[1]{\mathbb{L}_{#1}}

\newcommand{\eh}{e_{\phantom{h}}}

\newcommand{\EX}[1]{\textit{\normalsize
    e}^{\kern1pt\mbox{\standardfootnotesize$#1$}}_{\vphantom{x}}}
\newcommand{\ex}[1]{\textit{\normalsize
    e}^{\kern1pt\mbox{\standardfootnotesize$#1$}}}

\newcommand{\one}{\boldsymbol{1}}%{1\kern-4pt 1}

\newcommand{\K}{{\kern1pt\mathscr{K}}}

\newcommand{\myvect}{\mathsf} %{\boldsymbol}

\newcommand{\vect}[1]{\boldsymbol{#1}}
\newcommand{\Kvect}{\myvect{K}}

\newcommand{\Tr}{\mathop{\mathrm{Tr}}\nolimits}

\newcommand{\vt}[3]{\vartheta_{#1}(#2,#3)}

\newcommand{\bref}[1]{\textbf{\ref{#1}}}

\providecommand{\half}{\frac{1}{2}}

\newcommand{\fhalf}{\ffrac{1}{2}}

\newcommand\res{\mathop{\mathrm{res}}\limits}

\newtheorem{Thm}{Theorem}[section]
\newtheorem{Lemma}[Thm]{Lemma}

\theoremstyle{definition}

\newenvironment{prf}{%
  \noindent{\sc Proof.}}%
{\noindent{\nolinebreak\hfill\mbox{\rule{.5em}{.5em}}\,}\pagebreak[3]%
  \par\medskip}

\newcommand{\HyperF}[2]{\,{}^{\phantom{z}}_{#1}\kern-1pt F_{#2}^{\phantom{z}}}
\newcommand{\tSSL}[2]{\widehat{s\ell}(#1|#2)}

\newcommand{\hSSL}[2]{\widehat{s\ell}(#1|#2)}

\newcommand{\hSL}[1]{\widehat{s\ell}(#1)}
\newcommand{\SL}[1]{s\ell(#1)}

\newcommand{\hD}{\widehat{D}(2|1;\alpha)}

\newcommand{\N}[1]{N\!=\!#1}

\newcommand{\module}{\mathcal}
\def\mL{\module{L}}    %<----- \hSSL21 admissible representations
\def\mM{\module{M}}
\newcommand{\mP}{\module{P}}

\newcommand{\verma}[1]{\module{P}_{#1}}
\newcommand{\ket}[1]{|{#1}\rangle}
\newcommand{\ketplus}[1]{|{#1}\rangle^+}
\newcommand{\ketminus}[1]{|{#1}\rangle^-}

\newcommand{\Hminus}{H^-}
\newcommand{\Hplus}{H^+}

\newcommand{\oC}{\mathbb{C}}
\newcommand{\oN}{\mathbb{N}}

\newcommand{\oR}{\mathbb{R}}
\newcommand{\oZ}{\mathbb{Z}}

\def\cA{\mathcal{A}}

\def\cT{\mathcal{T}}
\def\cU{\mathcal{U}}

\def\bar{\overline}

\newcommand{\dd}{\partial}

\begin{document}

%% \runauthor{AM\,Semikhatov, IYu\,Tipunin and A\,Taormina}

%% \begin{frontmatter}

\title[Appell Functions and Characters]{%
  \vspace*{-2\baselineskip}{\mbox{}\hfill\texttt{\normalsize% 
        \lowercase{math}.QA/0311314}}\\[2\baselineskip]
    Higher-Level Appell Functions, Modular Transformations, and
    Characters}

\author[AM\,Semikhatov]{A.\kern2ptM.~Semikhatov} 

\address{AMS, IYuT: Lebedev Physics Institute, Russian Academy of
  Sciences}

\author[A\,Taormina]{A.~Taormina} 

\address{AT: Department of Mathematical Sciences, University of
  Durham}

\author[IYu\,Tipunin]{I.\kern2ptYu.~Tipunin}

\dedicatory{\textcyr{\itshape Bore Fei0ginu po sluchayu ego
    pyatidesyatiletiya}\rule{0pt}{16pt}}

\begin{abstract}
  We study modular transformation properties of a class of indefinite
  theta series involved in characters of infinite-dimensional Lie
  superalgebras.  The \textit{level-$\pell$ Appell
    functions}~$\K_{\pell}$ satisfy open quasiperiodicity relations
  with additive theta-function terms emerging in translating by the
  ``period.''  \ Generalizing the well-known interpretation of theta
  functions as sections of line bundles, the $\K_\pell$ function
  enters the construction of a section of a rank-$(\pell+1)$ bundle
  $\bundle{\pell,\tau}$.  We evaluate modular transformations of
  the~$\K_{\pell}$ functions and construct the action of an
  $SL(2,\oZ)$ subgroup that leaves the section of
  $\bundle{\pell,\tau}$ constructed from~$\K_{\pell}$ invariant.

  \smallskip
  
  Modular transformation properties of~$\K_{\pell}$ are applied to the
  affine Lie superalgebra $\hSSL21$ at rational level $k>-1$ and to
  the $\N2$ super-Virasoro algebra, to derive modular transformations
  of ``admissible'' characters, which are not periodic under the
  spectral flow and cannot therefore be rationally expressed through
  theta functions.  This gives an example where constructing a modular
  group action involves extensions among representations in a
  nonrational conformal model.
\end{abstract}


\maketitle
%% \end{frontmatter}
%%%%%%%%%%%%%%%%%%%%%%%%%%%%%%%%%%%%%%%%%%%%%%%%%%%%%%%%%%%%%%%%%%%%%
\thispagestyle{empty}

\section{Introduction}
%%%
We generalize some elements of the theta-function theory by studying
modular transformations of functions that are not doubly quasiperiodic
in a variable \hbox{$\mu\in\oC$}.  Such functions emerge in the study
of characters of representations in (nonrational) conformal field
theory models based on Lie superalgebras, which motivates
investigation of their modular properties.\enlargethispage{12pt}

A modular group representation associated with characters of a
suitable set of representations is a fundamental property of conformal
field theory models, related to the fusion algebra by the Verlinde
formula, via the argument traced to the consistency of gluing a
three-punctured sphere into a one-punctured torus\,---\,in fact, to
consistency of conformal field theory itself~\cite{[MS],[V]}.
Strictly speaking, this applies to rational conformal field theories,
where modular properties of the characters and the structure of the
Verlinde formula are known, at least in principle (for a discussion of
the modular transformation properties of characters and other
quantities and for further references,
see~\cite{[G],[By],[BK],[FS],[H]}).  Modular behavior of theta
functions can be considered a basic feature underlying good modular
properties in rational models (in particular, the well-known modular
group representation on a class of characters of affine Lie
algebras~\cite{[K-book]}); it is deeply connected with
quasiperiodicity of theta functions and hence of the characters in
rational models.

But characters that are not quasiperiodic (are not invariant under
lattice translations, often called ``spectral flows'' in that context)
cannot be rationally expressed through theta functions.  Such
characters often occur in nonrational conformal field theory models
(an infinite orbit of the spectral flow transform already implies that
the theory is nonrational).  Modular properties of such characters
present a problem both technically (the theta-function theory is of
little help) and conceptually (it is unclear what kind of modular
invariance is to be expected at all).

On the other hand, the paradigm that any consistent conformal field
theory must be related to a modular group representation, even beyond
the class of rational theories, motivates studying modular behavior of
nonrational characters and, on the technical side, seeking an adequate
``replacement'' of theta functions with some functions that are not
quasiperiodic but nevertheless behave reasonably under modular
transformations and can be used as ``building blocks'' of the
characters.  Such functions are to be found among indefinite theta
series (see~\cite{[GZ],[KP],[KW],[Pol-00]} and references therein).

In this paper, we study the modular (and other related) properties of
\textit{higher-level Appell functions}\,---\,a particular instance of
indefinite theta series, not-double-quasiperi\-odic functions involved
in the characters of modules of the $\hSSL21$ affine Lie
superalgebra~\cite{[ST]} and the $\N2$ and $\N4$ superextensions of
the Virasoro algebra (\cite{[FSST],[ET88]}).  Remarkably, the pattern
of modular behavior established for the Appell functions is then
reproduced by the characters.  For a positive integer $\pell$, we
define the level-$\pell$ Appell~function~as\enlargethispage{12pt}
\begin{gather}\label{Knew-exp}
  \K_\pell(\tau,\nu,\mu)=
  \sum_{m\in\oZ}
  \ffrac{
    \EX{i\pi m^2 \pell\tau\!+\!2i\pi m \pell\nu}}{
    1 - \EX{2i\pi(\nu\!+\!\mu\!+\!m\tau)}},
  \quad
  \begin{array}[t]{ll}
    \tau\in\upperH,~&\nu,\mu\in\oC,\\
    {}&\mu+\nu\notin\oZ\tau+\oZ.
  \end{array}
\end{gather}

\begin{Thm}\label{thm:TS}
  The level-$\pell$ Appell function $\K_\pell$ satisfies the relations
  \begin{gather}\label{Knew-Ttransf}
    \K_{\pell}(\tau+1,\nu,\mu)
    = \begin{cases}
      \K_{\pell}(\tau,\nu+\half,\mu-\half),&\pell\quad\text{odd},\\
      \K_{\pell}(\tau,\nu,\mu),&\pell\quad\text{even},
    \end{cases}
  \end{gather}
  and
  \begin{multline}\label{Knew-Stransf-fin}
    \K_\pell(-\ffrac{1}{\tau},\ffrac{\nu}{\tau},\ffrac{\mu}{\tau}) {}=
    \tau \smash{\EX{i\pi\pell\frac{\nu^2 - \mu^2}{\tau}}}\,
    \K_\pell(\tau,\nu,\mu)\\*
    \smash[b]{{}+
      \tau \sum_{a=0}^{\pell-1}
      \EX{i\pi\frac{\pell}{\tau}(\nu\!+\!\frac{a}{\pell}\tau)^2}}
    \Phi(\pell\tau,\pell\mu\!-\!a\tau)
    \vartheta(\pell\tau, \pell\nu\!+\!a \tau),
  \end{multline}
  where
  \begin{gather}\label{Phi-def}
    \Phi(\tau,\mu)=
    -\ffrac{i}{2\sqrt{-i\tau}} -\fhalf\int_{\oR}\!dx\,\EX{-\pi x^2}\,
    \smash[t]{
      \frac{\sinh\!\left(\pi
          x\sqrt{-i\tau}(1\!+\!2\frac{\mu}{\tau})\right)}{
        \sinh\!\left(\pi x\sqrt{-i\tau}\right)}}
  \end{gather}
\end{Thm}
\noindent
(we refer to~\eqref{theta10}--\eqref{thetaKac} for the theta-function
notation).\pagebreak[3]

The simplest, level-$1$ Appell function $\kappa(x,y; q) =
\K_1(q,y^{-1}q,xyq^{-1})$ has appeared
in~\cite{[Polisch],[KW-Appell]}, and its $S$-transformation properties
were formulated in~\cite{[Polisch]} as the statement that the
difference between~$\kappa$ and its $S$-transform is divisible by the
theta function~$\vartheta(\tau, \mu)$.  Theorem~\bref{thm:TS}
generalizes this to $\pell>1$ and in addition gives an integral
representation of the function $\Phi$ accompanying the
theta-functional terms in the modular transform.  This integral
representation allows studying the~$\Phi$~function, which is an
important ingredient of the theory of higher-level Appell functions,
similarly to Barnes-related functions arising
elsewhere~\cite{[B],[FK],[PT],[FKV],[JM],[StsKL]}.  We derive
functional equations satisfied by~$\Phi$ and the formula for its
$S$-transformation.

Already with theta functions, their modular properties are closely
related to (and can in fact be derived from) their quasiperiodicity
under lattice translations, which in geometric terms means that the
theta function $\vartheta(\tau,{}\cdot{})$ represents a section of a
line bundle over the torus determined by the modular parameter~$\tau$
(\textit{hence} the dependence on the second argument $\nu\in\oC$ is
doubly quasiperiodic).  With the $\K_\pell$ functions, which are no
longer doubly quasiperiodic, the geometric counterpart of ``open
quasiperiodicity'' (a[dditive]-quasiperiodicity in~\cite{[S]})
involves sections of rank-$(\pell+1)$ bundles.  The simplest Appell
function $\kappa$ in~\cite{[Polisch],[KW-Appell]} satisfies an open
quasiperiodicity relation with an additive \textit{theta-func\-tional}
term arising in shifting the argument by the ``period,''
\begin{equation}\label{kappa-open}
  \kappa(z q, a; q) = a\kappa(z, a; q) + \vartheta(q,z),
\end{equation}
which implies that $(\kappa({}\cdot{},a; q), \vartheta(q,{}\cdot{}))$
represents a section of a rank-$2$ bundle over the elliptic
curve~\cite{[Polisch]}.  Analogously, the higher-level Appell
functions, in the ``multiplicative'' notation,\footnote{We resort to
  the standard abuse of notation, by freely replacing functional
  arguments by their exponentials or conversely, by logarithms, via
  $q=e^{2i\pi\tau}$, $x=e^{2i\pi\nu}$, $y=e^{2i\pi\mu}$,~etc.} are
quasiperiodic under~$x\mapsto xq$ and satisfy an
%% \textit{inhomogeneous} finite-difference equation
open quasiperiodicity relation with the inhomogeneous terms involving
theta functions,
\begin{equation}\label{Knew-open}
  \K_\pell(q, x, y\,q)=q^{\frac{\pell}{2}}\,y^\pell
  \K_\pell(q, x, y)
  +\sum_{a=0}^{\pell-1} x^a\,y^a\,q^a\,\vartheta(q^\pell, x^\pell\,q^a).
\end{equation}
This generalizes~\eqref{kappa-open}, to which~\eqref{Knew-open}
reduces for $\pell=1$ (but there also exists a ``finer'' property for
$\pell>1$, see Sec.~\ref{sec:Appell-elementary}).
% \footnote{Incidentally, a certain version of the $\K_0$ function,
%   namely, the function $P_\lambda(x,q)=
%   \K_\pell(\tau,x^{\frac{1}{\pell}},\lambda\,x^{-1})
%   \bigr|_{\pell=0}$, related to an Eisenstein series, was introduced
%   in~\cite{[D]} for $\lambda$ a root of unity.}
The theta functions occurring in the right-hand side of the modular
transform~\eqref{Knew-Stransf-fin} are precisely those violating the
quasiperiodicity of~$\K_\pell$ in~\eqref{Knew-open}.
Together with $\K_\ell$, these theta functions enter the
construction of a section of a rank-$(\pell+1)$ bundle,
\begin{equation*}
  (\K_\pell(\tau,\nu,\mu),\;
  C_1\vartheta(\pell\tau,\pell\nu),\;\dots,\;
  C_{\pell}\vartheta(\pell\tau,\pell\nu +
  (\pell-1)\tau)),
\end{equation*}
with $C_a$ such that the entire vector is invariant under a subgroup
of lattice translations (see Lemma~\bref{Lemma:Kvect} in what
follows).  Moreover, it turns out that the action of a~subgroup of the
modular group can be defined on $(\pell+1)$-vectors
$F=(f_0(\tau,\nu,\mu), f_1(\tau,\nu),\dots, f_{\pell}(\tau,\nu))$ such
that the above bundle section is an invariant of this action
(Theorem~\bref{thm:abcd}, which gives a more ``invariant'' formulation
of the modular properties of the Appell functions).  This gives an
interesting realization of modular invariance with a \textit{matrix}
automorphy factor (cf.~\cite{[EhSk],[FHST]}).

As noted above, the Appell functions are a specific example of
indefinite theta series motivated by the study of characters.  We use
them to express the characters of ``admissible'' $\hSSL21$-modules at
rational level $k>-1$ and to study modular transformation properties
of the characters.  In this $\hSSL21$ example, indeed, the
higher-level Appell functions prove an adequate substitute for theta
functions; their modular behavior described in Theorem~\bref{thm:TS}
is essentially ``inherited'' by the admissible characters, whose
$S$-transform is given by
\begin{gather}\label{SandR}
 \chi_A( -\ffrac{1}{\tau},\ffrac{\nu}{\tau},
  \ffrac{\mu}{\tau})
  =\EX{i\pi k\frac{\nu^2 - \mu^2}{2\tau}}\sum_B S_{AB}
  \chi_B(\tau,\nu,\mu)+
  \sum_\alpha R_{A\alpha}(\tau,\nu,\mu)\,
  \Omega_\alpha(\tau,\nu,\mu),
\end{gather}
where $S_{AB}$ is a numerical matrix and the functions $R_{A\alpha}$
are expressed in terms of the above function $\Phi$, and
$\Omega_\alpha$ are some characters expressed through theta functions,
see Theorem~\bref{thm:S-chi} for the precise statement.  This shows a
triangular structure of the same type as in~\eqref{Knew-Stransf-fin}:
the additional elements $\Omega_\alpha$ occurring in the
unconventional $S$-transform formula of the~$\chi_A$ (which are not
quasiperiodic and hence cannot be rationally expressed through theta
functions) are expressible in terms of theta functions and are
therefore quasiperiodic under lattice shifts and carry a modular group
action.

The theta-functional terms $\Omega_\alpha$ turn out to be the
characters of certain \textit{extensions} among the admissible
$\hSSL21$-representations.  A key feature underlying most of the
unconventional properties of a number of nonrational conformal field
theories is that the irreducible representations allow nontrivial
extensions among themselves (by which we mean nonsplittable short
exact sequences, or actually the middle modules in such sequences).
Such extensions do not occur in rational conformal field theory
models.

Behavior of the admissible $\hSSL21$-characters under modular
transformations is related to their behavior under \textit{spectral
  flow transformations}, i.e., a representation of a lattice~$\Lambda$
whose elements $\vec\theta$ act via adjoint representation of the
group elements~$\exp(\vec\theta\cdot\vec h)$, where~$\vec h$ are
Cartan subalgebra elements.\footnote{The term spectral flow transform
  is taken over from the $N{=}2$ super-Virasoro algebra~\cite{[SS]}.}
The role of spectral flows appears to originate in the fact that the
fundamental group of the appropriate moduli space is not just
$\SLiiZ$, but its semidirect product with a lattice.  (Lattice
translations also require considering the so-called \chzerotext{} and
\chhalftext{} ``sectors'' and supercharacters.)  The admissible
$\hSSL21$-characters at the level $k\!=\!\frac{\pell}{u}-1$ acquire
additional theta-functional terms under the spectral flow transform
$\cU_\theta$ with $\theta\!=\!u$,
\begin{gather*}
  (\cU_u\chi_{A})(\tau,\nu,\mu)=\chi_A(\tau,\nu,\mu)
  +\sum_\alpha f_\alpha(\tau,\nu,\mu)\,\Omega_\alpha(\tau,\nu,\mu),
\end{gather*}
where $f_\alpha(\tau,\nu,\mu)$ are some trigonometric functions and
$\Omega_\alpha$ are the same as in the modular transform
formula~\eqref{SandR}.  This demonstrates an obvious similarity to the
properties of the Appell functions (the same theta functions occurring
in~\eqref{Knew-Stransf-fin} and~\eqref{kappa-open}).

\medskip

\noindent\textbf{This paper is organized as follows.} \
In Sec.~\ref{sec:APPELL}, we study the level-$\pell$ Appell functions
$\K_\pell$.  The basic quasiperiodicity and some other properties
are derived in Sec.~\ref{sec:Appell-elementary}.  The $\K_\pell$
functions satisfy several ``period multiplication'' formulas, which we
derive in Sec.~\ref{sec:scaling}.  Formulas relating the higher-level
Appell functions to theta functions are given in
Sec.~\ref{sec:to-theta}.  In Sec.~\ref{sec:integral-rep}, we then
derive modular transformation properties of the higher-level Appell
functions using their \textit{integral} representation through theta
functions.  The Appell functions of the lower levels are briefly
considered in Sec.~\ref{sec:lowest}.  In Sec.~\ref{sec:Phi}, we next
consider the $\Phi$ function arising in the modular transformation of
$\K_\pell$; its properties are in some respects analogous to the
properties of $\K_\pell$.  The geometric point of view on the
higher-level Appell functions is outlined in Sec.~\ref{sec:geometry}.
There, we first (in Sec.~\ref{sec:lattice}) consider how the
$\K_\pell$ function and the appropriate theta functions are
combined to produce a section of a rank-$(\pell+1)$ bundle.  In
Sec.~\ref{sec:modular-full}, we then derive the action of a subgroup
of $\SLiiZ$ on these sections (Theorem~\bref{thm:abcd}).

In Sec.~\ref{sec:sl21}, we use the established properties of the
higher-level Appell functions to evaluate modular transformation
properties of the class of ``admissible'' $\hSSL21$-characters.  The
main result (Theorem~\bref{thm:S-chi}) is formulated in
Sec.~\ref{sec:notaion}.  The characters are expressed through the
higher-level Appell functions in Sec.~\ref{subsec:q-periodicity}.
Their $S$-transformation formula is derived in Sec.~\ref{sec:sl21-S}.
Application of the higher-level Appell functions to $\N2$ and $\N4$
super-Virasoro theories is outlined in Sec.~\ref{sec:N24}.  

In Appendix~\ref{app:ab-int}, we evaluate several useful contour
integrals over the torus involving theta and Appell functions.  In
Appendix~\ref{app:sl21}, we recall the $\hSSL21$ affine Lie
superalgebra, consider its automorphisms (Sec.~\ref{app:the-algebra})
and define some of its modules (Sec.~\ref{app:sl21-modules}), and
finally give the admissible representation characters
(Theorem~\bref{thm:find-chars}, Sec.~\ref{app:sl21-adm}).  The
different ``sectors'' and the corresponding characters are given in
Sec.~\ref{sec:sectors}.

\medskip

\noindent
\textbf{Notation.} \ We let $\upperH$ denote the upper complex
half-plane.  The group $SL(2,\oZ)$ is generated by the two matrices
\begin{gather*}
  S= \mbox{\small$\begin{pmatrix}
    0~&-1\\
    1~&0
  \end{pmatrix}$}\!,\qquad
  T= \mbox{\small$\begin{pmatrix}
    1~ & 1\\
    0~ & 1
  \end{pmatrix}$}
\end{gather*}
with the relations
\begin{gather*}%%\label{eq:SL2relations}
  S^2=(ST)^3=(TS)^3=C,
\end{gather*}
where $C^2=\one$.  The standard $SL(2,\oZ)$ action
on~$\upperH\times\oC^2$ is
\begin{gather}\label{SL2-action}
  \gamma=
  \mbox{\small$\begin{pmatrix}
    a&b\\
    c&d 
  \end{pmatrix}$}
  {}:(\tau,\nu,\mu)\mapsto
  (\gamma\tau,\gamma\nu,\gamma\mu)=
  \bigl(
  \ffrac{a\tau+b}{c\tau+d},
  \ffrac{\nu}{c\tau+d},\ffrac{\mu}{c\tau+d}
  \bigr)
\end{gather}
(where the notation $\gamma\nu$ and $\gamma\mu$ is somewhat loose,
because this action depends on~$\tau$).

We use the classical theta functions%%\enlargethispage{24pt}
\begin{gather}
  \vt{1,0}{q}{z}= \sum_{m\in\oZ}^{}q^{\half(m^2 - m)} z^{-m}
  \label{theta10}
  = \prod_{m\geq0}\!(1\!+\!z^{-1}q^m)\prod_{m\geq1}\!(1\!+\!z q^m)
  \prod_{m\geq1}\!(1\!-\!q^m),\\
  \label{theta11}
  \vartheta_{1,1}(q, z)
  =  \sum_{m\in\oZ}q^{\half(m^2 - m)} (-z)^{-m}
  = \prod_{m\geq0}\!(1\!-\!z^{-1} q^m)
  \prod_{m\geq1}\!(1\!-\!z q^m)\prod_{m\geq1}\!(1\!-\!q^m),\\
  \label{plain-theta}
  \vartheta(q,z)\equiv\vartheta_{0,0}(q,z)
  =\sum_{m\in\oZ}q^{\frac{m^2}{2}}z^m
  =\vartheta_{1,0}(q,z\,q^{-\half}).
\end{gather}
Their $S$ transformations are given by
\begin{gather}\label{modular-theta}
  \begin{split}
    \vartheta_{1,1}(-\ffrac{1}{\tau},\ffrac{\nu}{\tau})&=
    -i\sqrt{-i\tau}\,
    \smash[t]{\EX{i\pi\nu\!+\!i\pi\frac{1}{\tau}(\nu\!-\!\half)^2
        \!+\!\frac{i\pi}{4}\tau}}\,
    \vartheta_{1,1}(\tau,\nu),\\ 
    \vartheta_{1,0}(-\ffrac{1}{\tau},\ffrac{\nu}{\tau})&=
    \sqrt{-i\tau}\,
    \EX{i\pi\frac{(\nu-\half)^2} \tau}
    \vartheta_{1,0}(\tau, \nu + \ffrac{1}{2} - \ffrac{\tau}{2}),
  \end{split}
\end{gather}
and
\begin{gather}
  \vartheta\smash{(-\ffrac{1}{\tau},\ffrac{\nu}{\tau})=
    \sqrt{-i\tau}\,
    \smash[t]{\EX{i\pi\frac{\nu^2}{\tau}}}
    \vartheta(\tau,\nu)}.
\end{gather}
The eta function
\begin{gather}\label{eta}
  \eta(q)  
  =q^{\frac{1}{24}}
  \smash{\prod\limits_{m=1}^\infty(1-q^{m})}
\end{gather}
transforms as 
\begin{gather}\label{modular-eta}
  \eta(\tau+1)=\smash{\EX{\frac{i\pi}{12}}\eta(\tau),\qquad
    \eta(-\ffrac{1}{\tau})=
    \sqrt{-i\tau}\,\eta(\tau)}.
\end{gather}

In a different theta-functional nomenclature, one introduces the
higher-level theta functions
\begin{gather}\label{thetaKac}
  \theta_{r,\pell}(q,z)=\sum_{j\in\oZ+\frac{r}{2\pell}}
  q^{\pell j^2}z^{\pell j}
  =z^{\frac{r}{2}}\,q^{\frac{r^2}{4\pell}}\,
  \vartheta(q^{2\pell},z^\pell q^r).
\end{gather}
Either $\theta$ or $\vartheta$ turn out to be more convenient
depending on circumstances.

For a positive integer $p$, we use $[x]_p$ to denote \ 
$x\,\mathrm{mod}\,p=x-p\floor{\frac{x}{p}}$, where $\floor{x}$ is the
greatest integer less than or equal to~$x$.

\section{Higher-Level Appell Functions}\label{sec:APPELL}

\subsection{Open quasiperiodicity and other basic properties}
\label{sec:Appell-elementary} For $\pell\in\oN$, the level-$\pell$
Appell function 
\begin{gather}\label{Knew}
  \K_\pell(q,x,y)=
  \sum_{m\in\oZ}\frac{\displaystyle
    q^{\frac{m^2 \pell}{2}} x^{m \pell}}{\displaystyle
    1 - x\,y\, q^{m}}, %\qquad~\Im\tau>0.
\end{gather}
generalizes the Appell function $\kappa$
in~\cite{[Polisch],[KW-Appell]}.  Along with~\eqref{Knew}, there is
another, double-series representation for $\K_\pell$
(cf.~\cite{[KW],[D],[Pol-00]})
\begin{gather}\label{K-mn}
  \K_\pell(q,x,y) =
  \Bigl(\sum_{m\geq0}\sum_{n\geq0}-\sum_{m\leq-1}\sum_{n\leq-1}
  \Bigr)
  q^{\frac{m^2 \pell}{2} + m n} x^{m \pell + n} y^{n},
\end{gather}
valid for $|q|<|x y|<1$.

The $\K_{\pell}$ functions satisfy an easily derived
quasiperiodicity property in the second argument,
\begin{gather} \label{quasi-periodicity}
  \K_\pell(q,xq^n,y) =
  q^{-\frac{n^2 \pell}{2}}\,x^{-n \pell}\K_\pell(q,x,y),\quad
  n\in\oZ,
\end{gather}
and an ``open quasiperiodicity'' relation along the antidiagonal with
respect to the second and third arguments,
\begin{equation}\label{diag-periodicity}
  \K_\pell(q,xq^{-\frac{n}{\pell}},yq^{\frac{n}{\pell}})
  = (xy)^{n}\K_\pell(q,x, y) +
  \begin{cases}
    \displaystyle
    \sum_{r=1}^{n}(xy)^{n-r}\vartheta(q^\pell,x^\pell q^{-r}),
    &n\in\oN,\\
    \displaystyle
    -\!\sum_{r=n+1}^{0}(xy)^{n-r}\!\vartheta(q^\pell,x^\pell q^{-r}),
    \kern-4pt
    &n\in-\oN.
  \end{cases}\kern-6pt
\end{equation}
These imply open quasiperiodicity in the third argument,
\begin{gather*}
  \K_\pell(q,x,yq^{n}) 
  = q^{\frac{n^2 \pell}{2}}y^{n\pell}\K_\pell(q,x,y)
    + \begin{cases}\displaystyle      
      \sum_{j=0}^{\pell n-1}
      x^j y^{j} q^{nj}\vartheta(q^\pell,x^\pell q^{j}), & n\in\oN,\\
      \displaystyle
      - \sum_{j=\pell n}^{-1}
      x^j\,y^{j}q^{nj}\vartheta(q^\pell,x^\pell q^{j}), & n\in-\oN,
    \end{cases}
\end{gather*}
or manifestly with only $\pell$ distinct theta functions in the
right-hand side,
\begin{multline}\label{veluti-periodicity}
  \K_\pell(q,x,yq^{n})
  = q^{\frac{n^2 \pell}{2}}y^{n\pell}\K_\pell(q,x,y)\\*
  {}+{}
  \begin{cases}\displaystyle
    \sum_{j=0}^{n-1}q^{\frac{j(2n-j)\pell}{2}}y^{j\pell}
    \sum_{r=0}^{\pell-1}x^r y^{r} q^{(n-j)r}
    \vartheta(q^\pell,x^\pell q^{r}),
    & n\in\oN,\\
    \displaystyle
    -\sum_{j=n}^{-1}q^{\frac{j(2n-j)\pell}{2}}y^{j\pell}
    \sum_{r=0}^{\pell-1}x^r\,y^{r}q^{(n-j)r}
    \vartheta(q^\pell,x^\pell q^{r}),
    & n\in-\oN.
  \end{cases}
\end{multline}
There is the easily derived ``inversion'' property
\begin{gather} \label{K-inverse}
  \K_\pell(q,x,y) =
  -\K_\pell(q,x^{-1},y^{-1})+\vartheta(q^\pell,x^\pell)
  =-y^{-1}x^{-1}
  \K_\pell(q,x^{-1}q^{\frac{1}{\pell}},y^{-1}q^{-\frac{1}{\pell}}).
\end{gather}

We also note that in the exponential notation, there are the obvious
relations
\begin{alignat}{2} \label{K-identity-triv}
  \K_{\pell}(\tau, \nu + m, \mu)&=
  \K_{\pell}(\tau, \nu, \mu)=
  \K_{\pell}(\tau, \nu, \mu + m),\quad&&m\in\oZ,\\
  \label{K-identity-integer}
  \K_{\pell}(\tau, \nu + \ffrac{m}{\pell}, \mu - \ffrac{m}{\pell})&=
  \K_{\pell}(\tau, \nu, \mu), &&m\in\oZ.
\end{alignat}

\subsection{``Scaling'' formulas}\label{sec:scaling} The scaling
(``period multiplication'') formulas in this subsection are useful in
studying modular transformations of functions expressed
through~$\K_\pell$.  

We first recall the elementary theta-func\-tion identity
\begin{gather}\label{theta-rewrite}
  \vartheta(q,z)=
  \sum_{s=0}^{p-1}
  q^{\frac{s^2}{2}}\,z^s
  \vartheta(q^{p^2},z^p q^{ps})
\end{gather}
and its version for $p=2\pell u$ with \textit{coprime} $\pell$
and~$u$,
\begin{gather}\label{theta-rs}
  \vartheta(q,z)=
  \sum_{r''=1}^{2\pell}\sum_{s''=1}^{u}
  q^{\half(u r'' - \pell(s''-1))^2}\,z^{u r'' - \pell(s''-1)}
  \vartheta(q^{(2u\pell)^2},z^{2u\pell}
  q^{u\pell(u r'' - \pell(s''-1))})
\end{gather}
(where we use double-primed variables in order to help identifying
them in more complicated formulas below).

Similarly, there is an elementary identity expressing
$\K_\pell(q,x,y)$ through~$\K_\pell$ with the
``period''~$q^{u^2}$ for an arbitrary positive integer~$u$,
\begin{gather} \label{elementary}
  \K_{\pell}(q,x,y) = \!\sum_{a=0}^{u-1}\!\sum_{b=0}^{u-1}
  q^{\frac{a^2 \pell}{2} + ab} x^{a \pell + b} y^b
  \K_{\pell}(q^{u^2}\!, x^u q^{ua+b\frac{u}{\pell}},
  y^{u}q^{-b\frac{u}{\pell}}),~u\in\oN.
\end{gather}
Whenever \textit{$u$ is coprime with $\pell$}, a formula relating the
$\K_\pell$ functions with the ``periods''~$q$ and~$q^{u^2}$
differently from~\eqref{elementary} is
\begin{multline*}
  \K_\pell(q, x, y)
  =\sum_{s=0}^{u-1}\sum_{\theta=0}^{u-1}
  x^{\pell s}\,y^{\pell\theta}\,q^{\frac{s^2-\theta^2}{2}\pell}\,
  \K_\pell(q^{u^2}\!, x^u\,q^{su}, y^u\,q^{-\theta u})\\*
  {} + \smash[t]{
    \mathop{\sum_{r=1}^{\pell-1}\sum_{s=1}^{u-1}}_{u r - \pell s\geq1}}
  x^{ur}\,y^{ur-\pell s}\,q^{urs - \frac{\pell s^2}{2}}
  \vartheta(q^\pell, x^\pell\,q^{u r}),
\end{multline*}
which is shown with the help of the identity
\begin{gather} \label{up-identity}
  \ffrac{1-q^{\pell u}}{(1-q^\pell)(1-q^u)} - \ffrac{1}{1-q}
  = -\mathop{\sum_{r=1}^{\pell-1}\sum_{s=1}^{u-1}}_{u r - \pell s\geq1}
  q^{u r - \pell s}
\end{gather}
for coprime positive integers $u$ and $\pell$.

The identity in the next lemma is crucial in Sec.~\ref{sec:sl21}.  For
$n\in\oZ$, let
\begin{gather*}%%%\label{rbar}
  \rbar{\pell, u}{n} = \floor{\ffrac{\pell n}{u}}.
\end{gather*}
\begin{Lemma}\label{Lemma:psi-blow-u}
  For coprime positive integers $\pell$ and $u$,
  \begin{multline} \label{eq:K-difference-new}
    \K_{2\pell}(q^{\frac{1}{u}},x^{\frac{1}{u}},y^{\frac{1}{u}})
    -\K_{2\pell}(q^{\frac{1}{u}},x^{-\frac{1}{u}},y^{\frac{1}{u}})
    ={}\\
    \shoveleft{{}
      =\sum_{s'=1}^{u}\sum_{b=0}^{u-1}
      x^{\frac{\pell}{u}(s'-1)}
      y^{\frac{\pell}{u}(s'+1+2b)}\,
      q^{-\frac{\pell}{u}(b+1)(b+s')}
      (x y q^{-b-1})^{-\rbar{\pell, u}{s' + 2b + 1}}}\\*
    \shoveleft{\phantom{{}={}}{}\times\bigl(\K_{2\pell}(q^{u}\!,
      x q^{\frac{s'-1}{2}-\frac{u}{2\pell}\rbar{\pell, u}{
          s' + 2b + 1}}\!,
      y q^{-\frac{s'+1}{2}-b+\frac{u}{2\pell}\rbar{\pell, u}{
          s' + 2b + 1}})}
    \\*
    \shoveleft{\qquad\quad{}-x^{2\rbar{\pell, u}{s' + 2b + 1}}
      q^{(s'-1)\rbar{\pell, u}{s' + 2b + 1}}}\\*
    {}\times\K_{2\pell}(q^{u}\!,
    x^{-1}q^{-\frac{s'-1}{2}-\frac{u}{2\pell}\rbar{\pell, u}{
        s' + 2b + 1}}\!,
    y q^{-\frac{s'+1}{2} - b+\frac{u}{2\pell}\rbar{\pell, u}{
        s' + 2b + 1}})
    \bigr).
  \end{multline}
\end{Lemma}
\begin{prf}
  The formula in the Lemma is equivalently rewritten
  as\enlargethispage{12pt}
  \begin{multline*}
    \K_{2\pell}(q^{\frac{1}{u}},x^{\frac{1}{u}},y^{\frac{1}{u}})
    -\K_{2\pell}(q^{\frac{1}{u}},x^{-\frac{1}{u}},y^{\frac{1}{u}})
    ={}\\
    \shoveleft{{}
      =\sum_{s'=1}^{u}\sum_{b=0}^{u-1}
      x^{-2\frac{\pell}{u}b
        +\frac{[\pell(s'+2b-1)]_u}{u}}\,
      y^{\frac{[\pell(s'+2b-1)]_u}{u}}\,
      q^{\frac{b}{u}(\pell b - [\pell(s'+2b-1)]_u)}
    }\\*
    \shoveleft{\phantom{{}={}}{}\times
      \bigl(\K_{2\pell}(q^{u}\!,
      x q^{-b+\frac{[\pell(s'+2b-1)]_u}{2\pell}}\!,
      y q^{-\frac{[\pell(s'+2b-1)]_u}{2\pell}})}
    \\*
    \shoveleft{\qquad\quad{}-x^{2\rbar{\pell, u}{s' - 1 + 2b}}
      q^{(s'-1)\rbar{\pell, u}{s' - 1 + 2b}}}\\*
    {}\times\K_{2\pell}(q^{u}\!,
    x^{-1}q^{-b-(s'-1)+\frac{[\pell(s'+2b-1)]_u}{2\pell}}\!,
    y q^{-\frac{[\pell(s'+2b-1)]_u}{2\pell}})
    \bigr).
  \end{multline*}
  Indeed, the summand here is mapped into that
  in~\eqref{eq:K-difference-new} by a redefinition of the $b$
  variable.  This changes the $b$ summation limits, but
  Eq.~\eqref{quasi-periodicity} shows that the summand actually
  depends on $b$ only $\mathrm{mod}\,u$, and hence the interval of $u$
  consecutive values of $b$ can be translated arbitrarily.
  %%   We thus obtain~\eqref{eq:K-difference-new}.
  But the last equation can be shown directly using the
  definition~\eqref{Knew} and the fact that for coprime $\pell$ and
  $u$, \ $[\pell s']_u$ takes all the values in $[0,\dots,u-1]$ as
  $s'$ ranges over the set of any $u$ sequential values.
\end{prf}

\subsection{Relations to theta functions}\label{sec:to-theta}  Some
special combinations of the Appell functions can be expressed through
theta functions.  We first note an identity showing that the
higher-level Appell functions are expressible through $\K_1$ modulo
a ratio of theta functions (cf.\ a more general statement
in~\cite{[Pol-00]}).%\enlargethispage{18pt}

%%%\medskip

\begin{Lemma}\label{lemma:favorite}For $\pell\geq2$,
\begin{multline*}
  \vartheta(q^\pell, x)\K_\pell(q,z,y) - \sum_{r=0}^{\pell-1} z^r y^r
  \vartheta(q^\pell,z^\pell q^r)
  \K_1(q^\pell,x^{-1}, y^{\pell}q^{-r})={}\\
  {}= -\vartheta(q^\pell, y^\pell z^\pell x^{-1})\,
  \frac{ \vartheta_{1,1}(q, z
    y^{1-\pell} x)\, q^{-\frac{1}{8}}\eta(q)^3 }{ \vartheta_{1,1}(q, z y)
    \vartheta_{1,1}(q, x y^{-\pell})}.
  \end{multline*}
\end{Lemma}
This can be proved either directly (using~\eqref{veluti-periodicity}
and~\eqref{K-mn}, via resummations similar to those in
Eqs.~\eqref{ilja-ident}--\eqref{combined} below) or by noting that in
view of the open quasiperiodicity formulas, the left-hand side is in
fact \textit{quasiperiodic} in $y$ (and obviously, in~$x$ and~$z$),
and is therefore expressible as a ratio of theta functions; the actual
theta functions in this ratio are found by matching the
quasiperiodicity factors, and the remaining $q$-dependent factor is
then fixed by comparing the residues of both
sides.%%\enlargethispage{20pt}

\begin{Lemma}
  For an even level~$2\pell$,
  \begin{gather}\label{even-identity}
    \smash[b]{\kern-2pt\sum_{b=0}^{\pell-1}}
    x^{2b}q^{\frac{b^2}{\pell}}
    \bigl(\K_{2\pell}(q,x q^{\frac{b}{\pell}}, y)-
    \K_{2\pell}(q,x^{-1} q^{-\frac{b}{\pell}}, y)\bigr)
    =-\frac{
      \vartheta_{1,1}(q^{\frac{1}{\pell}},x^2)
      q^{-\frac{1}{8\pell}}\,\eta(q^{\frac{1}{\pell}})^3}{
      \vartheta_{1,1}(q^{\frac{1}{\pell}},xy)
      \vartheta_{1,1}(q^{\frac{1}{\pell}},xy^{-1})}.\kern-4pt
  \end{gather}
\end{Lemma}
To prove this, we use the same strategy as above, the crucial point
being quasiperiodicity, which is shown as follows.  With
$\Delta_\pell\,f(q,x,y)$ used to temporarily denote
$f(q,x,yq)-q^{\frac{\pell}{2}}y^{\pell}\,f(q,x,y)$, it follows from
Eq.~\eqref{veluti-periodicity} that
\begin{multline}\label{even-derive2}
  \Delta_\pell\bigl(\K_{2\pell}(q,xq^{\frac{b}{\pell}},y)
  -\K_{2\pell}(q,x^{-1}q^{-\frac{b}{\pell}},y)\bigr)=\\
  {}= \sum_{a=1}^{\pell-1}
  x^{-a}q^{-\frac{ab}{\pell}+a}y^{a}
  (y^{2\pell-2a}q^{\pell-a}-1)\vartheta(q^{2\pell},x^{2\pell}q^{2b-a})\\*
  {}+
  \sum_{a=1}^{\pell-1} x^{a}q^{\frac{ab}{\pell}+a}y^a
  (1-y^{2\pell-2a}q^{\pell-a})\vartheta(q^{2\pell},x^{2\pell}q^{2b+a}).
\end{multline}
This also shows that $x^{2b}\,q^{\frac{b^2}{\pell}}
\Delta_\pell\K_{2\pell}(q,x^{-1} q^{-\frac{b}{\pell}}, y)$ depends
on~$b$ only modulo~$\pell$.  In applying $\sum_{b=0}^{\pell-1}
x^{2b}\,q^{b^2/\pell}$ to the second term in the right-hand side
of~\eqref{even-derive2}, we can therefore make the shift $b\mapsto
b-a$ without changing the summation limits for~$b$.  This readily
implies that the left-hand side of~\eqref{even-identity} is
quasiperiodic in~$y$.

\subsection{Integral representation and the $S$-transform
  of~$\K_\pell$}\label{sec:integral-rep} Although the Appell functions
cannot be rationally expressed through theta functions, they admit an
\textit{integral} representation through a ratio of theta functions.
This integral representation proves to be a useful tool, in particular
in finding modular transformations of~$\K_\pell$.  We give this
representation in~\eqref{K-integral} and then use it in the
calculation leading to Theorem~\bref{thm:TS}.

\begin{Lemma}\label{lemma:integral}
  The Appell function admits the integral representation
  \begin{gather}\label{K-integral}
    \K_\pell(\tau,\nu,\mu)
    =-\EX{-\frac{i\pi}{4}\tau}
    \int_0^1\!d\lambda\,
    \vartheta(\pell\tau, \pell\nu\!-\!\lambda)\,\frac{
      \vartheta_{1,1}(\tau, \nu\!+\!\mu\!+\!\lambda)\,
      \eta(\tau)^3}{
      \vartheta_{1,1}(\tau, \nu\!+\!\mu)
      \vartheta_{1,1}(\tau, \lambda\!+\!i0)},
  \end{gather}
  where $+i0$ specifies the contour position to bypass the
  singularities.
\end{Lemma}
\begin{prf}
  Starting with the easily derived identity
  \begin{gather}\label{ilja-ident}
    \sum_{m\in\oZ}\K_\pell(q,z,yq^{m})x^m =
    \vartheta(q^\pell, z^\pell x^{-1})\sum_{m\in\oZ}
    \frac{x^m}{1 - y z q^{m}}
  \end{gather}
  and combining it with the identity~\cite{[KW],[FSST]}
  \begin{gather}\label{truly}
    \sum\limits_{m\in\oZ}\frac{x^m}{1-yq^m}
    =-\frac{\vartheta_{1,1}(q,xy)\,\prod_{i\geq1}(1-q^i)^3}{
      \vartheta_{1,1}(q,y)\vartheta_{1,1}(q,x)},  
  \end{gather}
  we obtain
  \begin{gather}\label{combined}
    \sum_{m\in\oZ}\K_\pell(q,z,y q^{m})x^m
    =-\vartheta(q^\pell, z^\pell x^{-1})\,\frac{
      \vartheta_{1,1}(q, z y x)\,
      q^{-\frac{1}{8}}\eta(q)^3
    }{
      \vartheta_{1,1}(q, z y)
      \vartheta_{1,1}(q, x)}.
  \end{gather}
  The right-hand side of~\eqref{truly} is a meromorphic function of $x$
  with poles at $x=q^n$, $n\in\oZ$, but the identity holds in the
  annulus $\boldsymbol{A}_1=\{x\bigm||q|<|x|<1\}$, where the left-hand
  side converges.  We therefore temporarily assume that
  $x\in\boldsymbol{A}_1$ and then analytically continue the final
  result.  Integrating over a closed contour inside this annulus yields
  \begin{gather}\label{K-integral-z}
    \K_\pell(q,z,y)
    =-\frac{1}{2i\pi }\oint \frac{d x}{x}\,
    \vartheta(q^\pell, z^\pell x^{-1})\,\frac{
      \vartheta_{1,1}(q, z y x)\,
      q^{-\frac{1}{8}}\eta(q)^3
    }{
      \vartheta_{1,1}(q, z y)
      \vartheta_{1,1}(q, x)}.
  \end{gather}

  In the exponential notation $z\!=\!\ex{2i\pi\nu}$,
  $y\!=\!\ex{2i\pi\mu}$, $x\!=\!\ex{2i\pi\lambda}$, the annulus
  $\boldsymbol{A}_1%\{x\bigm| |q|<|x|<1\}
  $ is mapped into any of the parallelograms $\boldsymbol{P}_n$,
  $n\in\oZ$, with the vertices $(n, n+1, n+1 + \tau, n + \tau)$.  We
  choose $n=0$ in what follows.  In the exponential notation, the
  integration contour is then mapped into a contour in the interior of
  $\boldsymbol{P}_0$ connecting the points in a close vicinity of $0$
  and $1$ respectively.  Equation~\eqref{K-integral} thus follows.
\end{prf}

The integral representation in the Lemma allows us to find the
$S$-transform of $\K_\pell$.  (As regards the $T$ transformation, it
readily follows that Eq.~\eqref{Knew-Ttransf} holds for $\K_\pell$.)
For this, we use the known $S$-transformation properties of the $\eta$
and $\vartheta$ functions entering~\eqref{K-integral}, with the result
\begin{multline} \label{K-work-2}
  \K_\pell(-\ffrac{1}{\tau},\ffrac{\nu}{\tau},\ffrac{\mu}{\tau})
  =-\sqrt{\ffrac{-i\tau}{\pell}}\,\tau
  \int_0^1\!d\lambda\,
  \EX{i\pi(\frac{\pell\nu^2}{\tau} +
    \frac{\lambda^2\tau}{\pell} + 2\lambda\mu)}\\*
  {}\times
  \sum_{r=0}^{\pell-1}
  \EX{2i\pi r(\nu - \frac{\tau}{\pell}\lambda)
    +i\pi\frac{r^2}{\pell}\tau}
  \vartheta(\pell\tau,\pell\nu\!-\!\tau\lambda\!+\!r\tau)
  \frac{
    \vartheta_{1,1}(\tau, \nu\!+\!\mu\!+\!\tau\lambda)\,
    \ex{-\frac{i\pi\tau}{4}}\eta(\tau)^3}{
    \vartheta_{1,1}(\tau,\nu\!+\!\mu)
    \vartheta_{1,1}(\tau,\tau\lambda\!-\!\varepsilon)},
\end{multline}
where we also used~\eqref{theta-rewrite} (with $u=\pell$) to rewrite
the theta function $\vartheta(\frac{\tau}{\pell},\cdot)$ occurring in
the $S$-transform of~\eqref{K-integral} (and where infinitesimal
positive $\varepsilon$ specifies the contour position).
Equation~\eqref{K-work-2} allows applying Lemma~\bref{lemma:favorite}.
As a result, after some additional simple transformations
involving~\eqref{diag-periodicity} and~\eqref{quasi-periodicity}, the
second line in~\eqref{K-work-2} becomes
\begin{multline*}
  \sum_{r=0}^{\pell-1}
  \EX{i\pi\frac{r^2}{\pell}\tau
    \!-\! 2i\pi r\mu
    \!-\! 2i\pi\frac{r}{\pell}\tau\lambda}
  \Bigl(\vartheta(\pell\tau,\tau\lambda + \pell\mu - r\tau)
  \K_\pell(\tau,\nu,\mu)\\*
  {}- \sum_{a=0}^{\pell-1}
  \EX{2i\pi a(\nu\!+\!\mu)}
  \vartheta(\pell\tau,\pell\nu+a\tau)
  \K_1(\pell\tau,-\tau\lambda - \pell\mu + r\tau,
  \pell\mu - a\tau)\Bigr).
\end{multline*}
In the first term here, we next use~\eqref{theta-rewrite}, which gives
the integral in~\eqref{theta-b}, and in integrating the second term,
we change the integration variable as
$\lambda\to\frac{\lambda+r}{\tau}$, which then allows us to do the $r$
summation explicitly.  This gives
\begin{multline} \label{K-work-3}
  \K_\pell(-\ffrac{1}{\tau},\ffrac{\nu}{\tau},\ffrac{\mu}{\tau})=
  \tau\,\EX{i\pi\pell\frac{\nu^2-\mu^2}{\tau}}
  \K_\pell(\tau,\nu,\mu)\\*
  {}-\sqrt{\ffrac{-i\tau}{\pell}}  
  \sum_{a=0}^{\pell-1}
  \vartheta(\pell\tau,\pell\nu+a\tau)\,
  \EX{i\pi\pell\frac{(\nu+\frac{a}{\pell}\tau)^2}{\tau}}
  \int\limits_{-\tau-a\tau}^{(\pell-1-a)\tau}\!d\lambda\,
  \EX{i\pi\frac{\lambda^2}{\pell\tau}
    - 2i\pi\frac{\lambda(\pell\mu-a\tau)}{\pell\tau}}\\*
  {}\times\K_1(\pell\tau, \lambda - i0 -(\pell\mu-a\tau),
  \pell\mu - a\tau),  
\end{multline}
showing that the remaining integral is the one in~\eqref{Phi-K-shift}.
This leads to the sought equation expressing the $S$-transform of
$\K_\pell$ through a single $\K_\pell$ function and~$\pell$ theta
functions, Eq.~\eqref{Knew-Stransf-fin}, with
\begin{equation}\label{phi-through-K}
  \begin{split}
    \Phi(\tau,\mu)={}&\phi(\tau,\mu)
    -\ffrac{i}{2\sqrt{-i\tau}},\notag\\
    \phi(\tau,\mu)={}&\ffrac{i}{\sqrt{-i\tau}}\, \dashint_{0}^{\tau}
    d\lambda\, \EX{i\pi\frac{\lambda^2 - 2\lambda\mu}{\tau}}
    \K_1(\tau,\lambda-\mu,\mu)
  \end{split}
\end{equation}
involved in the theta-function terms (see Appendix~\ref{app:ab-int}
for the principal-value integral).  Equivalently, the $\phi$ function
can be rewritten as
\begin{align}
  \phi(\tau,\mu)
  ={}&-\fhalf\int_{\oR}\!dx\,\EX{-\pi x^2}\,
  \frac{\sinh\!\left(\pi
      x\sqrt{-i\tau}(1+2\frac{\mu}{\tau})\right)}{
    \sinh\!\left(\pi x\sqrt{-i\tau}\right)}\label{phi-def}\\
  \intertext{and $\Phi$ as} \Phi(\tau,\mu)
  ={}&-\int\limits_{\oR-i0}\!dx\,\EX{-\pi x^2}\,
  \frac{\ex{-2i\pi
      x\frac{\mu}{\sqrt{-i\tau}}}}{ 1-\ex{-2\pi
      x\sqrt{-i\tau}}}.\label{R-i0}
\end{align}
This proves the formula for the $S$-transform of~$\K_\pell$ in
Theorem~\bref{thm:TS}.  The integral is to be analytically continued
from $\tau=it$ with $t\in\oR_{>0}$.

\subsection{The lowest-level Appell functions}\label{sec:lowest} 
Appell functions of levels~$1$ and~$2$ have some special or simplified
properties.  For $\K_1$, Eq.~\eqref{Knew-Stransf-fin} becomes
(cf.~\cite{[Polisch]})
\begin{gather}
  \label{K-Stransf-p1}
  \K_1(-\ffrac{1}{\tau},\ffrac{\nu}{\tau},\ffrac{\mu}{\tau})
  = \tau \EX{i\pi\frac{\nu^2 - \mu^2}{\tau}}\,
  \K_1(\tau,\nu,\mu)
  + \tau\,\EX{i\pi\frac{\nu^2}{\tau}}
  \Phi(\tau,\mu)
  \vartheta(\tau, \nu).
\end{gather}
For $\pell=1$, the formula in Lemma~\bref{lemma:favorite} becomes an
identity in~\cite{[Polisch]},
\begin{gather*}
  \vartheta(q,z)\K_1(q,x,y) - \vartheta(q,x)\K_1(q,z,y)
  =\frac{\vartheta(q,xyz)\,\vartheta_{1,1}(q,x^{-1}z)
    \,q^{-\frac{1}{8}}\,\eta(q)^3}{
    \vartheta_{1,1}(q,x^{-1}y^{-1})\,\vartheta_{1,1}(q,yz)}.
\end{gather*}
For $\K_2$, Eq.~\eqref{even-identity} simplifies to
\begin{gather}\label{even-identity2}
  \K_{2}(q, x, y) - \K_{2}(q, x^{-1}, y)
  =-\frac{\vartheta_{1,1}(q,x^2)\,q^{-\frac{1}{8}}\,\eta(q)^3}{
    \vartheta_{1,1}(q,xy)\vartheta_{1,1}(q,xy^{-1})}.
\end{gather}

\subsection{The $\Phi$ function and its properties}\label{sec:Phi}
We now study the properties of the function $\Phi$ appearing in the
$S$ transform of higher-level Appell functions.  These properties
include open quasiperiodicity relations\,---\,which can be viewed as
functional equations satisfied by~$\Phi$\,---\,and a modular
transformation formula.  They already follow from
Eq.~\eqref{K-Stransf-p1}, or alternatively, can be derived from the
integral representation~\eqref{phi-def}, similarly to the study of
Barnes-like special functions arising in various
problems~\cite{[B],[FK]} (see also~\cite{[PT],[FKV],[JM],[StsKL]}).
Unlike Barnes-like functions, however, the $\Phi$ function cannot be
evaluated as a sum of residues (tentatively, at
$x_n=i\frac{n}{\sqrt{-i\tau}}$, $n\in\oZ_{\geq0}$) of the
integral~\eqref{phi-def}, because the Gaussian exponential causes the
sum to diverge.

\subsubsection{Open quasiperiodicity and related properties}
First, a simple calculation allows explicitly evaluating
$\phi(\tau,\mu)$ in~\eqref{phi-def} for $\mu=\frac{m\tau}{2}$,
$m\in\oZ$:
\begin{gather}\label{explicit-phi}
  \begin{aligned}
    \phi(\tau,\ffrac{m\tau}{2})&= -\fhalf
    \sum_{j=0}^m
    \ex{-i\pi\tau\frac{(m-2j)^2}{4}},\quad& m&{}\geq0,\\
    \phi(\tau,-\ffrac{m\tau}{2})&= \fhalf
    \sum_{j=1}^{m-1}
    \ex{-i\pi\tau\frac{(m-2j)^2}{4}},
    \quad& m&{}\geq1.
  \end{aligned}
\end{gather}
In particular, $\phi(\tau,-\tau/2)=0$.\enlargethispage{12pt}

Next, elementary transformations with the integral representation,
involving the identity
\begin{gather*}
  \ffrac{
    \EX{2\pi x m \sqrt{-i\tau}} - 1}{
    2\sinh(\pi x\sqrt{-i\tau})}=
  \sum_{j=0}^{m-1} \EX{\pi x\sqrt{-i\tau} (2j+1)},
  \qquad m\in\oN,
\end{gather*}
show that $\Phi$ satisfies the equations
\begin{gather}\label{phi-tau}
  \begin{split}
    \Phi(\tau,\mu\!+\!m\tau) &= \Phi(\tau,\mu)
    - \sum_{j=1}^{m}
    \ex{-i\pi\frac{(\mu+j\tau)^2}{\tau}},\\
    \Phi(\tau,\mu\!-\!m\tau) &= \Phi(\tau,\mu)
    + \sum_{j=0}^{m-1}
    \ex{-i\pi\frac{(\mu - j\tau)^2}{\tau}},
  \end{split}
  \qquad m\in\oN.
\end{gather}
Similarly to the equations for $\K_\pell$, these are open
quasiperiodicity relations.  They can be alternatively derived
from~\eqref{K-Stransf-p1} and the corresponding
property~\eqref{veluti-periodicity} of the Appell functions.  For
this, we evaluate the commutator of the $S$~transform
of~$\K_\pell(\tau,\nu,\mu)$ and the translation of the $\mu$ argument
by elements of the lattice generated by $(1,\tau)$; because modular
transformations act on lattice translations (thus forming the
semidirect product), this results in equations for~$\Phi$, equivalent
to~\eqref{phi-tau}.

\begin{figure}[tb]
  \unitlength=.8pt
  \begin{center}
  \begin{picture}(400,200)(0,0)
    \put(205,182){$i\frac{m}{\sqrt{t}}$}
    \put(205,150){$i\frac{m-1}{\sqrt{t}}$}
    \put(205,60){$i\frac{1}{\sqrt{t}}$}
%%%    \put(390,35){$x$}
    \dottedline[$\boldsymbol{\times}$]{30}(200,0)(200,180)
    \thicklines
    \dashline{5}(0,175)(400,175)
    \put(300,175){\vector(1,0){5}}
    \dashline{5}(0,25)(400,25)
    \put(300,25){\vector(1,0){5}}
    \Thicklines
    \put(0,30){\vector(1,0){400}}
    \put(200,-10){\vector(0,1){210}}
  \end{picture}
\end{center}
\caption[Integration contours]{\small Integration contours for
  $\Phi(it,\mu)$ (the lower dashed line) and
  $\Phi_{\frac{im}{\sqrt{t}}}(it,\mu)$ (the upper dashed line) in the
  complex $x$ plane and poles of the integrand (crosses).}
    \label{fig:residues}
\end{figure}
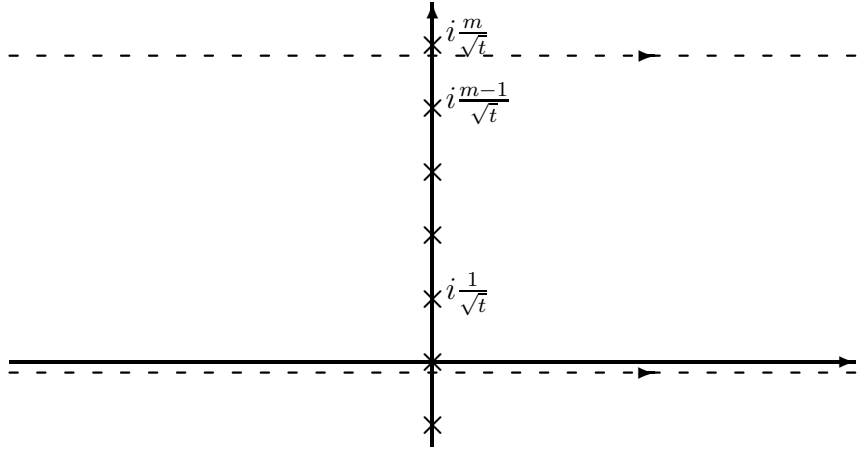

A slightly more involved calculation with the integral representation
leads to the ``dual'' open quasiperiodicity relations
\begin{gather}\label{phi-1}
  \begin{split}      
    \Phi(\tau, \mu\!+\!m) &=
    \EX{-i\pi\frac{m^2}{\tau}
      \!-\!2i\pi
      m\frac{\mu}{\tau}}\,\Phi(\tau,\mu)
    + \ffrac{i}{\sqrt{-i\tau}}
    \sum_{j=1}^m
    \EX{i\pi\frac{j(j-2m)}{\tau}
      \!-\!2i\pi j\frac{\mu}{\tau}},\\
    \Phi(\tau,\mu\!-\!m) &=
    \EX{-i\pi\frac{m^2}{\tau}
      \!+\!2i\pi
      m\frac{\mu}{\tau}}\,\Phi(\tau,\mu)
    - \ffrac{i}{\sqrt{-i\tau}}    
    \sum_{j=0}^{m-1}
    \EX{i\pi\frac{j(j-2m)}{\tau}
      \!+\!2i\pi j\frac{\mu}{\tau}},
  \end{split}
\end{gather}
where $m\in\oN$. \ To show this, we recall the analytic continuation
prescription and write $\Phi(it,\mu)$ with $t\in\oR_+$ as
in~\eqref{R-i0}.  We then consider $\Phi(it,\mu + m)$ with $m\in\oN$
and change the integration variable as $x=x' - i\frac{m}{\sqrt{t}}$.
This gives
\begin{gather*}
  \Phi(it,\mu\!+\!m) =
  \EX{-\pi\frac{m^2}{t} - 2\pi\frac{m\mu}{t}}
  \Phi_{\frac{im}{\sqrt{t}}}(it,\mu),
\end{gather*}
where $\Phi_{\frac{im}{\sqrt{t}}}(it,\mu)$ is given by the integral
along $\oR+i\frac{m}{\sqrt{t}}-i0$ of \textit{the same integrand} as
for $\Phi(it,\mu)$.  A residue calculation in accordance with
\begin{gather*}
  \Phi_{\frac{im}{\sqrt{t}}}(it,\mu) =
  \Phi(it,\mu) - 2i\pi\sum_{n=0}^{m-1}\res_{x=i\frac{n}{\sqrt{t}}}
  \biggl(-\EX{-\pi x^2}\,
  \ffrac{\EX{-2i\pi x\frac{\mu}{\sqrt{t}}}}{
    1-\EX{-2\pi x\sqrt{t}}}\biggr)
\end{gather*}
(see Fig.~\ref{fig:residues}) then yields the first equation
in~\eqref{phi-1}.
Alternatively, Eqs.~\eqref{phi-1} can also be deduced
from~\eqref{K-Stransf-p1} and the corresponding Appell function
property in Eq.~\eqref{K-identity-triv}.  

Next, a ``reflection property'' follows from~\eqref{K-inverse} (or can
be directly derived from~\eqref{phi-def}),
% (also consistent with~\eqref{explicit-phi})
\begin{align}\label{Phi-reflect}
  \Phi(\tau,-\mu) &=
  \ffrac{-i}{\sqrt{-i\tau}}
  - \EX{-i\pi\frac{\mu^2}{\tau}} - \Phi(\tau,\mu),\\
  \intertext{or equivalently,}
  \label{Phi-reflect2}
  \Phi(\tau,\mu+1) &=
  \smash[t]{-\EX{-2i\pi\frac{\mu + \half}{\tau}}} \Phi(\tau,-\mu-\tau).
\end{align}

A simple ``scaling law'' for $\Phi$ follows from~\eqref{elementary}
and~\eqref{K-Stransf-p1} (or directly from integral
representation~\eqref{phi-def}),
\begin{gather}\label{phi-u}
  \Phi(\tau,\mu) = \sum_{b=0}^{p-1}\Phi(p^2 \tau, p \mu - b p
  \tau),\qquad p\in\oN.
\end{gather}
In the case of ``scaling'' with an even factor, we have the following
two formulas.  As before, $[x]_p=x\,\mathrm{mod}\,p$.

\begin{Lemma}
  For any $\pell\in\oN$ and  $m\in\oZ$,
  \begin{gather}\label{phi-2p}
  \sum_{a=0}^{2\pell-1}\EX{i\pi\frac{am}{\pell}}
  \Phi(2\pell\tau,2\pell\mu\!-\!a\tau)
  =\EX{i\pi\frac{[m]_{2\pell}^2}{2\pell\tau}
    \!+\!2i\pi\frac{\mu}{\tau}[m]_{2\pell}}
  \Phi(\ffrac{\tau}{2\pell},\mu\!+\!\ffrac{[m]_{2\pell}}{2\pell}),\\
  \ffrac{1}{2\pell}\sum_{a=0}^{2\pell-1}
  \EX{
    2i\pi\frac{a\mu}{\tau}
    \!-\!i\pi\frac{am}{\pell}
    \!+\!i\pi\frac{a^2}{2\pell\tau}}
  \Phi(\ffrac{\tau}{2\pell},\mu\!+\!\ffrac{a}{2\pell})
  =\Phi(2\pell\tau,2\pell\mu\!-\![m]_{2\pell}\,\tau).
  \label{dual-scaling-m}
\end{gather}
\end{Lemma}
For $m=0$, Eq.~\eqref{phi-2p} reduces to~\eqref{phi-u}, and
Eq.~\eqref{dual-scaling-m} becomes
\begin{gather}\label{dual-scaling}
  \ffrac{1}{2\pell}\sum_{a=0}^{2\pell-1}
  \EX{2i\pi\frac{a\mu}{\tau}
    \!+\!i\pi\frac{a^2}{2\pell\tau}}
  \Phi(\ffrac{\tau}{2\pell},\mu\!+\!\ffrac{a}{2\pell})
  =\Phi(2\pell\tau,2\pell\mu).
\end{gather}
\begin{prf}
  Both formulas~\eqref{phi-u} and~\eqref{dual-scaling-m} follow from
  manipulations with the integral representation, but as we see
  shortly, these two equations are $S$-dual to each other, and it
  therefore suffices to prove any one of them.  We actually
  show~\eqref{dual-scaling}, from which the rest (all of the $m$
  dependence) follows via~\eqref{phi-tau}.  For this, we evaluate the
  left-hand side using the integral representation and change the
  integration variable $x$ in the $a$th term as
  $x=x'-i\frac{a}{2\pell\sqrt{-i\tau}}$; as before, we analytically
  continue from the positive part of the imaginary axis $\tau=it$.
  Then (omitting the prime at~$x$)
  \begin{gather*}
    \ffrac{1}{2\pell}\sum_{a=0}^{2\pell-1}
    \EX{i\pi\frac{a\mu}{\pell\tau}
      \!+\!i\pi\frac{a^2}{4\pell^2\tau}}
    \Phi(\tau,\mu\!+\!\ffrac{a}{2\pell})
    =-\ffrac{1}{2\pell}\sum_{a=0}^{2\pell-1}
    \int\limits_{\oR+\frac{i a}{2\pell\sqrt{t}}-i0}\!\!dx\,
    \EX{-\pi x^2}\frac{\ex{\pi x\sqrt{t}\,\frac{2\mu}{\tau}}}{
      1\!-\!\ex{i\pi\frac{a}{\pell}\!-\!2\pi x\sqrt{t}}}.
  \end{gather*}
  But the integration contour can be deformed to $\oR\!-\!i0$ in
  \textit{all} terms, which allows us to evaluate the sum over $a$ as
  \begin{gather*}
    \ffrac{1}{2\pell}
    \sum_{a=0}^{2\pell-1}\mfrac{1}{1-\ex{i\pi\frac{a}{\pell}}q} 
    =\mfrac{1}{1-q^{2\pell}},
  \end{gather*}
  which shows~\eqref{dual-scaling}.  
\end{prf}

We also note that combining Eqs.~\eqref{theta-rewrite}, \eqref{phi-u},
and~\eqref{phi-1} allows expressing the $S$-transform of $\K_\pell$,
Eq.~\eqref{Knew-Stransf-fin}, as
\begin{multline} \label{Knew-Stransf-fin2}
  \K_\pell(-\ffrac{1}{\tau},\ffrac{\nu}{\tau},\ffrac{\mu}{\tau}) =
  \tau\,\EX{i\pi\pell\frac{\nu^2-\mu^2}{\tau}}\K_\pell(\tau,\nu,\mu)
  \\*
  {}+\ffrac{\tau}{\pell}\sum_{a=0}^{\pell-1}
  \EX{i\pi\frac{\pell}{\tau}\nu^2\!+\!i\pi\frac{a^2}{\pell\tau}
    \!+\!2i\pi a\frac{\mu}{\tau}}
  \Phi(\ffrac{\tau}{\pell},\mu\!+\!\ffrac{a}{\pell})\,
  \vartheta(\ffrac{\tau}{\pell},\nu\!-\!\ffrac{a}{\pell}).
\end{multline}
With some care, Eq.~\eqref{Knew-Stransf-fin2} can be derived directly
via the Poisson resummation formula.

\subsubsection{Modular transformations of the $\Phi$ function}Modular
group relations impose constraints on the function $\Phi$ appearing in
the $S$-trans\-form of higher-level Appell functions.  With the action
of $C=S^2$ given by Eq.~\eqref{Phi-reflect}, we act with both sides of
the relation $(ST)^3=C$ on $\K_1$.  Comparing the results gives the
identity
\begin{multline}\label{Phi-T}
  i\sqrt{-i\tau}\,\EX{i\pi\frac{\mu^2}{\tau}}
  \Phi(\tau, \mu)
  + \EX{i\pi\frac{(\mu+\half)^2}{\tau-1}}
  \Phi(1\!-\!\ffrac{1}{\tau}, \ffrac{\mu}{\tau}\!+\!\fhalf)\\*
  {}- i\sqrt{-i\tau}\,\EX{i\pi\frac{(\mu+\half)^2}{\tau-1}}
  \Phi(\tau\!-\!1, \mu\!+\!\fhalf) = 0.
\end{multline}
Next, with the action of $S$ given by Eq.~\eqref{K-Stransf-p1}, the
$\SLiiZ$ relation $S^2=C$ results in the $S$-transformation formula
for $\Phi$,
\begin{gather}\label{Phi-S}
  \Phi(-\ffrac{1}{\tau},\ffrac{\mu}{\tau})
  = -i\sqrt{-i\tau}\bigl(
  \EX{i\pi\frac{\mu^2}{\tau}}
  \Phi(\tau,\mu) + 1\bigr)
  = -i\sqrt{-i\tau}\EX{i\pi\frac{\mu^2}{\tau}}
  \Phi(\tau,\mu\!-\!\tau).
\end{gather}
It now follows that identities~\eqref{phi-2p}
and~\eqref{dual-scaling-m} are the $S$-transform of each other.

It is instructive to verify (the first equality in)~\eqref{Phi-S} by
comparing the asymptotic expansions of the integral in~\eqref{phi-def}
as $-i\tau\!\searrow0$ and $\tau\to\infty$.  We first find the
asymptotic form of $\phi(it,\mu)$ for large positive $t$.  Writing
\begin{align*}
  \phi(it,iy)
  &\asymp-\int_{\epsilon}^{+\infty}\!dx\,\EX{-\pi x^2}\,
  \ffrac{\sinh\!\left(\pi x\sqrt{t}(1+2\frac{y}{t})\right)}{
    \sinh\!\left(\pi x\sqrt{t}\right)}\\
  &\stackrel{t\to+\infty}{\asymp} -2\sum_{m=0}^{\infty}
  \int_{\epsilon}^{+\infty}\!dx\,\EX{-\pi x^2}\,
  \sinh\Bigl(\pi x(\sqrt{t} + 2\ffrac{y}{\sqrt{t}})\Bigr)\,
  \EX{-\pi x(2m\!+\!1)\sqrt{t}},
\end{align*}
we readily obtain
\begin{multline*}  
  \phi(it,iy)\stackrel{t\to\infty}{\asymp}
  \fhalf\sum_{m=0}^\infty\EX{\pi\frac{((m+1) t + y)^2}{t}}
  \Erfc\Bigl(\sqrt{\ffrac{\pi}{t}}
  ((m\!+\!1) t\!+\!y\!+\!\sqrt{t}\epsilon)
  \Bigr)\\*
  {}-\fhalf\sum_{m=0}^\infty\EX{\pi\frac{(y - m t)^2}{t}}
  \Erfc\Bigl(\sqrt{\ffrac{\pi}{t}}(m t\!-\!y\!+\!\sqrt{t}\epsilon)
  \Bigr),
\end{multline*}
with the complementary error function
\begin{gather*}
  \Erfc(z) = 1 - \Erf(z)
  =1-\ffrac{2}{\sqrt{\pi}} \int_0^z dt \ex{-t^2}.
\end{gather*}
We isolate the $m=0$ term in the first sum (the only term where the
argument of $\Erfc$ is small as $t\to\infty$), rearrange the remaining
series, and use the asymptotic expansion
\begin{gather*}
  \Erfc(z) \asymp \ffrac{1}{\sqrt{\pi}}\,\EX{-z^2}
  \sum_{n\geq0}(-1)^n \ffrac{(2n\!-\!1)!!}{2^n\,z^{2n+1}}
\end{gather*}
for large positive $z$.  We can then set $\epsilon=0$, which gives
\begin{multline*}
  \phi(it,iy)\stackrel{t\to\infty}{\asymp}
  -\fhalf\, \EX{\pi\frac{y^2}{t}}
  \Erfc(-\sqrt{\ffrac{\pi}{t}}y)
  + \fhalf\bigl(\cot(\ffrac{\pi y}{t}) - \ffrac{t}{\pi y}\bigr)\\*
  + \smash[b]{\sum_{n=1}^\infty}
  \ffrac{(-1)^n\,(2 n\!-\!1)!!}{(2 \pi)^{n + 1}}\,
  \ffrac{\zeta(2 n\!+\!1, 1\!+\!\frac{y}{t})
    - \zeta(2 n\!+\!1, 1\!-\!\frac{y}{t})}{t^{n + \frac{1}{2}}},
\end{multline*}
where
\begin{align*}
  \zeta(s, a) &= \smash{\sum_{m = 0}^{\infty}}(m + a)^{-s}.\\
\intertext{Using}
  \zeta(n, 1 + x) &=
  \sum_{i\geq0}(-1)^i\,
  \fbinom{n + i - 1}{n - 1}
  x^i \zeta(n + i),\\
\intertext{where $\zeta(s)=\sum_{m\geq1}m^{-s}$, and expressing the
$\zeta$-function values at even positive integers through the
Bernoulli numbers $B_{2n}$ as}
  \zeta(2n) &=
  \pi^{2n}\,\ffrac{2^{n - 1} (-1)^{n+1}%%%{Mod[n, 2] + 1}
    }{n!\,(2n\!-\!1)!!}\,B_{2n},
\end{align*}
we then obtain the large-$t$ asymptotic expansion
\begin{gather}\label{asymp-infty}
  \phi(it,iy)\stackrel{t\to\infty}{\asymp}
  -\fhalf\, \EX{\pi\frac{y^2}{t}} \Erfc(-\sqrt{\ffrac{\pi}{t}}y)
  -\sum_{n=0}^{\infty}\sum_{j=0}^{\infty}
  \ffrac{(-4)^j \pi^{2 j + n + 1} 
    B_{2 (j\!+\!n\!+\!1)}}{(j\!+\!n\!+\!1)\,(2 j\!+\!1)!\,n!}\,
  \ffrac{y^{2 j + 1}}{t^{\frac{3}{2} + 2 j + n}}.
\end{gather}
The $\Erfc$ function in~\eqref{asymp-infty} can be expanded further,
with the result
\begin{gather*}
  \phi(it,iy)\stackrel{t\to\infty}{\asymp}
  -\fhalf\, \EX{\pi\frac{y^2}{t}}
  -\sum_{n=0}^{\infty}
  \ffrac{\pi^{n} 
    B_{2 n}}{n!\,t^{\frac{1}{2} + n}}\,
  y\HyperF11\bigl(1-n,\ffrac{3}{2},\ffrac{\pi y^2}{t}\bigr).
\end{gather*}

\smallskip

We next find the small-$t$ expansion.\footnote{We are grateful to
  V.I.~Ritus for the elegant derivation in
  Eqs.~\eqref{asymp-0-start}--\eqref{asymp-0}.}  Writing
\begin{gather}\label{asymp-0-start}
  \phi(it, iy)
  =-\fhalf\int_{-\infty}^{+\infty}
  \!dx\,\EX{-\pi x^2}
  \Bigl(
  \cosh\bigl(2\pi x\ffrac{y}{\sqrt{t}}\bigr)
  + \coth\bigl(\pi x\sqrt{t}\bigr)
  \sinh\bigl(2\pi x\ffrac{y}{\sqrt{t}}
  \bigr)\Bigr),
\end{gather}
we calculate the $\cosh$ integral and expand $\coth$, with the result
\begin{gather*}%%%\label{asymp-0-cont}
  \phi(it , iy)
  \stackrel{t\to0}{\asymp}
  -\fhalf\,\EX{\frac{\pi y^2}{t}} -
  \sum_{j=0}^\infty
  \ffrac{2^{2j} B_{2j}}{(2j)!}
  \int_0^\infty\!dx\,\EX{-\pi x^2}(\pi x\sqrt{t})^{2j-1}
  \sinh\bigl(2\pi x\ffrac{y}{\sqrt{t}}
  \bigr).
\end{gather*}
This involves the integrals
\begin{multline*}
  \begin{aligned}
    \int_0^\infty\!dx\,\EX{-\pi x^2}x^{-1}\sinh(\beta x)&=
    -\ffrac{i\pi}{2}\,\Erf\bigl(\ffrac{i\beta}{2\sqrt{\pi}}\bigr),
    \\
    \int_0^\infty\!dx\,\EX{-\pi x^2}x^{2j-1}\sinh(\beta x)&= (-1)^j
    \ffrac{i\sqrt{\pi}}{(4\pi)^j}\,
    \EX{\frac{\beta^2}{4\pi}}\,H_{2j-1}\bigl(
    \ffrac{i\beta}{2\sqrt{\pi}}\bigr),
    \quad j\geq1,
  \end{aligned}
\end{multline*}
where $H_m$ are the Hermite polynomials.  They can be written as
\begin{gather*}
  H_{m}(x) = (2x)^m\left(
    1
    - \fbinom{m}{2}\ffrac{1}{2x^2}
    + 1\cdot3 \fbinom{m}{4}\ffrac{1}{(2x^2)^2}
    - 1\cdot3\cdot5\fbinom{m}{6}\ffrac{1}{(2x^2)^3}+\dots
    \right),
\end{gather*}
which gives the small-$t$ expansion
\begin{multline}\label{asymp-0}
  \phi(it, iy)\stackrel{t\to0}{\asymp}
  -\fhalf\,\EX{\pi\frac{y^2}{t}}
  + \ffrac{i}{2\sqrt{t}}\,
  \Erf\Bigl(i\sqrt{\ffrac{\pi}{t}}\,y\Bigr)
  \\*  
  -\EX{\pi\frac{y^2}{t}}
  \sum_{n=0}^\infty\!\sum_{j=0}^\infty
  \ffrac{B_{2(j + n + 1)}\,4^{ j }\pi^{2j+n+1}}{(j\!+\!n\!+\!1)\,n!\,
    (2j\!+\!1)!}\,
  y^{2j + 1}\,t^n.
\end{multline}

To verify consistency with the $S$-transform formula~\eqref{Phi-S}, we
rewrite it as
\begin{gather*}
  \phi(\ffrac{i}{t},\ffrac{y}{t})
  + i\sqrt{t}\,\EX{-\pi\frac{y^2}{t}}\,\phi(i t, i y) =
  -\ffrac{i}{2}\,\sqrt{t} - \fhalf\,\EX{-\pi\frac{y^2}{t}}.
\end{gather*}
With the above asymptotic expansions, we then find that indeed,
\begin{gather*}
  \eqref{asymp-0}\Bigm|_{\substack{t\to\frac{1}{t}\\
      y\to-i\frac{y}{t}}}
  \!{}+{}
  i\sqrt{t}\,\EX{-\pi\frac{y^2}{t}}\cdot\eqref{asymp-infty}=
  -\ffrac{i}{2}\,\sqrt{t} - \fhalf\,\EX{-\pi\frac{y^2}{t}},
\end{gather*}
showing that the asymptotic expansions of the integral
in~\eqref{Phi-def} obey the $S$ transformation formula~\eqref{Phi-S}.

\section{Geometry and Further Properties of the Appell Functions}
\label{sec:geometry}
In this section, we consider some elements of the geometric
interpretation of higher-level Appell functions and then formulate
their modular properties in terms of the action of subgroups of
$SL(2,\oZ)$.  As noted in~\cite{[Polisch]},
$(\kappa({}\cdot{},a;q),\vartheta(q,{}\cdot{}))$ is a section of a
rank-$2$ vector bundle over the torus.  For $\K_\pell$ with
$\pell>1$, we unify $\K_\pell$ and the theta functions arising in
the open quasiperiodicity formula for $\K_\pell$ into a
vector~$\Kvect_\pell$ representing a section of a rank-$(\pell+1)$
bundle.  In the space of $(\pell+1)$-vectors, sections of this bundle
are selected by the \textit{invariance} condition with respect to the
action of a lattice in~$\oC^2$ (Lemma~\bref{Lemma:Kvect}).  Moreover,
there is a representation of a ($\pell$-dependent) subgroup of
$SL(2,\oZ)$ on $(\pell+1)$-vectors such that the
section~$\Kvect_\pell$ is also invariant under this subgroup action
(Theorem~\bref{thm:abcd}).  This description of the action of (a
subgroup of) $\SLiiZ$ as an invariance statement is in the spirit of
the well-known result for classic theta functions, which we quote in
Lemma~\bref{lemma:Mumford} below.  The ``essence'' of the modular
group action is then hidden in an automorphy factor involved in
defining this action.  With $\Kvect_\pell$, similarly, the invariance
statement in Theorem~\bref{thm:abcd} below involves a judiciously
chosen automorphy \textit{matrix}.

\subsection{Lattice translations and bundle
  sections}\label{sec:lattice} We begin with constructing a vector
bundle $\bundle{\pell,\tau}$, of the rank determined by the number of
terms in the right-hand side of the open quasiperiodicity formula
for~$\K_\pell$.  This bundle
\begin{gather} \label{eq:the-bundle}
  \begin{array}{c}
    \oC^2\times\oC^{\pell+1}\!\bigl/\mathscr{R}\\[2pt]
    \Bigm\downarrow\\
    \torus{\pell,\tau}
  \end{array}
\end{gather}
is defined as follows.  We take the $4$-dimensional
torus~$\torus{\pell,\tau}=\oC^2/\lattice{\pell,\tau}$, where
$\lattice{\pell,\tau}\subset\oC^2$ is the lattice generated by the
vectors
\begin{gather*}
  \gamma_1=\vect e,\quad
  \gamma_2=\ffrac{1}{\pell}(\vect e-\vect f),\quad
  \gamma_3=\tau\vect e,\quad
  \gamma_4=\ffrac{1}{\pell}\tau(\vect e-\vect f),
\end{gather*}
with~$\vect e$ and $\vect f$ being the standard basis in~$\oC^2$.
With~$\nu$ and~$\mu$ denoting the corresponding coordinates and
$(\nu,\mu,v)\in\oC^2\times\oC^{\pell+1}$, the relations $\mathscr{R}$ are
given by
\begin{multline*}
  \mathscr{R}=
  \bigl\{(\nu+1,\mu,v)\sim(\nu,\mu,v),~
  (\nu+\ffrac{1}{\pell},\mu-\ffrac{1}{\pell},v)\sim(\nu,\mu,A_2^{-1}
  v),
  \\
  (\nu+\tau,\mu,v)\sim(\nu,\mu,A_3(\tau,\nu)^{-1}v),~
  (\nu+\ffrac{\tau}{\pell},\mu-\ffrac{\tau}{\pell},v)
  \sim(\nu,\mu,A_4(\tau,\nu,\mu)^{-1}v) \bigr\}.
\end{multline*}
where the matrices~$A_2$, $A_3(\tau,\nu)$, and~$A_4(\tau,\nu,\mu)$ are
\begin{align*} %%\label{B-matrix}
  A_2={}&
  \begin{pmatrix}
    1& 0 & 0 & 0 &\dots& 0 \\
    0& 1 & 0 & 0 &\dots& 0 \\
    0& 0 & \ex{-2i\pi\frac{1}{\pell}}&0 &\dots&0\\
    0& 0 & 0 &\ex{-2i\pi\frac{2}{\pell}}&\dots&  0 \\
    \hdotsfor{6}\\
    0& 0 & 0 & 0 &\dots&\ex{-2i\pi\frac{\pell-1}{\pell}}
  \end{pmatrix}\!,\\[4pt]
%%  \label{A-matrix}
  A_3(\tau,\nu)={}&
  \ex{i\pi\pell\tau+2i\pi\pell\nu}\one_{(\pell+1)\times(\pell+1)},
  \\[4pt]
  A_4(\tau,\nu,\mu)={}&
  \ex{2i\pi\nu
    \!+\!i\pi\frac{1}{\pell}\tau}
  \begin{pmatrix}
    \ex{2i\pi\mu-i\pi\frac{1}{\pell}\tau}& \vect{v} \\
    \vect 0& B
  \end{pmatrix}\!,
\end{align*}
where the $\pell\times \pell$-matrix
\begin{gather*}
  B=
  \begin{pmatrix}
    0& 0 &0 &\dots& 0& 1\\
    1& 0 &0 &\dots& 0&0\\
    0& 1 &0 &\dots& 0&0\\
    \hdotsfor{6}\\
    0& 0 & 0 &\dots&1&0
  \end{pmatrix}
\end{gather*}
cyclically permutes the standard basis vectors, and
$\vect{v}=(1,\underbrace{0,0,\dots,0}_{\pell-1})$.  The projection
in~\eqref{eq:the-bundle} is given by~$(\nu,\mu,v)\mapsto(\nu,\mu,0)$.

We introduce the $\pell$-dimensional vector
$\boldsymbol{\theta}_\pell(\tau,\nu)=
\bigl(\theta_{2r,\pell}(\frac{\tau}{2},\nu)\bigr)_{0\leq r\leq
  \pell-1}$ (these are the same theta functions that appear in the
open quasiperiodicity formula for~$\K_\pell$ and
in~\eqref{Knew-Stransf-fin}).

\begin{Lemma} \label{Lemma:Kvect}
  The $(\pell+1)$-vector
  \begin{gather*}
    \Kvect_\pell(\tau,\nu,\mu) =
    \begin{pmatrix}
      \K_\pell(\tau,\nu,\mu)\\
    \boldsymbol{\theta}_\pell(\tau,\nu)
    \end{pmatrix}
  \end{gather*}
  is a section of the bundle~$\bundle{\pell,\tau}$.  
\end{Lemma}
\noindent
The Lemma is almost tautological in view of
Eqs.~\eqref{K-identity-triv}, \eqref{K-identity-integer},
\eqref{quasi-periodicity}, and~\eqref{diag-periodicity} and
quasiperiodicity of the theta functions.  We nevertheless note that
the assertion is conveniently formalized using the operators $U_i$,
$i=1,2,3,4$, acting on functions
$f(\tau,{}\cdot{},{}\cdot{}):\oC^2\to\oC^{\pell+1}$ as
\begin{align*}
  U_1 f(\tau,\nu,\mu)={}&f(\tau,\nu+1,\mu),\\
  U_2 f(\tau,\nu,\mu)={}&
  A_2f(\tau,\nu+\ffrac{1}{\pell},\mu-\ffrac{1}{\pell}),\\
  U_3 f(\tau,\nu,\mu)={}&A_3(\tau,\nu)f(\tau,\nu+\tau,\mu),\\
  U_4 f(\tau,\nu,\mu)={}&
  A_4(\tau,\nu,\mu)f(\tau,\nu+\ffrac{\tau}{\pell},
  \mu-\ffrac{\tau}{\pell}).
\end{align*}
Lemma~\bref{Lemma:Kvect} is then a reformulation of the following
easily verified fact.

\begin{Lemma}
  The operators $U_i$, $i=1,2,3,4$, pairwise commute, and hence
  $\gamma_i\mapsto U_i$ is a multiplicative representation of the
  Abelian group~$\lattice{\pell,\tau}$.  Sections
  of~$\bundle{\pell,\tau}$ can therefore be identified with invariants
  of~$\,\lattice{\pell,\tau}$ in this representation.
\end{Lemma}

\subsection{The action of a subgroup of the modular group}
\label{sec:modular-full} We now identify a subgroup in $SL(2,\oZ)$
and construct its action on
functions~$f:\upperH\times\oC^2\to\oC^{\pell+1}$ (where $\upperH$ is
the upper half-plane) such that the vector~$\Kvect_\pell$ defined in
Lemma~\bref{Lemma:Kvect} is \textit{invariant}.  The aim of this
subsection is to prove Theorem~\bref{thm:abcd}.  This action involves
a \textit{matrix} automorphy factor, which can be considered an
Appell-function analogue of the automorphy factor involved in the
classic statement in Lemma~\bref{lemma:Mumford}.

\subsubsection{Automorphy factors} We begin with recalling the
$SL(2,\oZ)$ action on~$\upperH\times\oC^2$, Eq.~\eqref{SL2-action}.
Possible $SL(2,\oZ)$ actions (actually, \textit{anti}representations)
on functions $f:\upperH\times\oC\to\oC$ are given by
\begin{gather}\label{act-scalar}
  (\gamma\acts f)(\tau,\nu) =j(\smatrix{a}{b}{c}{d};\tau,\nu)
  f(\gamma\tau, \gamma\nu),
\end{gather}
where $j$ is an \textit{automorphy factor} satisfying the standard
cocycle condition.  More generally, $SL(2,\oZ)$ actions on
functions~$f:\upperH\times\oC^2\to\oC^{\pell+1}$ are given by
 \begin{gather}\label{act-vector}
   \gamma\cdot f(\tau,\nu,\mu)=\mat{J}_\pell(\gamma;\tau,\nu,\mu)
   f(\gamma\tau,\gamma\nu,\gamma\mu),
\end{gather}
where $\mat{J}_\pell(\gamma;\tau,\nu,\mu)$ is the matrix automorphy
factor, a $(\pell+1)\times(\pell+1)$-matrix satisfying the cocycle
condition
\begin{gather}\label{eq:cocycle}
  \mat{J}_\pell(\alpha\beta;\tau,\nu,\mu)= 
  \mat{J}_\pell(\beta;\tau,\nu,\mu)
  \mat{J}_\pell(\alpha;\beta\tau,\beta\nu,\beta\mu),\qquad
  \mat{J}_\pell(\one;\tau,\nu,\mu)=\one.
\end{gather}

Let $\Gamma_{1,2}$ be the subgroup in $SL(2,\oZ)$ consisting of
matrices $\gamma=\smatrix{a}{b}{c}{d}$ with $a b \in2\oZ$ and $c d
\in2\oZ$.  We recall the following result about the invariance of
theta functions under~$\Gamma_{1,2}$.

\begin{Lemma}[\cite{[Mum]}]\label{lemma:Mumford}
  The theta function $\vartheta(\tau,\nu)$ is invariant under the
  action of $\Gamma_{1,2}$ \textup{(}see~\eqref{act-scalar}\textup{)}
  with the automorphy factor
  \begin{gather*}
    j(\smatrix{a}{b}{c}{d};\tau,\nu)=
    \zeta_{c,d}^{-1}(c\tau + d)^{-\half}
    \EX{-i\pi\frac{c\nu^2}{c\tau + d}}.
  \end{gather*}
\end{Lemma}
\noindent
Here,
\begin{gather*}
  \zeta_{c,d}=
  \begin{cases}
    \ex{i\pi\frac{d-1}{4}}\Jacobi{c}{|d|},& c~\text{even},\\
    \ex{-i\pi\frac{c}{4}}\Jacobi{d}{c},& c~\text{odd}
  \end{cases}
\end{gather*}
and $\Jacobi{c}{d}$ is defined as the quadratic residue for odd
positive prime~$d$ and multiplicatively extended to all~$d$;
see~\cite{[Koblitz]} for the details (we note that $c$ and $d$ are
coprime because of the determinant condition in $SL(2,\oZ)$, and we
assume~$\Jacobi{0}{\pm1}=1$).  As a simple corollary of the Lemma, we
have the formula
\begin{gather}\label{abcd-eta}
  \eta\bigl(\ffrac{a\tau + b}{c\tau + d}\bigr)^{\!3}
  {}=\zeta_{c,d}\,(c\tau+d)^{\frac{3}{2}}\,
  \EX{\frac{i\pi}{4}(a+c)(b+d) + \frac{i\pi}{4} - i\pi\frac{a+c}{2}}
  \eta(\tau)^3,\quad
  \smatrix{a}{b}{c}{d}\in\Gamma_{1,2},
\end{gather}
to be used in what follows.

\subsubsection{Matrix automorphy factors} The classic invariance
statement in Lemma~\bref{lemma:Mumford} extends to the theta-vector
$\boldsymbol{\theta}_\pell(\tau,\nu)$ introduced before
Lemma~\bref{Lemma:Kvect}.  This is shown in the next Lemma, but we
first set the necessary notation.

Let~$\Gamma_{1,2\pell}$ be the subgroup of~$SL(2,\oZ)$ consisting of
matrices $\smatrix{a}{b}{c}{d}$ such that $ab\equiv0\,{\rm
  mod}\,2\pell$ and $cd\equiv0\,{\rm mod}\,2\pell$.  For
$\smatrix{a}{b}{c}{d}\in\Gamma_{1,2\pell}$, let
$\gcd(\pell,a)=\pell_a$ and $\gcd(\pell,c)=\pell_c$ (hence,
$\pell_a\,\pell_c=\pell$).

We also need an $\pell$-dimensional representation $\vect{D}'_\pell$
of~$\Gamma_{1,2\pell}$ defined as
%   The $\pell\times \pell$ matrix
%   $\vect{D}_\pell(\gamma)'$ contains the
%   number~$\frac{1}{\sqrt{\pell_a}}\EX{-2i\pi\frac{b c r s}{\pell}}$
%   in the intersection of $s$-row and $((s a-r c)\,{\rm
%   mod}\,\pell)$-column  for~$0\leq r\leq \pell_a-1$ and~$0$
%   elsewhere.
\begin{gather*}
  \gamma=\smatrix{a}{b}{c}{d}\mapsto
  \vect{D}'_\pell(\gamma)=(d_{s,n})_{\substack{
      0\leq s\leq \pell-1\\
      0\leq n\leq \pell-1}},
\end{gather*}
where
\begin{gather*}
  d_{s,n} =
  \begin{cases}
    \frac{1}{\sqrt{\pell_a}}\,\ex{-2i\pi\frac{bc r s}{\pell}}&
    \parbox{.5\textwidth}{if $n\equiv (sa - r c)\,
      \mathrm{mod}\,\pell$\\
      for some integer $0\leq r\leq \pell_a-1$,}\\[10pt]
    0& \text{otherwise}.
  \end{cases}
\end{gather*}

\begin{Lemma} \label{lemma:theta-ell-matrix}
  The vector $\boldsymbol{\theta}_\pell(\tau,\nu)$ is invariant under
  the action of~$\Gamma_{1,2\pell}$ given by
  \begin{gather*}
    \gamma\acts f(\tau,\nu,\mu)=\vect{J}(\gamma;\tau,\nu,\mu)
    f(\gamma\tau,\gamma\nu,\gamma\mu),
  \end{gather*}
  with the $\pell\times \pell$ automorphy factor
  \begin{gather*}
    \vect{J}(\gamma;\tau,\nu,\mu)=
    k_\pell(\gamma;\tau,\nu)\vect{E}_\pell(\gamma;\tau),
  \end{gather*}
  where
  \begin{gather}\label{eq:automorphy-theta}
    k_\pell(\gamma;\tau,\nu)=\EX{-i\pi\frac{c\pell\nu^2}{c\tau+d}},
    \quad
    \vect{E}_\pell(\gamma;\tau)=\zeta_{\frac{c}{\pell_c},\pell_ad}^{-1}
    (c\tau + d)^{-\half}\,\vect{D}_\pell(\gamma),
  \end{gather}
  and $\vect{D}_\pell(\gamma)=\vect{D}'_\pell(\gamma)^{-1}$.
\end{Lemma}
The proof essentially reduces to the formula
\begin{gather}\label{theta-ell}
  \theta^{(\pell)}_s(\ffrac{a\tau+b}{c\tau+d},\ffrac{\nu}{c\tau+d})=
  \zeta_{\frac{c}{\pell_c},\pell_ad}\,
  \ffrac{\sqrt{c\tau + d}}{\sqrt{\pell_a}}\,
  \EX{i\pi\frac{c\pell\nu^2}{c\tau + d}}
  \sum_{r=0}^{\pell_a-1}\EX{-2i\pi\frac{bcrs}{\pell}}
  \theta^{(\pell)}_{sa-rc}(\tau,\nu),
\end{gather}
which can be verified for $\smatrix{a}{b}{c}{d}\in\Gamma_{1,2\pell}$
by direct calculation.

\subsubsection*{Example}
To illustrate the structure of the automorphy factor in the Lemma, we
consider two examples with the matrix $\gamma$ chosen as
$S=\smatrix{0}{-1}{1}{0}$ and $C=S^2=\smatrix{-1}{0}{0}{-1}$. For $S$,
the matrix elements $D_{nm}$ of $\vect{D}_\pell(S)$ are given by
\begin{gather*}
  D_{nm}=\ffrac{1}{\sqrt{\pell}}\,\EX{2i\pi\frac{n\,m}{\pell}},\qquad
  n,m=0,1,\dots,\pell-1.
\end{gather*}
For $C$, we have%\enlargethispage{12pt}
\begin{gather*}
  \vect{D}_\pell(C)=(\vect{D}_\pell(S))^2=
  \begin{pmatrix}
    1&0&\hdotsfor{2}&0\\
    0&\hdotsfor{3}&1\\
    0&\hdotsfor{2}&1&0\\
    \hdotsfor{5}\\
    0&1&0&\!\!\hdots\!\!&0
  \end{pmatrix}.
\end{gather*}

\subsubsection{Modular behavior of $\Kvect_\pell$} With these
ingredients, we now formulate an analogue of
Lemma~\bref{lemma:Mumford} for the vector $\Kvect_\pell$ in
Lemma~\bref{Lemma:Kvect}.\pagebreak[3] As a final preparation, we
define the automorphy factor before formulating the
result, because the corresponding formulas are somewhat bulky in view
of several cases that must be considered (the reader may first
concentrate on the ``basic'' case $c>0$).\pagebreak[3]

For $\gamma=\smatrix{a}{b}{c}{d}\in\Gamma_{1,2\pell}$,
let~$\mat{J}_\pell(\gamma;\tau,\nu,\mu)$ be the
$(\pell+1)\times(\pell+1)$ matrix defined~as
\begin{equation*}
  \mat{J}_\pell(\gamma;\tau,\nu,\mu)
  =\begin{cases}    
    k_\pell(\gamma;\tau,\nu)\!
    \begin{pmatrix}
      l_\pell(\gamma;\tau,\mu)&
      \vect{F}_\pell(\gamma;\tau,\mu)\vect{E}_\pell(\gamma;\tau)\\
      0&\vect{E}_\pell(\gamma;\tau)\!
    \end{pmatrix}\!,\kern-130pt&\kern100pt c>0,\\
    \begin{pmatrix}
      -1& \vect{v}\\
      \vect0&\vect{D}_\pell(C)
    \end{pmatrix}\!,
    &c=0\text{ and }\gamma=C=\smatrix{-1}{0}{0}{-1},\\
    \one,&c=0\text{ and }\gamma=\smatrix{1}{2\pell b}{0}{1},\\
    \mat{J}_\pell(C;\tau,\nu,\mu),&c=0\text{ and }
    \gamma=\smatrix{-1}{2\pell b}{0}{-1},\\
    \mat{J}_\pell(-\gamma;\tau,\nu,\mu)\mat{J}(C;\tau,\nu,\mu),
    &c<0,
  \end{cases}
\end{equation*}
where $k_\pell(\gamma;\tau,\nu)$ and $\vect{E}_\pell(\gamma;\tau)$ are
given in~\eqref{eq:automorphy-theta},
\begin{gather*}
  l_\pell(\gamma;\tau,\nu)
  =(c\tau+d)^{-1}\EX{i\pi\frac{c\pell\nu^2}{c\tau+d}},
\end{gather*}
and $\vect{F}_\pell(\gamma;\tau,\mu)=
(F^{(\pell)}_1,F^{(\pell)}_2,\dots,F^{(\pell)}_\pell)$ is the vector
with the components
\begin{multline*}
  F^{(\pell)}_r=-\zeta_{\frac{c}{\pell_c},\pell_ad}\,i\sqrt{-i}\,
  \sqrt{\ffrac{\pell_c}{c}}\,
  \EX{i\pi\frac{c\pell\mu^2}{c\tau+d}}\,
  \EX{i\pi\frac{rd(2\pell\mu-r\tau)}{\pell(c\tau+d)}}\\*
  {}\times \sum_{s=0}^{\frac{c}{\pell_c}-1}
  \EX{2i\pi (\pell\mu\!-\!r\tau)\frac{s d}{c\tau+d
    }\!-\!i\pi\tau\frac{s^2 d}{c\tau+d}}\,
  \Phi(\pell\tau\!+\!\pell\ffrac{d}{c},
  \pell\mu\!+\!s \pell\ffrac{d}{c}\!-\!r\tau).
\end{multline*}

\begin{Thm}\label{thm:abcd}
  The section of $\,\bundle{\pell,\tau}$ given by
  $\Kvect_\pell(\tau,\nu,\mu)$ is invariant under the action of
  $\;\Gamma_{1,2\pell}$ given by Eqs.~\eqref{act-vector} with the
  matrix automorphy factor~$\mat{J}_\pell$ defined above.
\end{Thm}

\subsubsection*{Example: the ${S}$ transformation} For the $S$ matrix,
the automorphy factor in the Theorem becomes
\begin{gather*}
  \mat{J}_\pell(S;\tau,\nu,\mu)=
  \EX{-i\pi\pell\frac{\nu^2}{\tau}}
  \begin{pmatrix}
    \tau^{-1}\eh^{i\pi\pell\frac{\mu^2}{\tau}}&
    -(-i\tau)^{-\half}\eh^{i\pi\pell\frac{\mu^2}{\tau}}
    \vect{\Psi}\vect{D}_\pell(S)\\
    \vect0&(-i\tau)^{-\half}\vect{D}_\pell(S)
  \end{pmatrix}\!,
\end{gather*}
where $\vect{D}_\pell(S)$ is given above and
$\vect{\Psi}=(\Psi_a)_{a=0,1,\dots,\pell-1}$ is the row vector with
the components $\Psi_a=\Phi(\pell\tau,\pell\mu-a\tau)$.

\subsubsection*{Example: ${\pell=1}$} The above formulas become
somewhat more transparent in the simplest case $\pell=1$.  For
$\smatrix{a}{b}{c}{d}\in\Gamma_{1,2}$ with $c>0$, we then
have
\begin{multline*}
  \K_1(\smatrix{a}{b}{c}{d} \tau,
  \ffrac{\nu}{c\tau+d},\ffrac{\mu}{c\tau+d})=
  (c\tau\!+\!d)\,\EX{i\pi\frac{c(\nu^2-\mu^2)}{c\tau+d}}
  \K_1(\tau, \nu, \mu)\\*
  {}+ \zeta_{c,d}\,i\sqrt{-i}\,
  \ffrac{c\tau+d}{\sqrt{c}}
  \EX{i\pi\frac{c\nu^2}{c\tau + d}} \sum_{\alpha=0}^{c-1}
  \EX{2i\pi\mu\frac{\alpha d}{c\tau+d} - i\pi\tau\frac{\alpha^2
      d}{c\tau+d}}\,
  \Phi(\tau\!+\!\ffrac{d}{c},
  \mu\!+\!\alpha\ffrac{d}{c})\, \vartheta(\tau,\nu).
\end{multline*}
% The case where $c<0$ is reduced to $c>0$ by multiplying the matrix
% by~$-1$; finally, $c=0$ with $d>0$ implies that $a=d=1$ and $b$ is
% even, and therefore, $\K_1(\tau+b, \nu, \mu)=\K_1(\tau, \nu,
% \mu)$.
For the section $\Kvect_1(\tau,\nu,\mu)$ in Lemma~\bref{Lemma:Kvect},
we therefore have
\begin{gather*}
  \smatrix{a}{b}{c}{d}
  \acts\Kvect_1=\Kvect_1,\qquad
  \smatrix{a}{b}{c}{d}\in\Gamma_{1,2},
\end{gather*}
where the action is defined as above %in~\eqref{act-vector}
with the $2\times2$ automorphy factor (for $c>0$)
\begin{multline*}
  \mat{J}_1(\smatrix{a}{b}{c}{d}\!;\tau,\nu,\mu)= (c\tau + d)^{-\half}
  \EX{-i\pi\frac{c\nu^2}{c\tau + d}}\\*
  {}\times
  \begin{pmatrix}
    (c\tau\!+\!d)^{-\half} \eh^{i\pi\frac{c\mu^2}{c\tau + d}}&
    -\frac{i\sqrt{-i}}{\sqrt{c}} \eh^{i\pi\frac{c\mu^2}{c\tau +
        d}} \sum\limits_{\alpha=0}^{c-1} \eh^{i\pi\alpha d\frac{2\mu
        -\tau\alpha}{c\tau+d}}
    \Phi(\tau\!+\!\frac{d}{c},\mu\!+\!\alpha\frac{d}{c})\\
    0 & \zeta_{c,d}^{-1}
  \end{pmatrix}.
\end{multline*}

\subsubsection{Proof of Theorem~\bref{thm:abcd}}
We first evaluate the integral representation~\eqref{K-integral} with
$(\tau,\nu,\mu)$ transformed by an element
$\smatrix{a}{b}{c}{d}\in\Gamma_{1,2}$; using that $a+c$ and $b+d$ are
then odd and applying~\eqref{abcd-eta}, we obtain
\begin{multline}\label{**}
  \K_\pell(\smatrix{a}{b}{c}{d}\tau,
  \ffrac{\nu}{c\tau\!+\!d},\ffrac{\mu}{c\tau\!+\!d})=
  \int\limits_0^{c\tau \!+\! d}
  d\lambda\,
  \vartheta\!\left(\pell\ffrac{a\tau \!+\! b}{c\tau \!+\! d},
    \ffrac{\pell\nu \!-\! \lambda}{c\tau \!+\! d}\right)
  \,\EX{2i\pi\frac{c}{c\tau\!+\!d}(\nu\!+\!\mu)\lambda}\\*
  {}\times
  \frac{\ex{-i\frac{\pi\tau}{4}}\,
    \vartheta(\tau, \nu\!+\!\mu\!+\!\lambda\!+\!\frac{\tau\!+\!1}{2})
    \eta(\tau)^{\!3}}{
    \vartheta(\tau, \nu \!+\! \mu \!+\! \frac{\tau\!+\!1}{2})
    \vartheta(\tau, \lambda \!+\! \frac{\tau\!+\!1}{2})}.
\end{multline}

We next assume that $ac$ is divisible by~$\pell$.
Equation~\eqref{theta-ell} then gives
\begin{multline*}
  \K_\pell(\smatrix{a}{b}{c}{d}\tau,
  \ffrac{\nu}{c\tau+d},\ffrac{\mu}{c\tau+d})={}\\
  {}=
  \zeta_{\frac{c}{\pell_c},\pell_ad}\,
  \ffrac{\sqrt{c\tau + d}}{\sqrt{\pell_a}}\,
  \EX{i\pi\frac{c\pell\nu^2}{c\tau+d}}
  \sum_{r=0}^{\pell_a-1}
  \int\limits_0^{c\tau + d} d\lambda\,
  \EX{i\pi\frac{c(\lambda^2\!+\!2\pell\lambda\mu)}{\pell(c\tau\!+\!d)}
    \!+\! 2i\pi rc(\mu\!+\!\frac{\lambda}{\pell})
    \!+\! i\pi\frac{r^2c^2}{\pell}\tau}\\*
  {}\times\vartheta(\pell\tau,\pell\nu\!-\!(\lambda\!+\!rc\tau))
  \frac{\ex{-i\frac{\pi\tau}{4}}\,
    \vartheta(\tau, \nu\!+\!\mu\!+\!\lambda\!+\!rc\tau
    \!+\!\frac{\tau\!+\!1}{2})
    \eta(\tau)^{\!3}}{
    \vartheta(\tau, \nu \!+\! \mu \!+\! \frac{\tau\!+\!1}{2})
    \vartheta(\tau, \lambda\!+\!rc\tau \!+\! \frac{\tau\!+\!1}{2})},
\end{multline*}
where we also shifted the theta function arguments using
quasiperiodicity in order to have $\lambda+rc\tau$ in the ratio of the
theta functions.  This allows us to apply Lemma~\bref{lemma:favorite},
with the result
\begin{multline} \label{abcd-work20}
  \K_\pell(\smatrix{a}{b}{c}{d}\tau,
  \ffrac{\nu}{c\tau+d},\ffrac{\mu}{c\tau+d})=\\
  \shoveleft{{}=\zeta_{\frac{c}{\pell_c},\pell_ad}\,
    \ffrac{\sqrt{c\tau + d}}{\sqrt{\pell_a}}\,
    \EX{i\pi\frac{c\pell\nu^2}{c\tau+d}}
    \sum_{r=0}^{\pell_a-1}
    \int\limits_0^{c\tau + d} d\lambda\,
    \EX{i\pi\frac{c(\lambda^2 + 2\pell\lambda\mu)}{\pell(c\tau+d)}
      \!+\!2i\pi rc(\mu+\frac{\lambda}{\pell})
      \!+\!i\pi\frac{r^2c^2}{\pell}\tau}}\\*
  \shoveleft{\qquad{}\times
    \Bigl(
    \vartheta(\pell\tau,\pell\mu\!+\!\lambda\!+\!rc\tau)
    \K_\pell(\tau,\nu,\mu)}\\*
  - \sum_{k=0}^{\pell-1} \EX{2i\pi k(\nu+\mu)}
  \vartheta(\pell\tau,\pell\nu+k\tau)
  \K_1(\pell\tau,
  -\pell\mu\!-\!(\lambda\!+\!rc\tau), \pell\mu\!-\!k\tau)\Bigr).
\end{multline}

This therefore consists of two terms, with the integration variable
$\lambda$ involved in the argument of the theta function in the first
term and in the argument of $\K_1$ in the second term.  In the first
term, we change the variable as
$\lambda\to\lambda-r(c\tau+d)-\pell\mu$ and use~\eqref{theta-cd}.  In
the second term, we change the variable as
$\lambda\to\lambda-r(c\tau+d)$ and then use~\eqref{K-cd} with
$c\to\frac{c}{\pell_c}$, $d\to \pell_ad$, $\tau\to \pell\tau$, and
$\mu\to \pell\mu-k\tau$.  This gives
\begin{multline*}
  \K_\pell(\smatrix{a}{b}{c}{d}\tau,
  \ffrac{\nu}{c\tau+d},\ffrac{\mu}{c\tau+d})=
  (c\tau\!+\!d)
  \EX{i\pi\frac{c\pell(\nu^2-\mu^2)}{c\tau+d}}
  \K_\pell(\tau,\nu,\mu)\\*
  \shoveleft{\quad
    {}+\zeta_{\frac{c}{\pell_c},\pell_ad}\,i\,\sqrt{-i}\,    
    \sqrt{\ffrac{\pell_c}{c}}\,(c\tau\!+\!d)\,
    \EX{i\pi\frac{c\pell\nu^2}{c\tau+d}}
    \sum_{k=0}^{\pell-1}
    \EX{2i\pi k\nu\!+\!i\pi\frac{k^2}{\pell}\tau}
    \vartheta(\pell\tau,\pell\nu\!+\!k\tau)}\\*
  {}\times 
  \EX{i\pi\frac{kd(2\pell\mu-k\tau)}{\pell(c\tau+d)}}
  \sum_{\alpha=0}^{\frac{c}{\pell_c}-1}
  \EX{2i\pi (\pell\mu\!-\!k\tau)\frac{\alpha d}{c\tau+d}
    \!-\!i\pi\tau\frac{\alpha^2 d}{c\tau+d}}\,
  \Phi(\pell\tau\!+\!\pell\ffrac{d}{c},
  \pell\mu\!+\!\alpha \pell\ffrac{d}{c}\!-\!k\tau).
\end{multline*}
This shows the desired behavior of the first element
in~$\Kvect_\pell(\tau,\nu,\mu)$.  The rest of the calculation leading to
the statement of Theorem~\bref{thm:abcd} involves only theta functions
and is therefore standard.

\section{Modular Transformations of $\protect\hSSL21$
  Characters}\label{sec:sl21} As an application of the higher-level
Appell functions, we consider the ``admissible'' representations of
the affine Lie superalgebra $\hSSL21$ at the level
$k=\frac{\pell}{u}-1$ with coprime positive integers $\pell$ and $u$.
For $\pell\geq2$, neither these representations nor their characters
are periodic under the spectral flow~\eqref{sl21-spectral}, the
characters cannot therefore be rationally expressed through theta
functions, and Appell functions enter the game.

The spectral flow is an action of the $\oZ$ lattice, and because it
acts on the admissible representation characters freely, there are
infinitely many representations involved and the theory is certainly
nonrational.  It might then be expected that defining the modular
group action would also require infinitely many characters.  But the
actual situation turns out to be somewhat closer to the case of
rational conformal field theories: if extensions among the
representations are taken into account, the spectral flow and the
modular group action can be defined on a finite number of characters.
For the spectral flow, this is shown by elementary manipulations, but
the calculation of the action of~$S\in\SLiiZ$ is more complicated.
The resulting formula for the $S$-transform of the characters in
Theorem~\bref{thm:S-chi} resembles that for the Appell functions:
although the characters are not closed under the $S$ transformation,
the offending terms are given by~$\Phi$ times theta-functional terms,
whose modular properties are already standard.

We have to consider a number of facts pertaining to the $\hSSL21$
representation theory.  We follow~\cite{[ST]},\footnote{See also
  \cite{[BT96],[BHT97]} for aspects of the $\hSSL21$-representation
  theory at fractional level.} with most of the representation-theory
part collected in Appendix~\ref{app:sl21}.  Calculations with the
Appell functions are given in Secs.~\ref{subsec:q-periodicity}
and~\ref{sec:sl21-S} below.

\subsection{Formulation of the main result}\label{sec:notaion} 
In Theorem~\bref{thm:find-chars}, we find the characters
$\chi_{r,s,\pell,u;\theta}$ of the admissible $\hSSL21$-representation
$\mL_{r,s,\pell,u;\theta}$ in the \chzerotext{} sector.  The
${}_{;\theta}$ notation is for the spectral flow transform,
see~\eqref{sl21-spectral}~\cite{[ST],[BFST]}.  The four different
sectors (\chzerotext, \chhalftext, super-\chzerotext, and
super-\chhalftext) are mapped under the $S$ and $T$ transformations as
indicated in~\eqref{ST-diagram}.  Any of the $S$-arrows
in~\eqref{ST-diagram} allows reconstructing any other, and it is
therefore a matter of taste which of these to evaluate explicitly.
The super-\chzerotext{} sector is chosen in the next theorem.  With
theta-functional terms inevitably occurring in the $S$-transform
of~$\chi_{r,s,\pell,u;\theta}$, such terms can be added to the $\chi$
characters from the start (as we see in
Sec.~\ref{subsec:q-periodicity}, the theta-functional terms in
question are actually the characters defined in~\eqref{Omega-short}).
It then turns out that to avoid redundancy, we can label the
admissible characters by $(s,\theta)$ with $s=1,\dots,u$ and
$\theta=0,\dots,u-1$ (see~\eqref{chi-new}).  We then have the
following result.

\begin{Thm}\label{thm:S-chi}
  At the level $k=\frac{\pell}{u}-1$ with coprime positive integers
  $\pell$ and $u$, the $S$-transform of the super-\chzerotext{}
  admissible $\hSSL21$-characters $\chi^{\schzero}_{(s;\theta)}$ is
  given by
  \begin{multline}\label{sR-transform}
    \chi^{\schzero}_{(s;\theta)}( -\ffrac{1}{\tau},\ffrac{\nu}{\tau},
    \ffrac{\mu}{\tau})
    =
    \EX{i\pi k\frac{\nu^2 - \mu^2}{2\tau}}
    \sum_{s'=1}^{u}\!\sum_{\theta'=0}^{u-1}\!
    S_{(s,\theta),(s',\theta')}^{\pell,u}    
    \chi^{\schzero}_{(s';\theta')}(\tau,\nu,\mu)\\*
%%%    
    -\EX{i\pi k\frac{\nu^2 - \mu^2}{2\tau}}
    \sum_{r'=1}^{\pell-1}\!\sum_{s'=1}^{u}
    R_{s,\theta,r',s'}(\tau,\mu)
    \Omega_{r',s'}(\tau,\nu,\mu\!+\!1),
  \end{multline}
  where the characters $\Omega_{r',s'}$ are defined
  in~\eqref{Omega-short},
  \begin{align*}
    S_{(s,\theta),(s',\theta')}^{\pell,u}
    &=\ffrac{1}{u}\,
    \EX{2i\pi\frac{\pell}{u}(s\!+\!s'\!+\!\theta\!+\!\theta'
      \!+\!s' \theta\!+\!s \theta'\!+\!2 \theta \theta')}
    \EX{i\pi\frac{[\pell(s'+1+2\theta')]_u
        - [\pell(s + 2 \theta + 1)]_u}{u}},
  \end{align*}
  and
  \begin{multline*}
    R_{s,\theta,r',s'}(\tau,\mu)
    =\ffrac{(-1)^{r's\!}}{u}\,
    \EX{
      i\pi\frac{\pell}{2u\tau}
      \left(\mu\!+\!\tau(s'\!+\!1\!-\!\frac{u}{\pell}r')\right)^2
      \!-\!i\pi\frac{[\pell(s + 2 \theta + 1)]_u}{u}
      \!+\!i\pi\frac{\pell}{u}(2s'\!+\!s\!-\!ss')
    }\\*
    \shoveleft{{}\times
      \sum_{b=0}^{u-1}
      \Bigl(
      \EX{i\pi(s\!+\!2\theta\!+\!1)
        \frac{2\pell b + [ur'+\pell(s'-1)]_{2\pell}}{u}}
      \Phi(2u\pell\tau,
      -\pell\mu\!-\!2b\pell\tau
      \!-\![ur'\!-\!\pell(s'\!-\!1)]_{2\pell}\tau)}\\*
    {}-\EX{i\pi(s\!+\!2\theta\!+\!1)
      \frac{2\pell (b+1)-[ur'+\pell(s'-1)]_{2\pell}}{u}}
    \Phi(2u\pell\tau,
    -\pell\mu\!-\!2(b\!+\!1)\pell\tau
    \!+\![ur'\!-\!\pell(s'\!-\!1)]_{2\pell}\tau)
    \Bigr).
  \end{multline*}
  \end{Thm}
We recall that $[x]_u$ denotes $x\,\mathrm{mod}\,u$.  The
theorem is proved by a calculation based on the properties of the
Appell functions established above.  The several-step derivation is
given in Sec.~\ref{sec:sl21-S}.

The $S$-transformation formula in the theorem has a triangular
structure similar to that for the Appell functions: the
$\chi_{(s;\theta)}$ characters are transformed through themselves
\textit{and} the additional characters $\Omega_{r',s'}$, while the
latter, being expressed through theta functions, transform through
themselves.
%% \footnote{A certain ``approximation'' to the exact
%%   formula~\eqref{sR-transform} is to consider its ``quotient'' with
%%   respect to theta-functional pieces, i.e., view it \textit{modulo}
%%   the terms involving~$\Omega$.  This may, e.g., lead to a rough
%%   fusion algebra, with the basis determined by only the
%%   $\chi_{(s;\theta)}$ characters and the extensions between the
%%   corresponding representations ignored.}
The $\Omega$ characters in the right-hand side of~\eqref{sR-transform}
are multiplied with the $\Phi$ functions, defined in~\eqref{Phi-def}
and studied in Sec.~\ref{sec:Phi}.  We note that the arguments of
$\Phi$ above depend on $s'$ ``weakly,'' in fact only on
$s'\,\mathrm{mod}\,2$.

\subsection{Relation to the Appell functions and open
  quasiperiodicity}\label{subsec:q-periodicity} We first express the
admissible representation characters~\eqref{eq:gen-char}
and~\eqref{eq:gen-char-bar} through the higher-level Appell functions:
the ``nontrivial'' part of the characters, Eq.~\eqref{def-psi}, is
expressed through the Appell functions with the \textit{even} level
$2\pell$ as
\begin{multline}\label{psi-through-K-exp}
  \psi_{r,s,\pell,u}(\tau,\nu,\mu)
  =\K_{2\pell}(u\tau,-\ffrac{\nu}{2}
  \!-\!\ffrac{u(r\!-\!1)}{2\pell}\tau
  \! +\! \ffrac{s\!-\!1}{2}\tau,
  \fhalf\!-\!\ffrac{\mu}{2}\!+\!\tau\ffrac{u(r\!-\!1)}{2\pell}
  - \tau\ffrac{s\!+\!1}{2})
  \\*[4pt]
  {}-\EX{2i\pi(s\!-\!1)(r\!-\!1)\tau\!-\!2i\pi(r\!-\!1)\nu}
  \kern40pt\\*
  {}\times
  \K_{2\pell}(u\tau, \ffrac{\nu}{2}\!-\!\ffrac{u(r\!-\!1)}{2\pell}\tau
  \!-\!\ffrac{s\!-\!1}{2}\tau, 
  \fhalf\!-\!\ffrac{\mu}{2}\!+\!\tau\ffrac{u(r\!-\!1)}{2\pell}
  - \tau\ffrac{s\!+\!1}{2}),
\end{multline}
and hence the characters are given by
\begin{multline}\label{chi-through-K-exp}
  \chi_{r,s,\pell,u;\theta}(\tau,\nu,\mu)
  =\THETA(\tau,\nu,\mu)\\*[4pt]
  {}\times
  \ex{2i\pi(\theta\!+\!1)(r\!-\!1\!-\!\frac{\pell}{u}(s\!+\!\theta))\tau
    \!+\!i\pi(r\!-\!1\!-\!\frac{\pell}{u}(s\!-\!1))\nu
    \!+\!i\pi(r\!-\!1\!-\!\frac{\pell}{u}(s\!+\!1\!+\!2\theta))\mu}
  \\*
  {}\times
  \Bigl(\!
  \K_{2\pell}(u\tau,-\ffrac{\nu}{2}\!-\!\ffrac{u(r\!-\!1)}{2\pell}\tau
  \! +\! \ffrac{s-1}{2}\tau,
  \ffrac{1\!-\!\mu}{2}\!+\!\tau\ffrac{u(r\!-\!1)}{2\pell}
  - \tau\ffrac{s\!+\!1\!+\!2\theta}{2})\\*[4pt]
  {}-\EX{2i\pi(s\!-\!1)(r\!-\!1)\tau\!-\!2i\pi(r\!-\!1)\nu}
  \kern115pt\\*
  \times
  \K_{2\pell}(u\tau, \ffrac{\nu}{2}\!-\!\ffrac{u(r\!-\!1)}{2\pell}\tau
  \!-\!\ffrac{s\!-\!1}{2}\tau, 
  \ffrac{1\!-\!\mu}{2}\!+\!\tau\ffrac{u(r\!-\!1)}{2\pell}
  \!-\!\tau\ffrac{s\!+\!1\!+\!2\theta}{2})\!\Bigr).
\end{multline}

We now discuss the range of the labels in these characters.  First,
$1\leq s\leq u$.  Second, the twist $\theta$ takes all integer values
in principle, but modulo addition of theta-functional terms to the
characters, the $\theta$ parameter can be restricted to $u$
consecutive values, because quasiperiodicity of $\K_\pell$ implies an
open quasiperiodicity property relating the admissible $\hSSL21$
characters with their spectral-flow transform by $\theta=u$
(see~\eqref{sl21-spectral}),
\begin{multline}\label{eq:open-periodicity}
  \chi_{r,s,\pell,u;\theta+u}
  = \chi_{r,s,\pell,u;\theta}
  -\sum_{a=1}^{r-1}(-1)^{a+r+1}
  \Omega_{a,s;\theta}\\*
  {}+\sum_{a=1}^{\pell-1}(-1)^{a+r+1}
  \bar\Omega_{a,s;\theta+s+1}
  -\sum_{a=r}^{\pell-1}(-1)^{a+r+1}
  \bar\Omega_{a,s;\theta+u},
\end{multline}
where $\Omega_{r,s,j}$ are the characters~\eqref{massive-char},
expressible in terms of theta functions.  Thus, although the
admissible $\hSSL21$-characters are not invariant under the spectral
flow~\eqref{sl21-spectral} with any $\theta$ (i.e., are mapped into
nonisomorphic representations), their characters are invariant modulo
theta-functional terms.\footnote{This is somewhat similar to the case
  with the admissible $\hSL2$ representations, which are not invariant
  under any spectral flow transformations, whereas their characters,
  given by quasiperiodic functions, are invariant under a certain
  sublattice of spectral flow translations.  (That is not a
  contradiction because the admissible representation characters are
  \textit{meromorphic} functions of the variable that is translated
  under the spectral flow, cf.~\cite{[FSST]}, and we are actually
  speaking of \textit{analytic continuation} of characters.)}
%% We note that $\Omega_{r,s,h}(q,x,y)$ defined in~\eqref{massive-char}
%% are characters of irreducible modules if
%% conditions~\eqref{cond-irreducible} are satisfied, but the characters
%% arising in~\eqref{eq:open-periodicity} are those of \textit{reducible}
%% modules, moreover, those where $\mL_{r',s,\pell,u}$ occur as
%% submodules and quotients.

\medskip

\noindent\textbf{Remark: $\boldsymbol{\pell=1}$.} \ For $\pell=1$,
Eq.~\eqref{even-identity2} allows us to express the characters in
terms of theta functions, which immediately shows that the characters
are periodic under the spectral flow with the period~$u$ and readily
leads to their modular transformation properties (the modular
transformations for $\pell=1$ were derived in~\cite{[John00]}).  This
originates in the fact that $\cU_\theta$ with $\theta=u$ acts as an
isomorphism on the representations $\mL_{1,s,1,u}$ (the
representations are invariant under the $\beta$
automorphism~\eqref{beta}, and accordingly, only one of the two values
of $r$ remains for $\pell=1$; for the characters, this can be easily
verified using formulas in Sec.~\ref{sec:Appell-elementary}).  There
remain $u^2$ representations $\mL_{1,s,1,u;\theta}$ with $1\leq s\leq
u$ and $0\leq \theta\leq u-1$~\cite{[BFST]}.

We assume $\pell\geq2$ in what follows and often abbreviate the
notation $\chi_{r,s,\pell,u}$, $\psi_{r,s,\pell,u}$, etc.\ to
$\chi_{r,s}$, $\psi_{r,s}$, etc.

Next, it follows from formulas in Sec.~\ref{sec:Appell-elementary}
that the dependence of the characters on~$r$ also amounts to additive
theta-functional terms; for $r\geq2$,
\ $\psi_{r,s%%,\pell,u 
}$ is related to $\psi_{1,s%%,\pell,u
}$ as
\begin{gather}\label{reduce-to-1}
  \psi_{r,s%%,\pell,u
  }(q,x,y)
  = (-x^{-\half}y^{-\half}q^{-1})^{r-1}
  \Bigl(\psi_{1,s%%,\pell,u
  }(q,x,y)
  + \Theta_{r,s%%,\pell,u
  }(q,x,y)\Bigr),
\end{gather}
where
\begin{equation*}
  \Theta_{r,s%%,\pell,u
  }(q,x,y)=
  \sum_{r'=1}^{r-1} (-1)^{r'}
  y^{\frac{r'}{2}}
  q^{-\frac{u{r'}^2}{4\pell} + \frac{r'(s+1)}{2}}
  \left(
    \theta_{r',\pell}(q^u,xq^{-(s-1)})
    -\theta_{-r',\pell}(q^u,xq^{-(s-1)})
  \right)  .
\end{equation*}
Consequently, the corresponding formula relating $\chi_{r,s}$ to
$\chi_{r',s}$ involves the $\Omega_{r,s,h}$
characters~\eqref{massive-char}:
\begin{multline}\label{r-relation}
  \chi_{r,s%%,\pell,u
    ;\theta}(q,x,y)=
  (-1)^{r-r'}\chi_{r',s%%,\pell,u
    ;\theta}(q,x,y)\\*
  {}+{}\begin{cases}
    \displaystyle
    \sum_{a=0}^{r-r'-1}(-1)^a
    \Omega_{r-1-a,s,\frac{r-1-a}{2}-\frac{\pell}{u}\frac{s+1+2\theta}{2}}
    (q,x,y),&r'<r,\\
    \displaystyle
    \sum_{a=0}^{r'-r-1}(-1)^a
    \Omega_{r+a,s,\frac{r+a}{2}-\frac{\pell}{u}\frac{s+1+2\theta}{2}}
    (q,x,y),
    &r'>r.
  \end{cases}
\end{multline}
This is another instance of a triangular structure, with the
characters being invariant only up to theta-functional terms, which
are already invariant.

To avoid redundancy, therefore, any fixed value of $r$ can be used,
for example only the characters $\chi_{1,s,\pell,u;\theta}$ can be
considered.  More generally, we can choose a unique $r$ for each
$(s,\theta)$.  As we see in what follows, a useful such choice is to
set $r=\floor{\frac{\pell}{u}(s+2\theta+1)}$+1.  We use the special
notation for these characters,
\begin{gather}\label{chi-new}
  \chi^{\vphantom{y}}_{(s;\theta)}=
  \chi^{\vphantom{y}}_{
    \floor{\frac{\pell}{u}(s+2\theta+1)}+1,s,\pell,u;\theta},
  \quad 1\leq s\leq u,~0\leq\theta\leq u-1.
\end{gather}

\subsection{Evaluating the $S$-transform of
  $\chi_{r,s,\pell,u;\theta}$}\label{sec:sl21-S} The $S$-transform of
$\chi_{r,s;\theta}(\tau,\nu,\mu)
\equiv\chi_{r,s,\pell,u;\theta}(\tau,\nu,\mu)$ is found in several
steps.

With the aim to use Eq.~\eqref{reduce-to-1}, which allows doing the
most difficult part of the calculation for $r=1$ only, we first
rewrite the character as
\begin{multline*}%%%\label{eq:gen-char-exp-mod}
  \chi_{r,s%%,\pell,u
    ;\theta}(\tau,\nu,\mu)
  = \ex{-2i\pi\frac{\pell}{u}(\theta\!+\!1)(s\!+\!\theta)\tau
    \!-\!i\pi\frac{\pell}{u}(s\!-\!1)\nu
    \!-\!i\pi\frac{\pell}{u}(s\!+\!1\!+\!2\theta)\mu}
  \\*
  {}\times(-1)^{r-1}
  \bigl(\psi_{1,s%%,\pell,u
  }
  + \Theta_{r,s%%,\pell,u
  }\bigr)(\tau,\nu,\mu\!+\!2\theta\tau)\,
  \THETA(\tau,\nu,\mu).
\end{multline*}
In evaluating $\chi_{r,s%%,\pell,u
  ;\theta}(-\ffrac{1}{\tau},\ffrac{\nu}{\tau}, \ffrac{\mu}{\tau})$, we
then use Eqs.~\eqref{theta10}--\eqref{modular-eta} to find the
$S$-transform of~$\THETA$, which gives
\begin{multline*}
  \chi_{r,s%%,\pell,u
    ;\theta}(-\ffrac{1}{\tau},\ffrac{\nu}{\tau},
  \ffrac{\mu}{\tau})
  =(-1)^{r-1}\ffrac{1}{\tau}\,
  \ex{2i\pi(\theta\!+\!1)\frac{\pell}{u}(s\!+\!\theta)\frac{1}{\tau}}\,
  \ex{-i\pi(s\!-\!1)\frac{\pell}{u}\frac{\nu}{\tau}}\,
  \ex{-i\pi\frac{\pell}{u}(s\!+\!1\!+\!2\theta)\frac{\mu}{\tau}}
  \\*
  {}\times  
  \bigl(\psi_{1,s%%,\pell,u
  }
  + \Theta_{r,s%%,\pell,u
  }\bigr)(-\ffrac{1}{\tau},\ffrac{\nu}{\tau},
  \ffrac{\mu\!-\!2\theta}{\tau})\,
  \EX{i\pi\frac{(\mu-\tau)^2-\nu^2}{2\tau}}\!
  \THETA(\tau,\nu,\mu\!-\!\tau\!+\!1),
\end{multline*}
where it remains to find the $S$-transform of $\psi_{1,s} +
\Theta_{r,s}$.  For $\Theta_{r,s}$, the calculation is again standard,
based on
\begin{gather*}
  \theta_{a,\pell}(-\ffrac{u}{\tau}, \ffrac{\nu}{\tau})
  -\theta_{-a,\pell}(-\ffrac{u}{\tau}, \ffrac{\nu}{\tau})
  =\sqrt{\ffrac{-i\tau}{2\pell u}}\,
  \EX{i\pi\frac{\pell\nu^2}{2u\tau}}
  \sum_{r'=1}^{2\pell-1}
  \EX{-i\pi\frac{a r'}{\pell}}
  \bigl(
  \theta_{r',\pell}(\ffrac{\tau}{u}, \ffrac{\nu}{u})
  -\theta_{-r',\pell}(\ffrac{\tau}{u}, \ffrac{\nu}{u})
  \bigr).
\end{gather*}
For $\psi_{1,s}$, we express through the Appell functions as
in~\eqref{psi-through-K-exp} and use the $S$-transform
formula~\eqref{Knew-Stransf-fin}, which we rewrite for the
level~$2\pell$,
\begin{gather*}%%%\label{Knew-Stransf-fin-2ell}
  \K_{2\pell}(-\ffrac{1}{\tau},\ffrac{\nu}{\tau},\ffrac{\mu}{\tau})
  {}=
  \tau \EX{2i\pi\pell\frac{\nu^2 - \mu^2}{\tau}}
  \K_{2\pell}(\tau,\nu,\mu)
  + \tau \EX{2i\pi\pell\frac{\nu^2}{\tau}}
  \sum_{r'=0}^{2\pell-1}\!
  \Phi(2\pell\tau,2\pell\mu\!-\!r'\tau)
  \theta_{r',\pell}(\tau, 2\nu).
\end{gather*}
This gives
\begin{multline*}
  \psi_{1,s%%,\pell,u
  }(
  -\ffrac{1}{\tau},\ffrac{\nu}{\tau},\ffrac{\mu}{\tau})
  = \ffrac{\tau}{u}\,
  \smash{\EX{i\pi\pell\frac{(\nu+s-1)^2 - (\mu-\tau-s-1)^2}{2u\tau}}}\\*
  \shoveright{{}\times\bigl(
    \K_{2\pell}(\ffrac{\tau}{u},-\ffrac{\nu\!+\!s\!-\!1}{2u},
    -\ffrac{\mu\!-\!\tau\!-\!s\!-\!1}{2u})
    -\K_{2\pell}(\ffrac{\tau}{u},\ffrac{\nu\!+\!s\!-\!1}{2u},
    -\ffrac{\mu\!-\!\tau\!-\!s\!-\!1}{2u})
    \bigr)}
  \\*
  {}-\ffrac{\tau}{u}\,
  \EX{i\pi\pell\frac{(\nu+s-1)^2}{2u\tau}}\sum_{r'=0}^{2\pell-1}
  \Phi(\ffrac{2\pell\tau}{u},-\pell\,\ffrac{\mu\!-\!\tau\!-\!s\!-\!1}{u} 
  - r'\ffrac{\tau}{u})\kern20pt\\*
  \qquad{}\times\bigl(\theta_{r',\pell}(\ffrac{\tau}{u},
  \ffrac{\nu\!+\!s\!-\!1}{u})
  - \theta_{-r',\pell}(\ffrac{\tau}{u},\ffrac{\nu\!+\!s\!-\!1}{u})
  \bigr).
\end{multline*}
Putting the $\psi$ and $\Theta$ parts together, we use the second
identity in~\eqref{phi-1} to obtain
\begin{multline}\label{relevant}
  (\psi_{1,s} + \Theta_{r,s})
  (-\ffrac{1}{\tau},\ffrac{\nu}{\tau},\ffrac{\mu}{\tau})
  =\ffrac{\tau}{u}\,
  \EX{i\pi\pell\frac{(\nu+s-1)^2 - (\mu-\tau-s-1)^2}{2u\tau}}\\*
  {}\times\Bigr(\K_{2\pell}(\ffrac{\tau}{u},
  -\ffrac{\nu\!+\!s\!-\!1}{2u},
  -\ffrac{\mu\!-\!\tau\!-\!s\!-\!1}{2u})
  -\K_{2\pell}(\ffrac{\tau}{u},\ffrac{\nu\!+\!s\!-\!1}{2u},
  -\ffrac{\mu\!-\!\tau\!-\!s\!-\!1}{2u})\\*
  -\boldsymbol{\mathscr{F}}_{r,s}(\tau,\nu,\mu\!-\!\tau)\!\Bigr),
\end{multline}
where
\begin{multline}\label{F-define}
  \boldsymbol{\mathscr{F}}_{r,s}(\tau,\nu,\mu)=
  \smash[t]{\EX{i\pi\frac{\pell}{2u\tau}(\mu\!-\!s\!-\!1
      \!+\!\frac{u}{\pell}(r\!-\!1))^2}
    \sum_{r'=1}^{2\pell-1}
    \EX{i\pi r'\frac{r-1}{\pell}}}\\*
  {}\times\Phi(\ffrac{2\pell\tau}{u},
  -\ffrac{\pell}{u}(\mu\!-\!s\!-\!1)
  \!-\!\ffrac{r'\tau}{u} \!-\! r\!+\!1)
  \bigl(
  \theta_{r',\pell}(\ffrac{\tau}{u}, \ffrac{\nu\!+\!s\!-\!1}{u})
  -\theta_{-r',\pell}(\ffrac{\tau}{u}, \ffrac{\nu\!+\!s\!-\!1}{u})
  \bigr).
\end{multline}

Therefore, $\chi_{r,s;\theta}(-\ffrac{1}{\tau},\ffrac{\nu}{\tau},
\ffrac{\mu}{\tau})$ is evaluated as
\begin{multline}\label{intermediate-1}
  \chi_{r,s;\theta}(-\ffrac{1}{\tau},\ffrac{\nu}{\tau},
  \ffrac{\mu}{\tau})=
  \ffrac{(-1)^{r-1}}{u}\,
  \EX{i\pi k\frac{\nu^2-(\mu-\tau)^2}{2\tau}
    -i\pi\frac{\pell}{u}(s\!+\!1\!+\!2\theta)}
  \THETA(\tau,\nu,\mu\!-\!\tau\!+\!1)
  \\*
  \shoveright{{}\times
    \Bigr(\!\K_{2\pell}(\ffrac{\tau}{u},-\ffrac{\nu\!+\!s\!-\!1}{2u},
    -\ffrac{\mu\!-\!\tau\!-\!s\!-\!2\theta\!-\!1}{2u})
    {}-\K_{2\pell}(\ffrac{\tau}{u},\ffrac{\nu\!+\!s\!-\!1}{2u},
    -\ffrac{\mu\!-\!\tau\!-\!s\!-\!2\theta\!-\!1}{2u})\!\Bigr)
  }\\*[4pt]
%%%
  {}+{}\ffrac{(-1)^{r}}{u}\,\EX{i\pi k\frac{\nu^2-(\mu-\tau)^2}{2\tau}
      -i\pi\frac{\pell}{u}(s\!+\!1\!+\!2\theta)}
  \boldsymbol{\mathscr{F}}_{r,s}(\tau,\nu,\mu\!-\!\tau\!-\!2\theta)
  \THETA(\tau,\nu,\mu\!-\!\tau\!+\!1),
\end{multline}
where $k=\frac{\pell}{u}-1$ is the $\hSSL21$ level.

The next step is to show that the $\K_{2\pell}$-terms
in~\eqref{intermediate-1} are expressible through the $\chi$
characters and the theta-functional terms are expressible through the
$\Omega$ characters.  We first show that the term involving
$\boldsymbol{\mathscr{F}}_{r,s}$, which has arisen in
form~\eqref{F-define}, can be expressed through the characters
$\Omega_{r',s'}$ introduced in~\eqref{Omega-short}.  Elementary
manipulations show that
\begin{multline*}
  \boldsymbol{\mathscr{F}}_{r,s}(\tau,\nu,\mu)=
  \ex{i\pi\frac{\pell}{2u\tau}(\mu\!-\!s\!-\!1
    \!+\!\frac{u}{\pell}(r\!-\!1))^2}\\*
  \qquad\qquad{}\times
  \ffrac{1}{2\pell}\sum_{a=r-2\pell}^{r-1}
  \sum_{r'=1}^{2\pell-1}
  \EX{i\pi r'\frac{r-1-a}{\pell}}\,
  \Phi(\ffrac{2\pell\tau}{u},
  -\ffrac{\pell}{u}(\mu\!-\!s\!-\!1)
  \!-\!\ffrac{r'\tau}{u} \!-\! r\!+\!1)\\*
  {}\times\bigl(
  \vartheta(\ffrac{\tau}{2\pell u},
  \ffrac{\nu\!+\!s\!-\!1}{2u}\!+\!\ffrac{a}{2\pell})
  - \vartheta(\ffrac{\tau}{2\pell u},
  \ffrac{\nu\!+\!s\!-\!1}{2u}\!-\!\ffrac{a}{2\pell})
  \bigr).
\end{multline*}
Here, the range of the $a$ summation ($2\pell$ consecutive values) can
be shifted arbitrarily, and it was chosen such that~\eqref{phi-2p}
becomes applicable with no remanining~$[\cdot]_{2\pell}$.
Applying~\eqref{phi-2p} to $\Phi$ and rearranging the $\vartheta$-part
in accordance with~\eqref{theta-rs}, we then have
\begin{multline*}
  \boldsymbol{\mathscr{F}}_{r,s}(\tau,\nu,\mu)=
  \ffrac{1}{2\pell}
  \!\sum_{a=r-2\pell}^{r-1}\!
  \EX{i\pi\frac{\pell}{2u\tau}(\mu\!-\!s\!-\!1\!+\!\frac{u}{\pell}a)^2}
  \!\Phi(\ffrac{\tau}{2\pell u},
  -\ffrac{\mu\!-\!s\!-\!1}{2u}\!-\!\ffrac{a}{2\pell})\\*
  \times{}\!
  \sum_{r''=1}^{2\pell-1}\!\sum_{s''=1}^{u}
  \EX{i\pi\frac{u}{\pell}a r''
    \!-\!i\pi\frac{\pell}{u}(s\!-\!1)(s''\!-\!1)
    \!+\!i\pi a(s''\!-\!1)
    \!+\!i\pi(s\!-\!1)r''}
  \varpi_{r'',s''}(\tau,\nu),
\end{multline*}
where 
\begin{gather*}
  \varpi_{r,s}(q,x)= x^{-\frac{\pell}{u}\frac{s-1}{2}}
  q^{\frac{\pell}{4u}(s-1)^2}\!\bigl( \theta_{r,\pell}(q^u,xq^{-(s-1)})
  -\theta_{-r,\pell}(q^u,xq^{-(s-1)})\bigr)
\end{gather*}
accumulates the $y$- and $h$- independent factors in the
character~\eqref{massive-char}.  {}From identity~\eqref{dual-scaling},
we now have
\begin{multline*}
  \boldsymbol{\mathscr{F}}_{r,s}(\tau,\nu,\mu)
  = \EX{i\pi\frac{u\pell}{2\tau}
    (\frac{\mu-s-1}{u}\!+\!\frac{r-1}{\pell})^2}
  \sum_{r'=1}^{2\pell-1}\!\sum_{s'=1}^{u}
  \EX{i\pi\frac{u}{\pell}
    \left(r\!-\!1\!-\!\frac{\pell}{u}(s-1)\right)
    \left(r'\!+\!\frac{\pell}{u}(s'\!-\!1)\right)}\\*
  \times
    \Phi(\ffrac{2\pell\tau}{u},-\ffrac{\pell}{u}(\mu\!-\!s\!-\!1)
    \!-\!r\!+\!1
    \!-\!\ffrac{[ur'\!-\!\pell(s'\!-\!1)]_{2\pell}}{u}\tau)
    \varpi_{r',s'}(\tau,\nu).
\end{multline*}
Next, identity~\eqref{phi-u} (with $p=u$) allows us to express the
$\Phi$ functions involved here through
$\Phi(2u\pell\tau,\,{\dots}{}-2\pell b\tau)$.  This produces the
integer $\pell(s+1)-u(r-1)$ in the argument of each $\Phi$, to which
we can further apply~\eqref{phi-1}, with the result
\begin{multline*}
  \boldsymbol{\mathscr{F}}_{r,s}(\tau,\nu,\mu)
  =\EX{i\pi\frac{\pell}{u}\frac{\mu^2}{2\tau}}
  \sum_{r'=1}^{2\pell-1}\!\sum_{s'=1}^{u}\!\sum_{b=0}^{u-1}
  \EX{-i\pi(s\!-\!1)(r'\!+\!\frac{\pell}{u}(s'\!-\!1))
  \!+\!2i\pi b\frac{\pell}{u}(s\!+\!1)}\\*  
  \shoveright{{}\times
    \EX{i\pi(s\!+\!1)\frac{[ur'+\pell(s'-1)]_{2\pell}}{u}}
    \Phi(2u\pell\tau,
    -\pell\mu\!-\![ur'\!-\!\pell(s'\!-\!1)]_{2\pell}\tau
    \!-\!2b\pell\tau)
    \varpi_{r',s'}(\tau,\nu)}\\*
  \shoveleft{\qquad{}+\EX{i\pi\frac{u\pell}{2\tau}\left(
        \frac{\mu-s-1}{u}\!+\!\frac{r-1}{\pell}\right)^2}
    \sum_{r'=1}^{2\pell-1}\!\sum_{s'=1}^{u}
    \EX{i\pi\frac{u}{\pell}
      \left(r\!-\!1\!-\!\frac{\pell}{u}(s\!-\!1)\right)
      \left(r'\!+\!\frac{\pell}{u}(s'\!-\!1)\right)}}\\*
     \times
    H_{\pell(s+1)-u(r-1),ur'-\pell(s'-1)}(\tau,\mu)
    \varpi_{r',s'}(\tau,\nu),
\end{multline*}
where 
\begin{gather*}
  H_{m,m'}(\tau,\mu)=
  \begin{cases}
    \ffrac{iu}{\sqrt{-2iu\pell\tau}}
    \displaystyle\sum_{\substack{n\geq1\\
        nu\leq m}}
    \EX{i\pi\frac{n^2u-2nm}{2\pell\tau}
      \!+\!i\pi n(\frac{\mu}{\tau}\!+\!\frac{m'}{\pell})}, & m\geq1,\\
    -\ffrac{iu}{\sqrt{-2iu\pell\tau}}
    \displaystyle\sum_{\substack{n\leq0\\
        nu\geq m+1}}
    \EX{i\pi\frac{n^2u-2nm}{2\pell\tau}
      \!+\!i\pi n(\frac{\mu}{\tau}\!+\!\frac{m'}{\pell})}, & m\leq-1.
  \end{cases}
\end{gather*}

We observe that the first (triple-sum) term in
$\boldsymbol{\mathscr{F}}_{r,s}$ is actually independent of~$r$.  The
second term does depend on~$r$, but $H_{m,m'}$ vanishes whenever
$1\leq m\leq u-1$, and because we consider such $m$ in what follows
(more precisely, $r$ such that the corresponding $m$ is within this
interval), we drop the second (double-sum) term
in~$\boldsymbol{\mathscr{F}}_{r,s}$.  More precisely,
Eq.~\eqref{intermediate-1} actually
involves~$\boldsymbol{\mathscr{F}}_{r,s}$ with the argument shifted by
$-2\theta$ with integer~$\theta$.  The integer in the argument of
$\Phi$ mentioned above, and hence the $m$ label in $H_{m,m'}$, is then
$\pell(s+2\theta-1)-u(r-1)$.  With $r$ chosen such that $1\leq m\leq
u-1$, it then follows that
$\boldsymbol{\mathscr{F}}_{s}(\tau,\nu,\mu-2\theta)
\equiv\boldsymbol{\mathscr{F}}_{r,s}(\tau,\nu,\mu-2\theta)$ is given
by
\begin{multline*}
  \boldsymbol{\mathscr{F}}_{s}(\tau,\nu,\mu-2\theta)
  =\EX{i\pi\frac{\pell}{u}\frac{\mu^2}{2\tau}}
  \sum_{r'=1}^{2\pell-1}\!\sum_{s'=1}^{u}\!\sum_{b=0}^{u-1}
  \EX{-i\pi(s\!-\!1)(r'\!+\!\frac{\pell}{u}(s'\!-\!1))
  \!+\!2i\pi b\frac{\pell}{u}(s\!+\!2\theta\!+\!1)}\\*  
  \times
    \EX{i\pi(s\!+\!2\theta\!+\!1)\frac{[ur'+\pell(s'-1)]_{2\pell}}{u}}
    \Phi(2u\pell\tau,
    -\pell\mu\!-\![ur'\!-\!\pell(s'\!-\!1)]_{2\pell}\tau
    \!-\!2b\pell\tau)
  \varpi_{r',s'}(\tau,\nu).
\end{multline*}
Using \eqref{omega+ell} and the fact that $[-a]_n=n-[a]_n$ for $a>0$,
we obtain
\begin{multline}\label{F-final}
  \boldsymbol{\mathscr{F}}_{s}(\tau,\nu,\mu\!-\!2\theta)
  =\EX{i\pi\frac{\pell}{u}\frac{\mu^2}{2\tau}}
  \sum_{r'=1}^{\pell-1}\!\sum_{s'=1}^{u}\!\sum_{b=0}^{u-1}
  \EX{-i\pi(s\!-\!1)(r'\!+\!\frac{\pell}{u}(s'\!-\!1))}\\*
  \shoveleft{{}\times\bigl(
    \EX{i\pi(s\!+\!2\theta\!+\!1)\frac{[ur'+\pell(s'-1)]_{2\pell}
        +2\pell b}{u}} \Phi(2u\pell\tau,
    -\pell\mu\!-\![ur'\!-\!\pell(s'\!-\!1)]_{2\pell}\tau
    \!-\!2b\pell\tau)}\\*
  {}-\EX{i\pi(s\!+\!2\theta\!+\!1)\frac{-[ur'+\pell(s'-1)]_{2\pell}
      +2\pell (b+1)}{u}} \Phi(2u\pell\tau,
  -\pell\mu\!+\![ur'\!-\!\pell(s'\!-\!1)]_{2\pell}\tau
  \!-\!2(b\!+\!1)\pell\tau)
  \bigr)\\*
  {}\times\varpi_{r',s'}(\tau,\nu).
\end{multline}

To finish the calculation, it remains to recognize the $\chi$
characters in the first term in~\eqref{intermediate-1}.
Lemma~\bref{Lemma:psi-blow-u} gives
\begin{multline}\label{intermediate-3}
  \chi_{r,s%%,\pell,u
    ;\theta}(-\ffrac{1}{\tau},\ffrac{\nu}{\tau},
  \ffrac{\mu}{\tau})
  = \ffrac{(-1)^{r-1}\!\!}{u}\,
  \EX{i\pi k\frac{\nu^2 - (\mu  - \tau)^2}{2\tau}
    -i\pi\frac{\pell}{u}(s \!+\! 1 \!+\! 2 \theta)}\\*
  {\times}\Bigl(\sum_{s'=1}^{u}\!
  \sum_{\theta'=0}^{u-1}\!
  \EX{i\pi\frac{\pell}{u}(s'\!+\!1\!+\!2\theta')
    (s\!+\!1\!+\!2\theta)
    \!-\!i\pi\frac{\pell}{u}(s'\!-\!1)(s\!-\!1)
    \!+\!i\pi\frac{[\pell(s'+1+2\theta')]_u}{u}}\!
  \chi^{\vphantom{y}}_{(s';\theta')}
  (\tau,\nu,\mu\!-\!\tau\!+\!1)
  \\*
  -\boldsymbol{\mathscr{F}}_{r,s}(\tau,\nu,\mu\!-\!\tau\!-\!2\theta)
  \THETA(\tau,\nu,\mu\!-\!\tau\!+\!1)\!\Bigr).
\end{multline}
The characters arising in the right-hand side are those
in~\eqref{chi-new}, and we can therefore restrict to the same
characters in the left-hand side, i.e., choose
$\chi^{\vphantom{y}}_{(s;\theta)}$ as ``representatives'' of the
$\chi_{r,s,\pell,u;\theta}(q,x,y)$ characters with different~$r$
(see~\eqref{r-relation}).  With $\chi^{\vphantom{y}}_{(s;\theta)}$ in
the left-hand side, it is then easy to see that the above $H_{m,m'}$
terms indeed vanish.

The $S$-transform formula~\eqref{intermediate-3} applies to the
characters $\chi_{r,s;\theta}\!\equiv\!\chi_{r,s,\pell,u;\theta}$ in
the \chzerotext{} sector, which are related to super-\chhalftext{}
characters under~$S$, see~\eqref{sR2n}.  Accordingly, we can rewrite
the above formula with the super-\chhalftext{} characters in the
right-hand side,
\begin{multline}\label{intermediate-4}
  \text{1st term in the rhs of~\eqref{intermediate-3}}=
  \ffrac{(-1)^{\rbar{\pell,u}{s+1+2\theta}+1}\!\!}{u}\,
  \EX{i\pi k\frac{\nu^2 - \mu^2}{2\tau}
    -i\pi\frac{\pell}{u}(s \!+\! 2 \!+\! 2 \theta)}\\*
  \times\!\sum_{s'=1}^{u}
  \sum_{\theta'=0}^{u-1}\!
  \EX{i\pi\frac{\pell}{u}(s'\!+\!1\!+\!2\theta')
    (s\!+\!1\!+\!2\theta)
    \!-\!i\pi\frac{\pell}{u}(s'\!-\!1)(s\!-\!1)
    \!+\!i\pi\frac{[\pell(s'+1+2\theta')]_u}{u}}
  \chi^{\schhalf}_{(s';\theta')}
  (\tau,\nu,\mu).
\end{multline}
The result given at the beginning of this section is for the
super-\chzerotext{} sector, which is mapped into itself.
Using~\eqref{sNS2n} and recalling~\eqref{F-final}, we immediately
rewrite~\eqref{intermediate-4} as~\eqref{sR-transform}.

\section{Super-Virasoro characters}
\label{sec:N24} The higher-level Appell functions also arise in
superconformal extensions of the Virasoro algebra.

\subsection{$\N2$ characters}
Following~\cite{[FSST]}, we consider the admissible representations of
the $\N2$ [super-Virasoro] algebra with central charge
$c=3(1-\frac{2}{t})$, $t=\frac{u}{\pell}$.  As in the $\hSSL21$ case, the
spectral flow transform acts freely on the representations, and
therefore the theory is nonrational.  The admissible representation
characters are given by
\begin{gather*}%%%\label{admissible-N2}
  \omega{}_{r,s,u,\pell }(q,z)=
  z^{-\frac{\pell}{u}(r-1)+s-1}\,
  \varphi_{r,s,u,\pell}(q,z)\,
  \frac{\vartheta_{1,0}(q,z)}{
    q^{-\frac{1}{8}}\eta(q)^3},
  \\
  1\leq r\leq u-1,\quad 1\leq s\leq \pell,
\end{gather*}
where
\begin{multline*}%%%\label{varphi}
  \varphi_{r,s,u,\pell}(q,z)=
  \K_{2\pell}(q^u,q^{\frac{r}{2}-\frac{u}{2\pell}(s-1)},
  -z^{-1}\,q^{-\frac{r}{2}+\frac{u}{2\pell}(s-1)})\\*
  {}-q^{r(s-1)}\,
  \K_{2\pell}(q^u,q^{-\frac{r}{2}-\frac{u}{2\pell}(s-1)},
  -z^{-1}q^{-\frac{r}{2}+\frac{u}{2\pell}(s-1)}).
\end{multline*}

The $\N2$ spectral flow acts on the character of any
$\N2$-module~$\mD$ as
\begin{gather*}%%%\label{spectral-n2-general}
  \charn{\mD}{;\theta}(q,z) = z^{-\frac{c}{3}\theta} 
  q^{\frac{c}{6}(\theta^2 - \theta)} 
  \charn{\mD}{}(q,z q^{-\theta}),
\end{gather*}
with $\theta\in\oZ$.  For $\omega{}_{r,s,u,\pell }$ above,
\textit{open quasiperiodicity} occurs for the spectral flow transform
with $\theta=u$~\cite{[FSST]}.
%% For $n\geq1$, we have
%% \begin{multline}\label{phi-identity}
%%   z^{2\pell n} q^{\pell n(r-u n) - u n(s-1)}
%%   \varphi_{r,s,u,\pell}(q,z q^{-u n})-\varphi_{r,s,u,\pell}(q,z)={}\\
%%   {}=\sum_{a=0}^{2\pell n-1} (-1)^{a+1} z^{a+1} q^{-u\pell/4}
%%   \Bigl(\vartheta_{1,0} (q^{2u\pell},q^{-u\pell + u(s+a) - \pell
%%   r})\qquad\\*
%%   {}-q^{r(s+a)}\vartheta_{1,0}(q^{2u\pell},q^{-u\pell+u(s+a)+\pell r})
%%   \Bigr),
%% \end{multline}
%% and for $n\leq-1$, similarly,
%% \begin{multline}\label{phi-identity2-}
%%   z^{2\pell n} q^{\pell n(r-u n) - u n(s-1)}
%%   \varphi_{r,s,u,\pell}(z q^{-u n},q)
%%   - \varphi_{r,s,u,\pell}(q,z)={}\\
%%   =\sum_{a=2\pell n}^{-1} (-1)^{a} z^{a+1}  q^{-u\pell/4}
%%   \Bigl(\vartheta_{1,0}(q^{2u\pell},q^{-u\pell + u(s+a) - \pell r})
%%   -q^{r(s+a)} \vartheta_{1,0}(q^{2u\pell},q^{-u\pell + u(s+a) + \pell r})
%%   \Bigr).
%% \end{multline}
%% and hence the entire character $\omega{}_{r,s,u,\pell }$ is shifted
%% by theta-functional terms.

The modular transformations of $\omega{}_{r,s,u,\pell;\theta}$ can be
derived either by repeating the calculations in Sec.~\ref{sec:sl21} in
the $\N2$ context or by noting that the $\N2$ characters follow by
taking residues of the appropriate $\hSSL21$-characters, and hence
Theorem~\bref{thm:S-chi} implies the $\N2$ modular transformation
formula.  Taking the residues amounts to using the formulas
\begin{gather*}%%%\label{theta1-prime}
  \vartheta_{1,1}(q,q^n)=0,\qquad
  \ffrac{\dd\vartheta_{1,1}(q,z)}{\dd z}\bigm|_{z=q^n}=
  (-1)^{n}\,q^{-\frac{1}{8}}
  \eta(q)^3\,q^{-\frac{n^2}{2} - \frac{3n}{2}},\quad
  n\in\oZ.
\end{gather*}
Noting that
\begin{gather*}
  \psi_{r,s,\pell,u}(q,q^{n},z^2\,q^{2\theta})
  =\varphi_{s-n-1,r,u,\pell}(q,z\,q^{\theta+1+\frac{n}{2}}),
\end{gather*}
we then immediately obtain that for even~$n$,
\begin{gather*}
  \res_{x=q^{n}}\chi_{r,s,\pell,u;\theta}(q,x,z^2)
  =z^{-k}\,
  q^{-k\frac{n^2}{4}+n}\,
  \frac{\vartheta_{1,0}(q,z)}{q^{-\frac{1}{8}}\eta(q)^3}
  \,\omega_{s-n-1,r,u,\pell;-(\theta+\frac{n}{2}+1)}(q,z).
\end{gather*}
The spectral flow transform by $-\theta$ (rather than~$\theta$) in the
right-hand side is due to oppositely chosen conventions for the $\N2$
algebra in~\cite{[FSST]}, which we reproduce for ease of comparison,
and for $\hSSL21$ in~\cite{[ST]}, which we follow here.
The $\hSSL21$-characters then correspond to $\N2$ characters in
accordance with the relation between the $\hSSL21$ and $\N2$
\textit{representations}~\cite{[Sl21-sing],[S-inv]} under the
Hamiltonian reduction~\cite{[BO],[BLNW],[IK]}.

\subsection{$\N4$ characters} Another application of the Appell
functions is to models of the $\N4$ super-Virasoro algebra.  To avoid
lengthening an already sufficiently long paper, we only note that
the\pagebreak[3] unitary irreducible $\N4$ characters at central
charge $c=6k$, $k \in \oN$, derived in~\cite{[ET88]}, can be expressed
through the level-$2(k+1)$ Appell functions as\footnote{An alternative
  form of the character follows by applying
  Eq.~\eqref{diag-periodicity} to each Appell function.}
\begin{multline*}
  \mathrm{Char}^{0}_{k,j}(q,z,y)
  =\frac{q^{j-\frac{2k-1}{8}}}{y^{-1}-y}
  \,\frac{\vartheta(q,zy)\vartheta(q,zy^{-1})}{
    \vartheta_{1,1}(q, z^2)\eta(q)^3}\\*
  {}\times\!\Bigl(yz^{2j}
  \K_{2(k+1)}(q,zq^{\frac{2j+1}{2(k+1)}},
  -yq^{\frac{k-2j}{2(k+1)}})
  -y^{-1}z^{2j}\K_{2(k+1)}(q,zq^{\frac{2j+1}{2(k+1)}},
  -y^{-1}q^{\frac{k-2j}{2(k+1)}})\\*
  \shoveleft{\qquad{}-yz^{-2j}q^{-1}
    \K_{2(k+1)}(q,zq^{\frac{-2j+1}{2(k+1)}},
  -yq^{\frac{2j-k-2}{2(k+1)}})}\\*
  {}+y^{-1}z^{-2j}q^{-1}\K_{2(k+1)}(q,zq^{\frac{-2j+1}{2(k+1)}},
  -yq^{\frac{2j-k-2}{2(k+1)}})
  \Bigr),\quad j =0,\fhalf,\dots,\ffrac{k}{2},
\end{multline*}
which reduces evaluation of the modular transform of the characters to
a calculation with Appell functions based on Theorem~\bref{thm:TS}.

\section{Conclusions}
%%%
We have investigated the modular properties of the Appell functions
and used these to calculate modular transformations of characters in
some nonrational conformal models.  Expressing representation
characters through higher-level Appell functions can be viewed as
going one step up in functional complexity compared with the
characters expressed through theta functions: while the characters are
not quasiperiodic, the quasiperi\-odicity-violating terms are still
given by theta functions.  Efficient manipulations with the $\K_\pell$
functions, as in the study of $\hSSL21$-characters, require using
properties of the $\Phi$ function defined by integral
representation~\eqref{Phi-def} (%%%or~\eqref{R-i0},
which is at the same time the $\boldsymbol{b}$-period integral of
$\K_\pell$, Eq.~\eqref{Phi-Kell}).  We have studied the properties of
$\Phi$ in some detail.\footnote{The function
  $f(\tau,\mu)=-\Phi(-\frac{1}{\tau},\frac{\mu}{\tau})$ has appeared
  in~\cite{[Zh]}, where the role of the integral representation was to
  give a solution of finite-difference equations~\eqref{phi-tau}
  and~\eqref{phi-1}, in a context not unrelated to the present one.}

There are many rational models of conformal field theory, but
nonrational models are also interesting.  The theory of nonrational
models is still in its infancy, however.  The axioms of rational
conformal field theory can be relaxed to different degrees, which in
some cases gives ``almost rational'' theories whose structure may be
worth studying, but difficulties in treating them in the same spirit
as truly rational theories emerge at full scale in calculating the
modular group representation on characters.  The characters of
nonrational models are usually not expressible in terms of theta
functions;
%% whose good modular behavior is a basic property in rational models.
going beyond rational conformal field theories requires an adequate
replacement of theta functions with some functions that are not
quasiperiodic but nevertheless behave reasonably under modular
transformations.

In the examples in this paper, the spectral flow transform action
leads to infinite proliferation of representations, and at first sight
also of characters to be involved in modular transform formulas.  But
the deviation from rational theories may be expected to be ``soft''
because the spectral-flow-transformed representations, although
nonisomorphic, have ``the same'' structure.  It turns out that at the
expense of including extensions among the representations, a modular
group action can be defined on a finite number of characters.
Technically, this was achieved by first studying modular properties of
the Appell functions, which demonstrate a triangular structure in
their behavior under both lattice translations and modular
transformations.

More specifically, we investigated the properties of $\hSSL21$ models
based on the set of admissible representations at rational level.  A
crucial property of these representations is that they allow
nontrivial extensions among themselves.  Such extensions do not occur
in rational theories but are typical of logarithmic conformal field
theories (see~\cite{[Gur],[GK1],[FHST],[flohr-0111228],[G-alg]} and
extensive bibliography therein; such extensions of representations
have been known to play an important role in the derivation of the
modular transformations of $N=4$ superconformal characters since
\cite{[ET88]}, where the corresponding characters are called ``massive
at the unitarity bound.''). It might therefore be expected that the
theory can be consistently formulated as a logarithmic one (i.e.,
further extension of modules results in modules where $L_0$ and$/$or
some Cartan generators act nondiagonally).  There also arises a very
general problem of defining a reasonable class of nonrational
conformal field theory models, where by ``reasonable'' we mean that
the properties known in the rational case are \textit{modified}, but
\textit{not dropped} in going beyond the rational models.  Good
examples are the $(1,p)$ Virasoro models~\cite{[FHST]} and, probably,
the logarithmic extensions of all the $(p',p)$ models.  We hope that
some features of this class of nonrational models have also been
captured in this paper.

We have found the $S$-transform~\eqref{sR-transform} of admissible
characters in the sense that we expressed
$\chi(-\frac{1}{\tau},\frac{\nu}{\tau}, \frac{\mu}{\tau})$ in terms of
the $\chi(\tau,\nu,\mu)$ and $\Omega(\tau,\nu,\mu)$ characters.  At
the next step, we face the ``$S(\tau)$ problem,'' which is a typical
difficulty encountered in nonrational conformal field theories: the
matrix representing the action of $\left(\begin{smallmatrix}0~&-1\\
    1~&0\end{smallmatrix}\right)\in\SLiiZ$ on the characters acquires
dependence on coordinates on the moduli space,
\begin{gather*}
  \boldsymbol{\chi}(\left(\begin{smallmatrix}
      0~&-1\\
      1~&0
    \end{smallmatrix}
  \right)\ldot x)=S(x)\boldsymbol{\chi}(x),
\end{gather*}
where $\boldsymbol{\chi}(x)$ is a vector whose entries are the
characters and $x$ denotes coordinates on the moduli space
($x=(\tau,\nu,\dots)$, with $\gamma\ldot(\tau,
\nu,\dots)=(\frac{a\tau+b}{c\tau+d},\frac{c\nu}{c\tau+d},\dots)$ for
$\gamma=\left(\begin{smallmatrix}
    a~&b\\
    c~&d
  \end{smallmatrix}\right)
\in\SLiiZ$).  The problem is that the matrix $S(x)$ depends on $x$ in
general, making the Verlinde formula in its standard form
inapplicable.

The general strategy to deal with the ``$S(x)$ problem'' was outlined
in~\cite{[FHST]}.  The modular group action on the characters is to be
redefined as
\begin{gather*}
  \gamma*
  \boldsymbol{\chi}(x)=J(\gamma;x)\boldsymbol{\chi}(\gamma\ldot x),
\end{gather*}
with a matrix automorphy factor $J(\gamma;x)$ such that the matrix
\begin{gather}\label{new-S}
  S=J(\left(\begin{smallmatrix}
      0~&-1\\
      1~&0
    \end{smallmatrix}\right);x)S(x)
\end{gather}
is numerical (independent of the coordinates on the moduli
space).\footnote{We note that the $(\tau,\nu,\mu)$-dependence through
  the standard scalar factor $e^{i\pi k\frac{\nu^2 - \mu^2}{2\tau}}$
  is in fact eliminated similarly, with the standard scalar automorphy
  factor.  Any automorphy factor must satisfy the cocycle equation;
  also see~\cite{[EhSk]} for the matrix case.}  Then $S$ defined
in~\eqref{new-S} can be used in a Verlinde-like formula (an example of
successful application of this ideology is given in~\cite{[FHST]}).
For the $\hSSL21$-characters, with the vector $\boldsymbol{\chi}$
composed of the $\chi$ \textit{and} $\Omega$ characters, the most
essential part of the ``$S(x)$ problem'' is the $\tau$- and
$\mu$-dependence in the $\Phi$ functions arising in the $S$-transform.
The $\Phi$ functions, studied in Sec.~\ref{sec:Phi}, are a
characteristic element of the Appell function theory.  We leave this
problem for a future work.

We finally note that the derivation of modular transformation
properties of $\K_\pell$ and the characters given above may not be
``optimized'' --- apart from technical improvements, a more
``conceptual'' derivation must exist, possibly applicable to more
general indefinite theta series.

\subsubsection*{Acknowledgments}
We are grateful to B.L.~Feigin for interesting discussions, to
J.~Fuchs for a useful suggestion, and to V.I.~Ritus for his help with
the small-$t$ asymptotic expansion.  AMS acknowledges support from the
Royal Society through a grant RCM/ExAgr and the kind hospitality in
Durham.  AT acknowledges support from a Small Collaborative Grant of
the London Mathematical Society that made a trip to Moscow possible,
and the warm welcome extended to her during her visit.  AMS\;\&\;IYuT
were supported in part by the grant LSS-1578.2003.2, by the Foundation
for Support of Russian Science, and by the RFBR Grant
%% ~01-01-00906, superseded by 
04-01-00303.  IYuT was also supported in part by the RFBR Grant
03-01-06135 and the INTAS Grant~00-01-254.

\appendix

\section{The $\acycle$- and $\bcycle$-Cycle Integrals on the Torus}
\label{app:ab-int}
We first evaluate the integral along the $\bcycle$ cycle on the torus,
\begin{gather*}
  \oint_{\bcycle} d\lambda\,
  \EX{i\pi\frac{\lambda^2}{\tau}}\,
  \vartheta(\tau,\lambda)
  = \int\limits_{0}^{\tau}
  d\lambda\,\EX{i\pi\frac{\lambda^2}{\tau}}\,
  \vartheta(\tau,\lambda)
  =\sum_{m\in\oZ}\int\limits_{m\tau}^{(m+1)\tau}
  \!d\lambda\,\EX{i\pi\frac{\lambda^2}{\tau}},
\end{gather*}
where we shifted the integration variable as $\lambda\to\lambda -
m\tau$ in each term of the $\vartheta$-series.  For $\Im\tau>0$, the
integrals are defined by analytic continuation from $\tau=i t$ with
$t\in\oR_{>0}$, and therefore
\begin{gather} \label{theta-b}
  \oint_{\bcycle} d\lambda\,\EX{i\pi\frac{\lambda^2}{\tau}}\,
  \vartheta(\tau,\lambda)
  =i\int_{\oR} dx\,\EX{-\pi
    \frac{x^2}{t}}\Bigm|_{t=-i\tau}= i\sqrt{-i\tau}.
\end{gather}
The ``dual'' integral is, obviously,
\begin{gather} \label{theta-a}
  \oint_{\acycle} d\lambda\, \vartheta(\tau,\lambda)
  \stackrel{\mathrm{def}}{=}
  \int_{0}^{1}\!d\lambda\, \vartheta(\tau,\lambda) = 1.
\end{gather}

Somewhat more generally than in~\eqref{theta-b}, we can evaluate the
integral 
\begin{gather*}
  \int_{\mu}^{\mu+c\tau+d}\!d\lambda\,
  \EX{i\pi\frac{c\lambda^2}{c\tau+d}}\,
  \vartheta(\tau,\lambda),\qquad c\in\oN,\quad d\in\oZ,
  \quad cd\in2\oZ,
\end{gather*}
(with an arbitrary~$\mu$) by shifting the integration variable as
$\lambda\to\lambda-m(\tau+\frac{d}{c})$ in each term of the theta
series and then summing over $m$ as $\sum_{m\in\oZ}f(m)
=\sum_{\alpha=1}^{c}\sum_{m\in\oZ}f(c m + \alpha)$.  With even $cd$,
it then follows that $m$ drops from the exponentials, and we readily
obtain
\begin{align}\label{theta-cd}
  \int\limits_{\mu}^{\mu+c\tau+d}\!\!d\lambda\,
  \EX{i\pi\frac{c\lambda^2}{c\tau+d}}\,
  \vartheta(\tau,\lambda)
  ={}&\sum_{m\in\oZ}\!\!
  \int\limits_{\alpha(\tau+\frac{d}{c}) + m(c\tau+d)}^{
    \alpha(\tau+\frac{d}{c}) + (m+1)(c\tau+d)}\!\!
  d\lambda\,
  \EX{i\pi\frac{c\lambda^2}{c\tau+d}}\cdot
  \sum_{\alpha=1}^{c}\EX{-i\pi\alpha^2\frac{d}{c}}\\
  {}={}&i\sqrt{\ffrac{-i(c\tau+d)}{c}}\cdot
  \sqrt{-i c}\,\zeta_{c,d}^{-1}
  =(c\tau+d)^{\half}\,\zeta_{c,d}^{-1},\notag
\end{align}
where the two factors are a Gaussian integral, calculated by analytic
continuation (for $\tau$ such that $\Im\tau>0$) from the integral
over~$\oR$, and a Gaussian sum, see, e.g.,~\cite{[Mum]}.

Remarkably, much similarity is preserved if theta functions are
replaced with Appell functions in the above integrals.  We first
consider the corresponding analogue of~\eqref{theta-b},
\begin{align*}
  \oint_{\bcycle} d\lambda\,
  \EX{i\pi\frac{\lambda^2 - 2\lambda\mu}{\tau}}
  \K_1(\tau,\lambda\!-\!\mu,\mu)&{}
  \stackrel{\mathrm{def}}{=}
  \int_{0}^{\tau} d\lambda\,
  \EX{i\pi\frac{\lambda^2 - 2\lambda\mu}{\tau}}\,
  \K_1(\tau,\lambda+\varepsilon-\mu,\mu),\\
  \intertext{where an infinitesimal positive real $\varepsilon$
    specifies the prescription to bypass the singularities.  Again
    continuing from $\tau=it$ and $\mu=iy$ with positive real $t$ and
    real~$y$, we have}
  &{}=i\!\sum_{m\in\oZ}\int_0^{t}\!\!dx\,
  \EX{-\pi\frac{x^2 - 2 xy}{t}}\,
  \ffrac{
    \EX{-\pi t^2 m^2\!-\!2\pi m(x{-}y)}}{
    1\!-\!\ex{-2\pi(x {+} m t)}\!-\!i0}
  \biggm|_{\substack{
      t=-i\tau\\
      y=-i\mu}}\!\!.
\end{align*}
Making the same substitution $\lambda\to\lambda - m\tau$ as above, or
$x\to x - mt$, and using that
\begin{gather*}
  \int_{-\infty}^{+\infty} dx\,
  \ffrac{f(x)}{x - i0}=
  \dashint_{-\infty}^{+\infty} dx\,
  \ffrac{f(x)}{x}
  + i\pi\res_{x=0} f(x),
\end{gather*}
we obtain
\begin{gather}\label{Phi-K}
  \oint_{\bcycle} d\lambda\,
  \EX{i\pi\frac{\lambda^2 - 2\lambda\mu}{\tau}}
  \K_1(\tau,\lambda-\mu,\mu)
  = -i\sqrt{-i\tau}\,\Phi(\tau,\mu),
\end{gather}
with $\Phi$ defined in~\eqref{Phi-def}.  The derivation shows that the
same result is valid for the ``$\bcycle$''-integral with a translated
contour:
\begin{gather}\label{Phi-K-shift}
  \int_{\alpha\tau}^{\tau+\alpha\tau} d\lambda\,
  \EX{i\pi\frac{\lambda^2 - 2\lambda\mu}{\tau}}\,
  \K_1(\tau,\lambda+\varepsilon-\mu,\mu)
  = -i\sqrt{-i\tau}\,\Phi(\tau,\mu),\qquad\alpha\in\oR.
\end{gather}

The version of~\eqref{Phi-K} for $\K_\pell$ is given
by%%\enlargethispage{12pt}
\begin{gather}\label{Phi-Kell}
  \oint_{\bcycle} d\lambda\,
  \EX{i\pi\pell\frac{\lambda^2 - 2\lambda\mu}{\tau}}
  \K_\pell(\tau,\lambda-\mu,\mu)
  = -i\sqrt{\ffrac{-i\tau}{\pell}}\,\Phi(\ffrac{\tau}{\pell},\mu).
\end{gather}
More generally than in~\eqref{Phi-K}, we can evaluate the
integral
\begin{multline} \label{K-cd}
  \int_{\alpha\tau}^{c\tau+d+\alpha\tau} \!d\lambda\,
  \EX{i\pi\frac{\lambda^2 + 2\lambda\mu}{\tau+\frac{d}{c}}}
  \K_1(\tau,-\lambda+\varepsilon-\mu,\mu)\\*
  =\smash[b]{-i\sqrt{-i(\tau\!+\!\ffrac{d}{c})}
    \sum_{r=0}^{c-1}
    \EX{2i\pi\mu\frac{r d}{c\tau+d}
      \!-\!i\pi\tau\frac{r^2 d}{c\tau+d}}\,
    \Phi(\tau\!+\!\ffrac{d}{c},\mu\!+\!r\ffrac{d}{c})},
  \quad\alpha\in\oR
\end{multline}
for $c\in\oN$, $d\in\oZ$, and~$cd\in2\oZ$.

Similarly to~\eqref{theta-a}, we have the dual, $\acycle$-cycle
integral
\begin{gather} \label{K-a}
  \oint_{\acycle} d\lambda\K_1(\tau,\lambda - \mu, \mu)
  \stackrel{\mathrm{def}}{=}
  \int_{0}^{1}\!
  d\lambda\K_1(\tau,\lambda - \mu +i0, \mu) =
  1.
\end{gather}

\section{$\protect\widehat{s\ell}(2|1)$ Algebra, its Automorphisms and
  Modules}\label{app:sl21}
\subsection{The algebra and automorphisms}\label{app:the-algebra} The
affine Lie superalgebra $\hSSL21$ is spanned by four bosonic currents
$E^{12}$, $H^-$, $F^{12}$, and $H^+$, four fermionic ones, $E^1$,
$E^2$, $F^1$, and $F^2$, and the central element (which we identify
with its eigenvalue~$k$).  The $\hSL2$ subalgebra is generated by
$E^{12}$, $H^-$, and $F^{12}$, and it commutes with the~$u(1)$
subalgebra generated by $H^+$.  The fermions $E^1$ and $F^2$ on one
hand and $F^1$ and $E^2$ on the other hand form $\SL2$ doublets.  The
nonvanishing commutation relations are
\begin{gather}\label{sl21}
  \begin{alignedat}{2}
    {[}H^-_m, E^{12}_n] &= E^{12}_{m+n},&
    {[}H^-_m, F^{12}_n] &= -F^{12}_{m+n},\\
    {[}E^{12}_m, F^{12}_n] &= m \delta_{m+n, 0} k + 2
    H^-_{m+n},\quad& 
    {[}H^\pm_m, H^\pm_n] &= \mp\fhalf m \delta_{m+n, 0} k,\\
    {[}F^{12}_m, E^2_n] &= F^1_{m+n},&
    {[}E^{12}_m, F^2_n] &= -E^1_{m+n},\\
    {[}F^{12}_m, E^1_n] &= -F^2_{m+n},&
    {[}E^{12}_m, F^1_n] &= E^2_{m+n},\\
    {[}H^\pm_m, E^1_n] &= \fhalf  E^1_{m+n},&
    {[}H^\pm_m, F^1_n] &= -\fhalf  F^1_{m+n},\\
    {[}H^\pm_m, E^2_n] &= \mp\fhalf  E^2_{m+n},&
    {[}H^\pm_m, F^2_n] &= \pm\fhalf  F^2_{m+n},\\
    {[}E^1_m, F^1_n] &= -m \delta_{m+n, 0} k +
      H^+_{m+n} -
      H^-_{m+n},\kern-200pt\\
    {[}E^2_m, F^2_n] &= m \delta_{m+n, 0} k + 
      H^+_{m+n} +
      H^-_{m+n},\kern-200pt\\
    {[}E^1_m, E^2_n] &= E^{12}_{m+n},&
    {[}F^1_m, F^2_n] &= F^{12}_{m+n}.
  \end{alignedat}
\end{gather}
The Sugawara energy-momentum tensor is given by
\begin{gather}\label{Tsug-sl21}
  T_{\mathrm{Sug}} =
  \ffrac{1}{k +1}\bigl(\Hminus \Hminus - \Hplus \Hplus +
  E^{12} F^{12} + E^1 F^1 - E^2 F^2 \bigr).
\end{gather}

There are involutive algebra automorphisms
\begin{align}
  \alpha:{}&
  \begin{alignedat}{3}
    E^1_n&\mapsto F^2_n,\qquad&E^2_n&\mapsto F^1_n,\qquad&E^{12}_n&
    \mapsto
    F^{12}_n,\\
    F^1_n&\mapsto E^2_n,&F^2_n&\mapsto E^1_n,&F^{12}_n&\mapsto
    E^{12}_n,\\
    H^+_n&\mapsto H^+_n,\qquad&H^-_n&\mapsto- H^-_n,&&
  \end{alignedat}\\
  \intertext{and}
  \beta:{}&
  \begin{alignedat}{3}
    E^1_n&\mapsto E^2_n,\qquad&E^2_n&\mapsto E^1_n,\qquad&E^{12}_n&
    \mapsto
    E^{12}_n,\\
    F^1_n&\mapsto-F^2_n,&F^2_n&\mapsto-F^1_n,&F^{12}_n&\mapsto
    F^{12}_n,\\
    H^+_n&\mapsto-H^+_n,\qquad& H^-_n&\mapsto H^-_n,&&
  \end{alignedat}
  \label{beta} \\
  \intertext{and a $\oZ$ subgroup of automorphisms called the spectral
    flow,}
  \cU_\theta:{}&
  \begin{alignedat}{2}
    E^1_n&\mapsto E^1_{n-\theta},\qquad&E^2_n&\mapsto E^2_{n+\theta},
    \\
    F^1_n&\mapsto F^1_{n+\theta},&F^2_n&\mapsto F^2_{n-\theta},
  \end{alignedat}
  \quad H^+_n\mapsto H^+_n + k\theta\delta_{n,0},
  \label{sl21-spectral}
\end{align}
where $\theta\in\oZ$ (and the~$\hSL2$ subalgebra remains invariant).
We note the relations
\begin{gather}
  \alpha^2=1,\qquad\beta^2=1,\qquad(\alpha\beta)^4=1,\qquad
  \alpha\cU_\theta=\cU_\theta\alpha,\qquad(\beta\cU_\theta)^2=1.
\end{gather}

Another $\oZ$ algebra of automorphisms (a spectral flow affecting the
$\hSL2$ subalgebra, cf.~\cite{[FSST]}) acts as
\begin{gather}\label{sl21-spectral2}
  \cA_\eta:{}
  \begin{gathered}
    \begin{alignedat}{3}
      E^1_n&\mapsto  E^1_{n+\eta},&\qquad
      E^2_n&\mapsto E^2_{n+\eta},&\qquad
      E^{12}_n&\mapsto E^{12}_{n+2\eta},
      \\
      F^1_n&\mapsto  F^1_{n-\eta},&F^2_n&\mapsto F^2_{n-\eta},&
      F^{12}_n&\mapsto F^{12}_{n-2\eta},
    \end{alignedat}\\
    H^-_n\mapsto H^-_n + k\eta\delta_{n,0},\qquad
    H^+_n\mapsto H^+_n
  \end{gathered}
\end{gather}
There also exists the automorphism
$\gamma=\cU_{\half}\circ\cA_{-\half}$ (while $\cU_{\half}$ and
$\cA_{-\half}$ are not automorphisms, but rather mappings into an
isomorphic algebra, their composition is).  For $\theta\in\oZ$, its
powers $\cT_\theta=\gamma^\theta$ map the generators as
\begin{gather}\label{sl21-spectral3}
  \cT_\theta:{}
  \begin{aligned}%%%{6}{}
    {}E^1_n&\mapsto E^1_{n-\theta},~
    &E^2_n&\mapsto E^2_{n},~
    &E^{12}_n&\mapsto E^{12}_{n-\theta},~
    &H^-_n&\mapsto H^-_n - \ffrac{k}{2}\theta\delta_{n,0},
    \\
    F^1_n&\mapsto  F^1_{n+\theta},&F^2_n&\mapsto F^2_{n},&
    F^{12}_n&\mapsto  F^{12}_{n+\theta},
    &H^+_n&\mapsto  H^+_n + \ffrac{k}{2}\theta\delta_{n,0}.
  \end{aligned}
\end{gather}

Spectral flow transform~\eqref{sl21-spectral}, affecting the fermions
and leaving the $\hSL2$ subalgebra invariant, plays an important role
in the study of $\hSSL21$-representations~\cite{[ST]}.  We use the
notation
\begin{gather*}
  \mP_{;\theta}\equiv\cU_\theta\mP
\end{gather*}
for the action of spectral flow transform on any $\hSSL21$-module
$\mP$.  Obviously, $\mP_{;0}\equiv\mP$.  For a module $\mP$, we let
\begin{gather}
  \chi[\mP](q,x,y)=\Tr^{\phantom{y}}_{\mP}
  \bigl(
  q^{L_0}\,x^{H^-_0}\,y^{H^+_0}
  \bigr)
\end{gather}
be its character.  The character $\chi[\mP]_{;\theta}$ of the
spectral-flow transformed module ${\mP}_{;\theta}$ is expressed
through the character of $\mP$ as
\begin{gather}\label{sl21-sf}
  \chi[\mP]_{;\theta}(q,x,y)=
  y^{-k\theta}\,q^{-k\theta^2}\,\chi[\mP](q,x,y\,q^{2\theta}).
\end{gather}

%% Writing ${\mP}_{|\theta}\equiv\cT_{\theta}\mP$ for the action
%% of~\eqref{sl21-spectral3} on modules and $\chi[\mP]_{|\theta}$ for
%% its action on characters, we have
%% \begin{gather}
%%   \chi[\mP]_{|\theta}(q,x,y)
%%   =x^{\frac{k\theta}{2}}\,y^{-\frac{k\theta}{2}}\,
%%   \chi[\mP](q,x q^\theta,y q^\theta)
%% \end{gather}
%% for any $\hSSL21$-module~$\mP$.

\subsection{Highest-weight conditions and
  modules}\label{app:sl21-modules} A significant role in the $\hSSL21$
representation theory is played by the spectral flow
transform~\eqref{sl21-spectral}, which is a family of $\hSSL21$
automorphisms.  Applying algebra automorphisms to modules gives
nonisomorphic modules in general.  The (upper) triangular subalgebra
is also mapped under the action of automorphisms, and the annihilation
conditions satisfied by highest-weight vectors change accordingly.
Thus, the existence of an automorphism group leads to a freedom in
choosing the type of annihilation conditions imposed on highest-weight
vectors in highest-weight representations (in particular, Verma
modules).  We consider the family of annihilation
conditions
\begin{gather}\label{twisted-hwcond}
  \begin{aligned}
    E^1_{\geq-\theta}\approx{}&0,& E^2_{\geq\theta}\approx{}&0,&&\\
    F^1_{\geq\theta+1}\approx{}&0,\quad&
    F^2_{\geq1-\theta}\approx{}&0, \quad&F^{12}_{\geq1}\approx{}&0,
  \end{aligned}
  \qquad
  \theta\in\oZ
\end{gather}
that are an orbit of~$\cU_\theta$
%% for each fixed $\theta\in\oZ$ 
(the $\approx$ sign means that the left-hand sides must be applied to
a vector; at the moment, we are interested in the list of annihilation
operators, rather than in the vector, hence the notation).  These
annihilation conditions are called the \textit{twisted highest-weight
  conditions} in general.  By the $\hSSL21$ commutation relations, the
conditions explicitly written in~\eqref{twisted-hwcond} imply that
$E^{12}_{\geq0}\approx0$ and $H^{\pm}_{\geq1}\approx0$.  This is
understood in similar relations in what follows.

Accordingly, a twisted Verma module $\verma{h_-,h_+,k;\theta}$ over
the level-$k$ $\hSSL21$ algebra is freely generated by
$E^1_{\leq-\theta-1}$, $E^2_{\leq\theta-1}$, $E^{12}_{\leq-1}$,
$F^1_{\leq\theta}$, $F^2_{\leq-\theta}$, $F^{12}_{\leq0}$,
$H^-_{\leq-1}$, and $H^+_{\leq-1}$ from the \textit{twisted
  highest-weight state}~$\ket{h_-,h_+,k;\theta}$ satisfying
annihilation conditions~\eqref{twisted-hwcond} and additionally fixed
by the eigenvalue relations
\begin{align}\label{twisted-eigen-minus}
  H^-_0\,\ket{h_-,h_+,k;\theta} = {}&h_-\,\ket{h_-,h_+,k;\theta},\\
  H^+_0\,\ket{h_-,h_+,k;\theta} = {}&(h_+ -
  k\theta)\,\ket{h_-,h_+,k;\theta}.
  \label{twisted-eigen-plus}
\end{align}
With the parameterization of the $H^+_0$ eigenvalue chosen
in~\eqref{twisted-eigen-plus}, we have
\begin{gather}
  \cU_{\theta'}\ket{h_-,h_+,k;\theta}=\ket{h_-,h_+,k;\theta+\theta'}
\end{gather}
and, obviously, $\cU_{\theta'}\verma{h_-,h_+,k;\theta}=
\verma{h_-,h_+,k;\theta+\theta'}$ for the Verma modules.  This simple
behavior of~$\ket{h_-,h_+,k;\theta}$ under the spectral flow explains
the subtraction of $k\theta$ in~\eqref{twisted-eigen-plus}.

The character of $\verma{h_-,h_+,k;\theta}$ is given by
\begin{gather}\label{sl21-Verma-char}
  \chi[\verma{}]_{h_-,h_+,k;\theta}(q,x,y)
%%   \equiv
%%   \Tr^{\phantom{y}}_{\verma{h_-,h_+,k;\theta}}
%%   \bigl(q^{\cL^{\text{Sug}}_0}\,x^{H^-_0}\,y^{H^+_0}\bigr)=
%%   \\
  = x^{h_-}\,y^{h_+-(k+1)\theta}\, q^{\frac{h_-^2 - h_+^2}{k+1} +
    2\theta h_+ -(k+1)\theta^2}\,
  \THETA(q,x,y),  
\end{gather}
where
\begin{gather}\label{THETA}
  \THETA(q,x,y)=
  \smash[t]{\frac{    
      \vartheta_{1,0}(q, x^{\half}y^{\half})
      \vartheta_{1,0}(q, x^{\half}y^{-\half})
    }{
      \vartheta_{1,1}(q,x)\,
      q^{-\frac{1}{8}}\,\eta(q)^3}}.
\end{gather}

We let $\ketminus{h_-,k;\theta}$ denote a state satisfying the
highest-weight conditions
\begin{gather}\label{hwtop-minus}
  \begin{aligned}
    E^1_{\geq-\theta}\approx{}&0,& \qquad 
    E^2_{\geq\theta}\approx{}&0,\\
    F^1_{\geq\theta}\approx{}&0,& 
    F^2_{\geq1-\theta}\approx{}&0
  \end{aligned}
\end{gather}
and the eigenvalue relations
\begin{gather}\label{hwtop-eigen}
  H^-_0\ketminus{h_-,k;\theta}= h_-\ketminus{h_-,k;\theta}.
\end{gather}
Conditions~\eqref{hwtop-minus} are stronger
than~\eqref{twisted-hwcond}, and as a result, the eigenvalue of
$H^+_0$ is no longer an independent parameter.  We let
$\vermaNminus{h_-,k;\theta}$ denote the (twisted) \textit{narrow Verma
  module}\,---\,the module freely generated by $E^1_{\leq-\theta-1}$,
$E^2_{\leq\theta-1}$, $E^{12}_{\leq-1}$, $F^1_{\leq\theta-1}$,
$F^2_{\leq-\theta}$, $F^{12}_{\leq0}$, $H^-_{\leq-1}$, and
$H^+_{\leq-1}$ from $\ketminus{h_-,k;\theta}$.\footnote{The name
  \textit{Verma} module is a (very convenient) abuse of terminology.
  The $\vermaNminus{}$ modules, as well as $\vermaNplus{}$ introduced
  momentarily, occur as submodules generated from a charged singular
  vector in the proper Verma modules $\verma{}$~\cite{[ST]}.  The
  modules are called \textit{narrow} for the reason explained
  in~\cite{[ST]} (essentially because they \textit{are} narrow
  compared with the proper Verma modules).}

Similarly, let $\ketplus{h_-,k;\theta}$ denote the states satisfying a
different set of the highest-weight conditions
\begin{gather}\label{hwtop-plus}
  \begin{aligned}
    E^1_{\geq-\theta}\approx{}&0,& \qquad 
    E^2_{\geq\theta}\approx{}&0,\\
    F^1_{\geq\theta+1}\approx{}&0,& 
    F^2_{\geq-\theta}\approx{}&0
  \end{aligned}
\end{gather}
(which are again stronger than~\eqref{twisted-hwcond}) and the
eigenvalue relations
\begin{gather}
  H^-_0\ketplus{h_-,k;\theta}= h_-\ketplus{h_-,k;\theta}.
\end{gather}
We write $\vermaNplus{h_-,k;\theta}$ for the corresponding (twisted)
narrow Verma module freely generated by $E^1_{\leq-\theta-1}$,
$E^2_{\leq\theta-1}$, $E^{12}_{\leq-1}$, $F^1_{\leq\theta}$,
$F^2_{\leq-\theta-1}$, $F^{12}_{\leq0}$, $H^-_{\leq-1}$, and
$H^+_{\leq-1}$ from $\ketplus{h_-,k;\theta}$.

The characters of $\vermaNminus{h_-,k;\theta}$ and
$\vermaNplus{h_-,k;\theta}$ are given by
\begin{gather*}%\label{char-minus}
  \chi[\vermaNminus{}]_{h_-,k;\theta}(q,x,y)=
  \frac{\chi[\verma{}]_{h_-,h_-,k;\theta}(q,x,y)}
  {1 + q^{-\theta} x^{-\half}y^{-\half}},
  \\
  \chi[\vermaNplus{}]_{h_-,k;\theta}(q,x,y)=
  \frac{\chi[\mP{}]_{h_-,-h_-,k;\theta}(q,x,y)}
  {1 + q^{\theta} x^{-\half}y^{\half}},
\end{gather*}
where $\chi[\mP{}]$ is defined in~\eqref{sl21-Verma-char}.  The
twisted narrow Verma modules are convenient in constructing a
resolution of the admissible representations $\mL_{r,s}$, see the next
subsection.

\subsection{Admissible $\hSSL21$
  representations~$\mL_{r,s,\pell,u;\theta}$}\label{app:sl21-adm} The
admissible $\hSSL21$-representa\-tions, which belong to the class of
irreducible highest-weight representations characterized by the
property that the corresponding Verma modules are \textit{maximal}
elements with respect to the (appropriately defined) Bruhat order,
have arisen in a vertex-operator extension of two $\hSL2$ algebras
with the ``dual'' levels $k$ and $k'$ such that
$(k+1)(k'+1)=1$~\cite{[BFST]}; via this extension
$\hSL2_{k}\oplus\hSL2_{k'}\to\hSSL21_{k}$, the admissible $\hSSL21$
representations are related to the admissible $\hSL2$
representations~\cite{[KW0]}.  We fix the $\hSSL21$ level~as
\begin{gather*}
  k=\ffrac{\pell}{u}-1
\end{gather*}
with coprime positive integers $\pell$ and $u$.

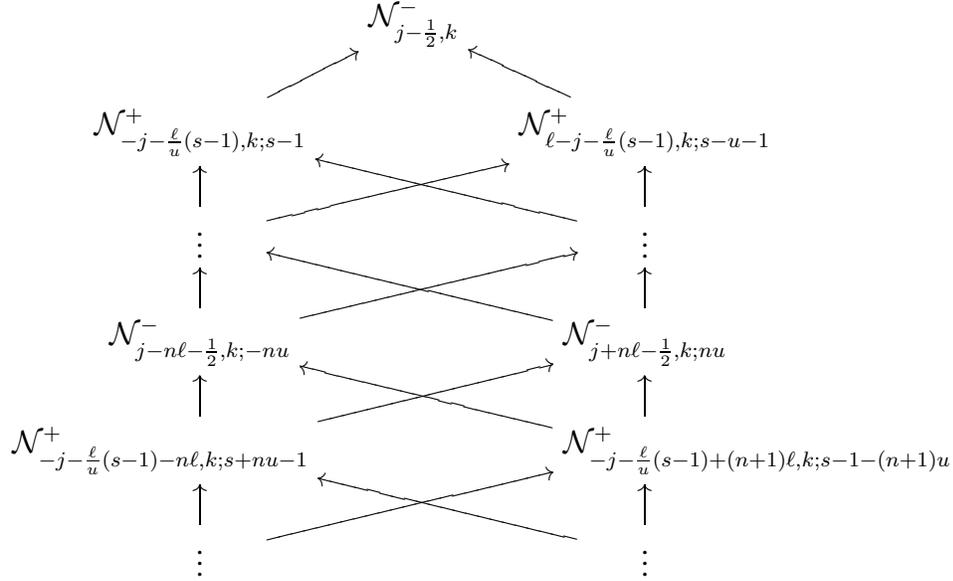
\begin{figure}[tb]
  \mbox{}\qquad\quad\xymatrix@=15pt{
    &\vermaNminus{j -\half,k}&\\
    {\vermaNplus{-j-\frac{\pell}{u}(s-1), k;s-1}}\ar[-1,1]
    &&{\vermaNplus{\pell-j-\frac{\pell}{u}(s-1),k;s-u-1}}\ar[-1,-1]\\
    {\kern20pt\vdots\kern20pt}\ar[-1,0]\ar[-1,2]
    &&{\kern20pt\vdots\kern20pt}\ar[-1,0]\ar[-1,-2]\\
    {\vermaNminus{j - n\pell -\half,k;-n u}}\ar[-1,0]\ar[-1,2]
    &&{\vermaNminus{j+n\pell-\half,k;n u}}\ar[-1,0]\ar[-1,-2]\\
    {\kern-30pt\vermaNplus{-j-\frac{\pell}{u}(s-1) - n\pell, k;s + n u
        - 1}}\ar[-1,0]\ar[-1,2]
    &&{\vermaNplus{-j-\frac{\pell}{u}(s-1)+(n+1)\pell,k;s-1-(n+1)u}
      \kern-85pt}\ar[-1,0]\ar[-1,-2]\\
    {\kern20pt\vdots\kern20pt}\ar[-1,0]\ar[-1,2]
    &&{\kern20pt\vdots\kern20pt}\ar[-1,0]\ar[-1,-2]\\
  }
  \caption[Mappings between narrow Verma modules]{Mappings between
    narrow Verma modules.}\label{fig:grid}
\end{figure}

For $1 \leq r \leq \pell$ and $1 \leq s \leq u$, the admissible $\hSSL21$
representations $\mL_{r,s,\pell,u;\theta}$ is the irreducible quotient
of the Verma module $\verma{\frac{r}{2}-\frac{\pell}{u}\frac{s-1}{2},
  \frac{r}{2}-\frac{\pell}{u}\frac{s+1}{2},
  \frac{\pell}{u}-1;\theta}$.  We omit the level $k=\frac{\pell}{u}-1$
in $\verma{h_-,h_+,k;\theta}$ and similar notation in what follows.
In this Verma module, there is the charged singular vector given by
$E^1_{-\theta-1}F^2_{-\theta}$ acting on the twisted highest-weight
vector $\ket{\frac{r}{2}-\frac{\pell}{u}\frac{s-1}{2},
  \frac{r}{2}-\frac{\pell}{u}\frac{s+1}{2};\theta}$.  The
corresponding submodule is the narrow Verma module
$\vermaNminus{\frac{r}{2}-\frac{s-1}{2}\frac{\pell}{u}; \theta+1}$ and
the quotient is the narrow Verma module $\vermaNminus{\frac{r-1}{2} -
  \frac{\pell}{u}\frac{s-1}{2};\theta+1}$,
\begin{equation}\label{exact}
  0\to
  \vermaNminus{\frac{r}{2}-\frac{s-1}{2}\frac{\pell}{u};
    \theta+1}\to
  \verma{\frac{r}{2}-\frac{\pell}{u}\frac{s-1}{2},
    \frac{r}{2}-\frac{\pell}{u}\frac{s+1}{2};\theta}\to
  \vermaNminus{\frac{r-1}{2} - \frac{\pell}{u}\frac{s-1}{2},
    \frac{\pell}{u}-1;\theta+1}\to
  0.
\end{equation}
The admissible representation $\mL_{r,s,\pell,u;\theta}$ is therefore
the irreducible quotient of the latter narrow Verma module.  Combining
the canonical mapping $\vermaNminus{\frac{r-1}{2} -
  \frac{\pell}{u}\frac{s-1}{2};\theta+1}\to\mL_{r,s,\pell,u;\theta}$
with the mappings in Fig.~\ref{fig:grid}, taken from~\cite{[ST]},
with $j-\half=\frac{r-1}{2} - \frac{\pell}{u}\frac{s-1}{2}$, we obtain
a resolution of the admissible representation.  The resolution readily
implies a character formula.

\begin{Thm}\label{thm:find-chars}
  For $1 \leq r \leq \pell$, $1 \leq s \leq u$, and $\theta\in\oZ$, the
  character of $\mL_{r,s,\pell,u;\theta}$ is given by
  \begin{multline}\label{eq:gen-char}
    \chi_{r,s,\pell,u;\theta}(q,x,y)=
    q^{(\theta+1)(r-1-\frac{\pell}{u}(s+\theta))}\,
    x^{\frac{r-1}{2} - \frac{s-1}{2}\frac{\pell}{u}}\,
    y^{\frac{r-1}{2}-\frac{s+1}{2}\frac{\pell}{u} -
      \theta\frac{\pell}{u}}
    \\*
    {}\times  
    \psi_{r,s,\pell,u}(q,x,y q^{2\theta})\,
%%     \frac{\vartheta_{1,0}(q,x^{\half}y^{\half})
%%       \vartheta_{1,0}(q,x^{\half}y^{-\half})}{
%%       \vartheta_{1,1}(q,x)\prod_{i\geq1}(1-q^i)^3},
    \THETA(q,x,y),
  \end{multline}
  where
  \begin{multline}\label{def-psi}
    \psi_{r,s,\pell,u}(q,x,y)
    = \sum_{m\in\oZ} q^{m^2 u \pell - m u(r-1)}
    \Bigl(\frac{q^{m\pell(s-1)}x^{-m\pell}}{
      1+y^{-\half}x^{-\half}q^{mu-1}}\\*
    - q^{(s-1)(r-1)}\,x^{1-r}\frac{q^{-m\pell(s-1)}x^{m\pell}}{
      1+y^{-\half}x^{\half}q^{mu-s}} \Bigr)
  \end{multline}
  \textup{(}and $\THETA$ is defined in~\eqref{THETA}\textup{)}.
\end{Thm}
\noindent
For $r<p$, this is proved by a straightforward summation of the
characters of the twisted narrow Verma modules involved in the
resolution.  For $r=p$, the above character formula follows from a
somewhat different resolution, with a special role played by one of
the charged singular vectors (see~\cite{[ST]}).

Unless $\pell=1$, the discrete $\hSSL21$ automorphism $\beta$ (see
Appendix~\ref{app:sl21}) maps the $\mL_{r,s,\pell,u;\theta}$
representations into nonisomorphic ones,
$\bar\mL_{r,s,\pell,u}=\beta\,\mL_{r,s,\pell,u}$, which gives the
second half of the admissible $\hSSL21$-representations.  It follows
that the character of $\bar\mL_{r,s,\pell,u}$ can be written as
\begin{gather}\label{def-bar}
  \bar\chi_{r,s,\pell,u}(q,x,y)=\chi_{r,s,\pell,u}(q,x,y^{-1}),
\end{gather}
which after a simple calculation gives
\begin{gather}\label{chi-bar}
  \bar\chi_{r,s,\pell,u}=-\chi_{1-r,s,\pell,u;-s-1},
\end{gather}
and hence $\bar\chi_{r,s,\pell,u;\theta}=
-\chi_{1-r,s,\pell,u;\theta-s-1}$,
\begin{multline}\label{eq:gen-char-bar}
  \bar\chi_{r,s,\pell,u;\theta}(q,x,y)
  =-q^{(\theta-s)(-r-\frac{\pell}{u}(\theta-1))}\,
  x^{-\frac{r}{2} - \frac{s-1}{2}\frac{\pell}{u}}\,
  y^{-\frac{r}{2} + \frac{s+1}{2}\frac{\pell}{u} -
    \theta\frac{\pell}{u}}\\*
  {}\times
  \psi_{1-r,s,\pell,u}(q, x, y q^{2(\theta-s-1)})\,
%%   \frac{\vartheta_{1,0}(q,x^{\half}y^{\half})
%%     \vartheta_{1,0}(q,x^{\half}y^{-\half})}{
%%     \vartheta_{1,1}(q,x)\prod_{i\geq1}(1-q^i)^3}.
  \THETA(q,x,y).
\end{multline}

% In view of~\eqref{chi-bar}, both the $\mL_{r,s,\pell,u}$ characters
% and their $\beta$-images can be labeled by $\mL_{r,s,\pell,u}$ with
% \begin{gather*}
%   1-\pell\leq r\leq \pell,\qquad 1\leq s\leq u.
% \end{gather*}
%% A simple calculation also shows that under $\cT_\eta$, the characters
%% behave as
%% \begin{gather}
%%   \chi_{r,s,\pell,u
%% ;\theta|\eta}=(-1)^\eta\chi_{r,s-\eta,\pell,u
%% ;\theta+\eta}.
%% \end{gather}
%% We note that (in the exponential notation) the above characters have
%% the obvious property
%% \begin{gather*}
%%   \chi_{r,s,\pell,u;\theta}(\tau,\nu,\mu+2n)
%%   =\EX{-i\pi\frac{\pell}{u}(s+1+2\theta)n}
%%   \chi_{r,s,\pell,u;\theta}(\tau,\nu,\mu)
%%   ,\quad n\in\oZ.
%% \end{gather*}

% We note that %$\cU_{\half}\circ\cA_{-\half}\circ\alpha$ 
% $\gamma\circ\alpha$ maps the twisted highest-weight
% conditions~\eqref{twisted-hwcond} into~\eqref{second-hwcond}.

We finally consider the Verma modules $\verma{h_-, h_+;\theta}$ with
the same $h_-=\frac{r}{2}-\frac{\pell}{u}\frac{s-1}{2}$ as in the Verma
module involved in the construction of $\mL_{r,s,\pell,u;\theta}$.  Let
$\mM_{r,s,h_+;\theta}\equiv\mM_{r,s,h_+,\pell,u;\theta}$ be the quotient
with respect to the MFF singular vectors defined in~\cite{[ST]}.  In
the case where no charged singular vectors exist in the above
$\verma{h_-, h_+;\theta}$, i.e., for
\begin{gather}\label{cond-irreducible}
  \ffrac{u r}{2\pell}-\ffrac{s-1}{2} \pm \ffrac{u}{\pell}\,h_+
  \notin\oZ,
\end{gather}
the modules $\mM_{r,s,h_+;\theta}$ are irreducible and
$\mM_{r,s,h_+;\theta}\simeq\mM_{r,s,h_+ - \frac{\pell}{u}\theta;0}$.
Another straightforward calculation shows that the character of
$\mM_{r,s,h_+;\theta}$ is given~by
\begin{multline}\label{massive-char}
  \Omega_{r,s,h}(q,x,y)=
  y^{h}\,
  x^{-\frac{\pell}{u}\frac{s-1}{2}}\,  
  q^{\frac{\pell}{4u}(s-1)^2-\frac{u}{\pell}h^2}\\*
  {}\times\left(
    \theta_{r,\pell}(q^u,xq^{-(s-1)})
    -\theta_{-r,\pell}(q^u,xq^{-(s-1)})
  \right)
%%   \frac{\vartheta_{1,0}(q,x^{\half}y^{\half})\,
%%     \vartheta_{1,0}(q,x^{\half}y^{-\half})}{
%%     \vartheta_{1,1}(q,x)\prod_{i\geq1}(1-q^i)^3}.
  \THETA(q,x,y),
\end{multline}
where $h=h_+-(k+1)\theta$.  It follows that
\begin{gather}\label{omega+ell}
  \begin{gathered}
    \Omega_{r+2n\pell,s,h}(q,x,y)=\Omega_{r,s,h}(q,x,y),\qquad
    \Omega_{n\pell,s,h}(q,x,y)=0,\quad n\in\oZ,\\
    \Omega_{-r,s,h}(q,x,y)=-\Omega_{r,s,h}(q,x,y).
  \end{gathered}
\end{gather}

We use the special notation for the \textit{reducible} modules
$\mM_{r,s,h_+}$ with $h_+$ such that the left-hand side
of~\eqref{cond-irreducible} \textit{is} an integer, namely with
$h_+=\pm(\frac{r}{2}-\frac{\pell}{u}\frac{s+1}{2})$:\ 
$\mM_{r,s}=\mM_{r,s,\frac{r}{2}-\frac{\pell}{u}\frac{s+1}{2}}$ and
$\bar\mM_{r,s}=\mM_{r,s,-\frac{r}{2}+\frac{\pell}{u}\frac{s+1}{2}}$.
The respective characters of $\mM_{r,s}$ and $\bar\mM_{r,s}$ are given
by
\begin{gather}\label{Omega-short}
  \begin{split}
    \Omega_{r,s}(q,x,y)
    &=\Omega_{r,s,\frac{r}{2}-\frac{\pell}{u}\frac{s+1}{2}}(q,x,y),
    \\
    \bar\Omega_{r,s}(q,x,y)
    &=\Omega_{r,s,-\frac{r}{2}+\frac{\pell}{u}\frac{s+1}{2}}(q,x,y).
  \end{split}
\end{gather}
Reducibility of these modules can be expressed as the exact sequences
\begin{gather*}
  \begin{aligned}
    0&\to\mL_{r+1,s}\to
    \mM_{r,s}\to\mL_{r,s}\to0,
    \\
    0&\to\bar\mL_{r+1,s}\to
    \bar\mM_{r,s}\to\bar\mL_{r,s}\to0,
  \end{aligned}
  \qquad 1\leq r\leq\pell-1.
\end{gather*}
Consequently,
\begin{gather}\label{eq:reducible}
  \begin{split}
    \chi_{r+1,s}+\chi_{r,s}&=\Omega_{r,s},\\
    \bar\chi_{r+1,s}+\bar\chi_{r,s}&=\bar\Omega_{r,s},
  \end{split}
  \qquad 1\leq r\leq\pell-1,
\end{gather}
which can also be easily verified directly
using~\eqref{chi-through-K-exp} and~\eqref{diag-periodicity}.
%% Consequently,
%% \begin{gather*}
%%   \chi_{r,s,\pell,u;\theta}(q,x,y) +
%%   \chi_{r+1,s,\pell,u;\theta}(q,x,y)=
%%   \Omega_{r,s,\frac{r}{2}-\frac{\pell}{u}\frac{s+1+2\theta}{2}}(q,x,y)
%% \end{gather*}
We note that at the same time, $\Omega_{r,s}$ and $\bar\Omega_{r,s}$
are essentially (apart from their $y$-dependence) the admissible
characters of a level-$(k\!-\!1)$ \ $\hSL2$ algebra obtained from
$\hSSL21_k$ by the \textit{reduction} with respect to only the
fermionic generators; this reduction therefore sends reducible
$\hSSL21$-modules into irreducible $\hSL2$-ones, and the cohomology of
the complex associated with the reduction is certainly not
concentrated at one term.

\subsection{\chzerotext{} and \chhalftext{} characters and
  supercharacters}\label{sec:sectors} The $\hSSL21$-characters
introduced above are in the so-called \chzerotext{} sector.  We also
introduce \chhalftext{} characters and supercharacters in both sectors
as follows.  The \chhalftext{} characters are simply the
$\theta=-\half$ spectral flow transformations of the above
(\chzerotext) characters: for the (twisted) character $\chi$ of any
$\hSSL21$-module in the \chzerotext{} sector, the corresponding
(twisted) \chhalftext{} character is
$\chi^{\chhalf}_{;\theta}=\chi_{;\theta-\half}$, and therefore
(see~\eqref{sl21-sf})
%% \begin{multline}
%%   \chi^\chhalf_{r,s,\pell,u;\theta}(q,x,y)=
%%   y^{\frac{k}{2}}\,q^{-\frac{k}{4}}
%%   \chi_{r,s,\pell,u;\theta}(q,x,y\,q^{-1})={}\\
%%   {}=\EX{i\pi k\mu - i\pi k\frac{\tau}{2}}
%%   \chi_{r,s,\pell,u;\theta}(\tau,\nu,\mu-\tau).
%% \end{multline}
\begin{gather}
  \chi^{\chhalf}_{;\theta}(\tau,\nu,\mu)
  =\EX{i\pi k\mu\!-\!i\pi k\frac{\tau}{2}}
  \chi_{;\theta}(\tau,\nu,\mu-\tau).
\end{gather}

With the supercharacter of a module $\mP$ defined as
\begin{gather*}
  \sigma[\mP](\tau,\nu,\mu)=\Tr^{\phantom{y}}_{\mP}
  \bigl(
  \EX{2i\pi L_0\tau\!+\!2i\pi H^-_0\nu\!+\!2i\pi H^+_0(\mu\!+\!1)}
  \bigr),
\end{gather*}
we then have the \chzerotext{} supercharacters
$\chi^{\schzero}_{;\theta}\equiv\sigma^{\chzero}_{;\theta}$ and the
\chhalftext{} supercharacters
$\chi^{\schhalf}_{;\theta}\equiv\sigma^{\chhalf}_{;\theta}$ (for
$\theta\in\oZ$)
\begin{align}\label{sR2n}
  \chi^{\schhalf}_{;\theta}(\tau,\nu,\mu)
  &=
  \EX{i\pi k\mu\!-\!i\pi\frac{k}{2}\tau\!+\!i\pi k}
  \chi_{;\theta}(\tau,\nu,\mu\!-\!\tau\!+\!1),\\
  \chi^{\schzero}_{;\theta}(\tau,\nu,\mu)
  &=
  \chi_{;\theta}(\tau,\nu,\mu\!+\!1),
  \label{sNS2n}
\end{align}
expressed in terms of the \chzerotext{}
character~$\chi$.

Schematically, behavior of characters in the different sectors under
the $S$ and $T$ modular transformations can be summarized in the
diagram
\begin{gather}\label{ST-diagram}
  \xymatrix@=40pt{%
    \chzero\ar@/^/[0,1]|{S}
    \ar@{->}@(lu,ru)[]|{T}
    &\schhalf\ar@/^/[0,-1]|{S}
    \ar@/^/[1,0]|{T}
    \\
    \schzero\ar@{->}@(lu,ru)[]|{T}
    \ar@{->}@(ld,rd)[]|{S}
    &\chhalf\ar@/^/[-1,0]|{T}
    \ar@{->}@(ld,rd)[]|{S}
  }
\end{gather}

\def\NPB{Nucl.\ Phys.\ B}
\def\PLB{Phys.\ Lett.\ B}
\def\MPLA{Mod.\ Phys.\ Lett.\ A}
\def\CMP{Commun.\ Math.\ Phys.}
\def\IJMPA{Int.\ J.\ Mod.\ Phys.\ A}

\end{document}